\documentclass[a4paper, 10pt]{amsart}
\usepackage[utf8]{inputenc}
\usepackage[T1]{fontenc}
\usepackage{amsmath}
\usepackage{amssymb}
\usepackage{mathtools}
\usepackage{array}

\usepackage[svgnames]{xcolor}
\definecolor[named]{LightGray}{rgb}{0.95,0.95,0.96}
\definecolor{codegreen}{rgb}{0,0.6,0}
\definecolor{codegray}{rgb}{0.5,0.5,0.5}
\definecolor{sageblue}{RGB}{94, 94, 255}
\definecolor{codepurple}{rgb}{0.58,0,0.82}
\definecolor{backcolour}{rgb}{0.98,0.98,0.96}

\usepackage{listings}
\lstdefinestyle{mystyle}{
    frame=tb,
    language={Python},
    keywordstyle=\color{magenta},
    escapeinside={@}{@},
    backgroundcolor=\color{backcolour},   
    commentstyle=\color{codegreen},
    numberstyle=\tiny\color{codegray},
    stringstyle=\color{codegreen},
    basicstyle=\ttfamily\footnotesize,
    breakatwhitespace=false,         
    breaklines=true,                 
    captionpos=b,                    
    keepspaces=true,                 
    numbersep=5pt,                  
    showspaces=false,                
    showstringspaces=false,
    showtabs=false,                  
    tabsize=2,
    upquote=true,
}
\lstdefinestyle{pythonconsole}{
    style=mystyle,
    moredelim=[il][\color{codegray}]{PYOUT:},
	literate={·}{{$\cdot$}}1 {↦}{{$\mapsto$}}1,
}
\usepackage{tikz}
\usetikzlibrary{arrows.meta}

\theoremstyle{plain}
\newtheorem{theorem}{\bf Theorem}[section]
\newtheorem{proposition}[theorem]{\bf Proposition}
\newtheorem{lemma}[theorem]{\bf Lemma}

\theoremstyle{definition}

\theoremstyle{definition}
\newtheorem{claim}{}[theorem]

\newcommand{\N}{\mathbb N}
\newcommand{\Z}{\mathbb Z}

\newcommand{\Q}{\mathbb Q}
\newcommand{\F}{\mathbb F}

 \DeclareMathOperator{\ord}{ord}
 
 \DeclareMathOperator{\im}{im}
 \DeclareMathOperator{\supp}{supp}

 \DeclareMathOperator{\AG}{AG}
\DeclareMathOperator{\Aut}{Aut} \DeclareMathOperator{\GL}{GL}
\DeclareMathOperator{\Stab}{Stab} \DeclareMathOperator{\lex}{lex}
\DeclareMathOperator{\orb}{Orb}

\newcommand{\und}{\;\mbox{ and }\;}

\newcommand{\red}{{\text{\rm red}}}

\newcommand{\la}{\langle}
\newcommand{\ra}{\rangle}

\newcommand{\code}[1]{\texttt{\detokenize{#1}}}

\numberwithin{equation}{section}

\usepackage[colorlinks]{hyperref}
\definecolor{linkblue}{RGB}{1,1,190}
\definecolor{citered}{RGB}{190,1,1}
\hypersetup{
	linkcolor=linkblue,
	urlcolor=linkblue,
	citecolor=citered
}

\renewcommand{\t}{\, | \,}

\begin{document}

\title[The system of  length sets of transfer Krull monoids over $C_3 \oplus C_3 \oplus C_3$]{The system of  length sets of \\ transfer Krull monoids over $C_3 \oplus C_3 \oplus C_3$}

\author{Alfred Geroldinger and  Benjamin Hackl and Qinghai Zhong}

\address{University of Graz, NAWI Graz\\
Department of Mathematics and Scientific Computing\\
Heinrichstra{\ss}e 36\\
8010 Graz, Austria (A.~Geroldinger, B.~Hackl)}
\urladdr{https://geroldinger.github.io/, https://benjamin-hackl.at/}

\address{Xidian University \\
Institute of Mathematics and Interdisciplinary Sciences \\
Xi'an, China (Q. Zhong)}
\urladdr{https://addcombin.github.io/zhong/}

\subjclass[2020]{20M13, 13A05, 11B30, 11R27}
\keywords{Transfer Krull monoids, length sets, zero-sum sequences, sets of distances}

\begin{abstract}
Let $H$ be a Krull monoid with class group $G$ and suppose that each class contains a prime divisor. There are precisely nine groups (up to isomorphism) for which the set of distances of $H$ is contained in $\{1, 2\}$. For all these groups, except for $C_3 \oplus C_3 \oplus C_3$, the system of length sets has been written down explicitly. In this paper we close this gap and determine the system of length sets when the class group is isomorphic to $C_3 \oplus C_3 \oplus C_3$.
\end{abstract}

\maketitle

\smallskip
\section{Introduction} \label{1}
\smallskip

Let $H$ be a Krull monoid (e.g., the monoid of nonzero elements of a Krull domain) with class group $G$ and suppose that each class contains a prime divisor (classic examples in ring and number theory have this property; they, together with the more general class of transfer Krull monoids, will be discussed in Section \ref{2}). 
Every non-unit $a \in H$ can be written as a finite product of irreducible elements (atoms), and factorization is unique for all elements if and only if the class group is trivial.  If  $a = u_1 \cdot \ldots \cdot u_k$ is a factorization into atoms $u_1, \ldots, u_k$, then $k$ is a factorization length of $a$ and $\mathsf L (a) = \{k \colon k \ \text{is a factorization length of $a$} \} \subset \N$ is the {\it length set of $a$}. If $a \in H$ is invertible, then we set $\mathsf L (a) = \{0\}$. All length sets are finite and the {\it system of length sets}
\[
\mathcal L (H) = \{ \mathsf L (a) \colon a \in H\}
\]
is  a main  invariant describing the arithmetic of $H$. 

The system $\mathcal L (H)$ just depends on the class group $G$ (see Subsection \ref{2.4}). If $G$ is infinite, then every finite nonempty subset of $\N_{\ge 2}$ occurs as a length set. Suppose that $G$ is finite. 
The Structure Theorem for Length Sets says that  length sets  are almost arithmetic multiprogressions with global bounds for all parameters (\cite[Section 7.3]{Ge-Gr-Zh26a}), and this description is best possible (\cite{Sc09a}). 
It is classic that all length sets are singletons if and only if the class group $G$ has at most two elements. Furthermore, the class groups are characterized for which all length sets are arithmetic progressions, or almost arithmetic progressions (\cite[Corollaries 8.2.21 and 9.2.7]{Ge-Gr-Zh26a}). These special cases occur for a very small number of small groups. 
However, in general there is   only limited information (in terms of the class group $G$) on the parameters controlling the structure of length sets. 

A key parameter is the set of distances of $H$. To recall its definition, let
$L = \{m_0, \ldots, m_k\} \subset \Z$ be a finite nonempty subset,  with $k \in \N_0$ and with $m_0 < \ldots < m_k$. Then  $\Delta (L) = \{m_i - m_{i-1} \colon i \in [1, k] \} \subset \N$ is the set of distances of $L$, and 
\[
\Delta (H) = \ \bigcup_{L \in \mathcal L (H)} \Delta (L) \ \subset \ \N
\]
is the {\it set of distances} of $H$.
If $\Delta (H) \ne \emptyset$, then the set of distances $\Delta (H)$ is an interval with $\min \Delta (H) = 1$ (\cite[Theorem 6.1.5]{Ge-Gr-Zh26a}), and $\Delta (H)$ is finite if and only if $G$ is finite. More precisely, we have
\begin{itemize}
\item[(a)] $\Delta (H) = \N$ if and only if $G$ is infinite (\cite[Theorem 6.4.2]{Ge-HK06a}).

\item[(b)] $\Delta (H) = \emptyset$ if and only if $|G| \le 2$ (\cite[Corollary 3.4.12]{Ge-HK06a}).

\item[(c)] $\Delta (H) = \{1\}$ if and only if $G$ is isomorphic to one of the following groups (\cite[Corollary 6.2.4]{Ge-Gr-Zh26a})
    \[
    C_3, \ C_2 \oplus C_2, \ C_3 \oplus C_3 \,.
    \]

\item[(d)] $\Delta (H) = \{1,2\}$ if and only if $G$ is isomorphic to one of the following groups (\cite[Corollary 6.2.14]{Ge-Gr-Zh26a})
    \[
    C_4, \ C_2 \oplus C_4, \ C_2 \oplus C_2 \oplus C_2, \ C_3 \oplus C_3 \oplus C_3 \,.
    \]
\end{itemize}

There is no such characterization for any further set of distances $\Delta (H)$, apart from the listed ones. For all groups -- apart from $C_3 \oplus C_3 \oplus C_3$ --  listed in (c) and (d), the systems of length sets have been written down (see \cite[Theorem 7.3.2]{Ge-HK06a}, \cite[Proposition 3.12]{Ge-Sc16b}, and \cite[Theorem 4.5]{Ge-Sc-Zh17b}).

The goal of this paper is to close this gap and to write down the system of length sets when the class group is isomorphic to $C_3 \oplus C_3 \oplus C_3$. We formulate our main result.

\begin{theorem} \label{1.1}
Let $H$ be a transfer Krull monoid over the group $C_3 \oplus C_3 \oplus C_3$. 
Then, for the system of length sets $\mathcal L (H)$, we have
$\mathcal L (H) = \mathcal L_1\sqcup \mathcal L_2\sqcup \mathcal L_3\sqcup \mathcal L_4\sqcup \mathcal L_5$, where
\begin{align*}
\mathcal L_1&= \big\{\{y\}\colon y\in \N_0  \big\}\,,\\
\mathcal L_2&=\big\{ [2k,\nu] \colon k \in \N, \nu\in [2k+1, 7k]  \big\} \ \sqcup \ \big\{ [2k+1, \nu] \colon  k \in \N, \nu \in [2k+2, 7k+3] \big\}\,, \\
\mathcal L_3
		&=\big\{ 3k+2\cdot[0,\nu] \colon k \in \N, \nu\in [1, 2k] \big\} \ \sqcup \ \big\{3k+1+2\cdot[0, \nu] \colon  k \in \N, \nu \in [1, 2k] \big\}  \\
		&\phantom{ =\big\{3k+2\cdot[0,\nu] \colon k \in \N, \nu\in [1, 2k] \big\} \cup }  \sqcup \ \big\{3k+2+2\cdot[0, \nu] \colon  k \in \N_0, \nu \in [1, 2k+1] \big\}\,,  \\
\mathcal L_4 & =  \big\{ [4,10]\cup\{12\},  [6, 14]\cup\{16\} \big\} \ \sqcup	\\
	& \quad\quad \big\{[3k+2, \nu]\cup\{\nu+2\} \colon k \in \N_0, \nu\in [3k+3, 7k+5] \big\} \ \sqcup \\
	  & \quad\quad \big\{[3k+3, \nu]\cup\{\nu+2\} \colon k \in \N_0, \nu\in [3k+4, 7k+6] \big\} \ \sqcup\\
	 & \quad\quad  \big\{[3k+4, \nu]\cup\{\nu+2\} \colon k \in \N_0, \nu\in [3k+5, 7k+9] \big\}\,, \quad \text{and}  \\           
\mathcal L_5&=\big\{\{3k+2\}\cup[3k+4,\nu] \colon k \in \N_0, \nu\in [3k+5, 7k+5] \big\} \ \sqcup \\
&\quad\quad  \big\{\{3k+3\}\cup[3k+5,\nu] \colon k \in \N_0, \nu\in [3k+6, 7k+7] \big\} \ \sqcup\\
&\quad\quad \big\{\{3k+4\}\cup[3k+6,\nu] \colon k \in \N_0, \nu\in [3k+7, 7k+9] \big\}\,.
	\end{align*}
\end{theorem}

\smallskip
Let $H$ be as in Theorem \ref{1.1} and let $G = C_3^3$. Factorization Theory provides some structural information on $\mathcal L (H)$. 
We have $\mathcal L (H) = \mathcal L (G)$ (see Subsection \ref{2.4}), $\Delta (H) = \{1,2\}$, the elasticity $\rho (H) = 7/2$, and length sets having maximum elasticity are arithmetic progressions with difference $1$ (see Subsection \ref{sec:basic}). Moreover, 
length sets are arithmetic progressions with difference $1$ or $2$ with possible gaps at the initial and end parts. There is a global bound for the initial and end parts of length sets where gaps can arise, but there  is no information on the size of this global bound. However, these results alone are enough to show that $\mathcal L (G)$ is not equal to $\mathcal L (G')$ for any finite abelian group $G'$ that is not isomorphic to $G$.

On the other hand, there is no precise information on the structure of length sets $L$ with $\min L$ small, in particular not for $\min L = 2$. We handle the length sets, whose minimum is small,  with computational methods. To do so, we establish a Python software package which allows to study length sets of a product of a small number of (partly very special) atoms. In this computational part we also worked with AI assistance. All this is described in  Section \ref{software}, and the associated software packages and databases can be found at \cite{Ge-Ha-Zh27b}. These computational results are the starting point for  precise structural results on general length sets (Section \ref{5}). The proof of Theorem \ref{1.1} will then be given in Section \ref{6}.

\section{Background on Transfer Krull Monoids and Factorizations} \label{2}

Let $\N$ be the set of positive integers and let $\N_0 = \N \cup \{0\}$. For real numbers $a$ and $b$, we denote by $[a, b] = \{ x \in \Z \colon a \le x \le b\}$ the discrete interval between $a$ and $b$.
For subsets $A,B\subset\Z$,  $A+B=\{a+b\colon a\in A,b\in B\}$ is the {\it sumset} of $A$ and $B$, and for $k\in\Z$, we set $k+A=\{k\}+A$. For $k\in\N$, let $A_{\geq k}=\{a\in A \colon a\geq k\}$ and let $k\cdot A=\{ka\colon a\in A\}$ be the {\it dilation} of $A$ by $k$. If $A=\{m_0,\ldots, m_k\}\subset\Z$, with $k\in\N_0$ and $m_0<\ldots<m_k$, then $\Delta(A)=\{m_i-m_{i-1}\colon i\in [1,k]\}$ is the {\it set of distances} of $A$. 
If $A \subset \N$, then $\rho (A) = m_k/m_0 \in \Q_{\ge 1}$ is the {\it elasticity} of $A$.

 For a set $P$, let $\mathcal F (P)$ denote the free abelian monoid with basis $P$. An element $a \in \mathcal F (P)$ will be written in the form
\[
a = \prod_{p \in P} p^{\mathsf v_p (a)} \,,
\]
where $\mathsf v_p \colon \mathcal F (P) \to \N_0$ is the $p$-adic valuation, and we call
\begin{itemize}
\item $|a|_{\mathcal F (P)} = |a|= \sum_{p \in P} \mathsf v_p (a) \in  \N_0$ the {\it length} of $a$, and

\item $\supp_{\mathcal F (P)}(a) = \supp (a) = \{p \in P \colon \mathsf v_p (a) > 0 \} \subset P$ the {\it support} of $a$.
\end{itemize}
Let $S$ be a multiplicatively written commutative semigroup with identity element. We denote by $S^{\times}$ its group of invertible elements and we say that $S$ is reduced if $S^{\times} = \{1\}$. Then  $S_{\red} = \{ aS^{\times} \colon a \in S\}$ is  the {\it associated reduced semigroup} of $S$.  The semigroup $S$ is
\begin{itemize}
\item {\it cancellative} if $a,b,u \in S$ and $au=bu$ implies that $a=b$.

\item {\it unit-cancellative} if $a, u \in S$ and $a=au$ implies that $u \in S^{\times}$.
\end{itemize}

\smallskip
\noindent
\subsection{Arithmetic of Monoids.} \label{2.1}
A monoid means a commutative unit-cancellative semigroup with identity element.
Let $H$ be a monoid. An element $u \in H$ is said to be  {\it irreducible} (or an {\it atom}) if $u \notin H^{\times}$ and, for all $a, b \in H$, $u=ab$ implies that $a \in H^{\times}$ or $b \in H^{\times}$. We denote by $\mathcal A (H)$ the set of atoms of $H$.
We say that $H$ is {\it atomic} if every  $a \in H \setminus H^{\times}$ can be written as a finite product of atoms.

The free abelian monoid $\mathsf Z (H) = \mathcal F \big(\mathcal A (H_{\red}) \big)$ is the {\it factorization monoid} of $H$ and the homomorphism
\begin{equation} \label{factorization-hom}
\pi_H = \pi \colon \mathsf Z (H) \to H_{\red}, \quad \text{defined by} \quad \pi (u) = u \ \text{for all $u \in \mathcal A (H_{\red})$}
\end{equation}
is the {\it factorization homomorphism} of $H$. For $z \in \mathsf Z (H)$, $|z|_{\mathsf Z (H)} = |z| \in \N_0$ is the length of the factorization. For each $a \in H$,
\begin{itemize}
\item $\mathsf Z_H (a) = \mathsf Z (a) = \pi^{-1} (aH^{\times}) \subset \mathsf Z (H)$ is the {\it set of factorizations} of $a$, and

\item $\mathsf L_H (a) = \mathsf L (a) = \{|z| \colon z \in \mathsf Z (a) \} \subset \N_0$ is the {\it length set} of $a$.
\end{itemize}
Thus, the elements of $\mathsf Z (H)$ are formal products of atoms from $H_{\red}$.   Note that $H$ is atomic if and only if $\mathsf Z (a) \ne \emptyset$ for all $a \in H$. For $k \in \N_0$, we set
\[
\mathsf Z_k (H) = \{ z \in \mathsf Z (H) \colon |z|=k\} \quad \text{and} \quad \mathsf Z_k (a) = \mathsf Z (a) \cap \mathsf Z_k (H) \,.
\]
Let $z, z' \in \mathsf Z (H)$. Then we can write
\[
z = u_1\cdot\ldots\cdot u_{\ell} v_1\cdot\ldots\cdot v_m\quad
\text{and}\quad z' = u_1\cdot\ldots\cdot u_{\ell} w_1\cdot\ldots
\cdot w_n\,,
\]
where $\ell,\,m,\,n\in\N_0$ and $u_1,\ldots,u_\ell,\,v_1,
\ldots,v_m,\,w_1,\ldots,w_n\in\mathcal A(H_\red)$ are such that
\[
\{v_1,\ldots,v_m\}\cap\{w_1,\ldots,w_n\} =\emptyset\,.
\]
Then $\gcd(z,z')=u_1\cdot\ldots\cdot u_{\ell}$, and we call
\[
\mathsf d (z,z') =\max\{m,\,n\} =\max\{ |z\gcd (z,z')^{-1}|,
|z'\gcd (z,z')^{-1}|\}\in\N_0
\]
the {\it distance} between $z$ and $z'$.

Let $a\in H$ and $N\in\N_0$. A finite sequence $z_0,\ldots,
z_k\in\mathsf Z (a)$ is called an {\it $N$-chain} of factorizations of $a$ if
$\mathsf d (z_{i-1},z_i)\le N$ for all $i\in [1,k]$. We denote by
$\mathsf c (a)$ the smallest $N\in\N _0\cup\{\infty\}$ such
that any two factorizations $z,\,z'\in\mathsf Z (a)$ can be
concatenated by an $N$-chain. Then
\[
\mathsf c(H) =\sup\{\mathsf c(b)\colon b\in H\}\in\N_0\cup
\{\infty\}
\]
denotes the {\it catenary degree} of $H$. If $H$ is atomic and all length sets are finite (which will be the case for all monoids under consideration), then 
\[
\Delta (H) = \bigcup_{L \in \mathcal L (H)} \Delta (L) \ \subset \N
\]
is the {\it set of distances} of $H$.

\smallskip
\noindent
\subsection{Krull Monoids.} \label{2.2}
We gather the basic properties of Krull monoids (for details see \cite[Chapter 2]{Ge-HK06a}). A monoid is said to be a {\it Krull monoid} if it is cancellative and satisfies one of the following equivalent conditions.
\begin{itemize}
\item[(a)] $H$ is completely integrally closed and satisfies the ACC on divisorial ideals.

\item[(b)] $H$ has a divisor theory.

\item[(c)] There is a divisor homomorphism $\varphi \colon H \to D$, where $D$ is a free abelian monoid.
\end{itemize}

A commutative domain $D$ is Krull if and only if its multiplicative monoid $D\setminus \{0\}$ of nonzero elements is Krull  if and only if its monoid $\mathcal I^* (D)$ of invertible ideals is a Krull monoid. Holomorphy domains in global fields are Dedekind domains (and hence Krull domains) with finite class group and infinitely many prime ideals in each class. Moreover, for every finite abelian group $G$, there are an algebraic number field $K$ and a finite set $S$ of places such that $G$ is isomorphic to the class group of the ring of $S$-integers of $K$ (\cite{Pe99}).
Thousands  of algebraic number fields with class group being isomorphic to $C_3 \oplus C_3 \oplus C_3$ can be found on the website \cite{lmfdb} (note that, according to Conjecture 1.3 of \cite{H-J-K-M-P19}, there is no imaginary quadratic number field with class group isomorphic to $C_3 \oplus C_3 \oplus C_3$).

To give a further example of Krull domains having prime divisors in all classes, we consider monoid algebras (for more examples, see \cite{Ge-HK06a, Ge-Zh20a}). Let $D$ be a Krull domain and let $H$ be a Krull monoid such that $H^{\times}$ satisfies the ACC on cyclic subgroups. Then $D[H]$ is a Krull domain with class group being isomorphic to $\mathcal C (D) \oplus \mathcal C (H)$ and every class contains infinitely many prime divisors (\cite{Fa-Wi22b}).

\smallskip
\noindent
\subsection{Monoids of Zero-Sum Sequences.} \label{2.3}
Let $G$ be an additive abelian group and let $G_0 \subset G$ be a nonempty subset. Then $\langle G_0 \rangle$ denotes the subgroup generated by $G_0$, and we set $G_0^{\bullet}= G_0 \setminus \{0\}$. For every $n \in \N$, $C_n$ denotes an additive cyclic group of order $n$.
In additive combinatorics, a {\it sequence} (over $G_0$) means a finite sequence of terms from $G_0$, where repetition is allowed and the order of the elements is disregarded. As usual, we consider sequences as elements of the free abelian monoid with basis $G_0$. Let
\[
S = g_1 \cdot \ldots \cdot g_{\ell} = \prod_{g \in G_0} g^{\mathsf v_g (S)} \in \mathcal F (G_0)
\]
be a sequence over $G_0$. We  call
\begin{itemize}
\item $ \supp (S)  = \{g \in G \colon \mathsf v_g (S) > 0 \} \subset G \ \text{the \ {\it support} \ of \ $S$}\,$,
\item $|S|_{\mathcal F (G)} =  |S|  = \ell = \sum_{g \in G} \mathsf v_g (S) \in \mathbb N_0 \ \text{the \ {\it length} \ of \ $S$} \,$,
\item $\mathsf h (S) = \max \{ \mathsf v_g (S) \colon g \in G_0\} \in \N_0 \ \text{the \ {\it maximum of the multiplicities} \ of \ $S$}\,$,
\item $ \sigma (S)  = \sum_{i = 1}^{\ell} g_i \ \text{the \ {\it sum} \ of \
$S$} \,$,
\item $\Sigma (S)  = \Big\{ \sum_{i \in I} g_i \colon \emptyset \ne I \subset [1,\ell] \Big\} \ \text{ the \ {\it set of
subsequence sums} \ of \ $S$} \,$, \ and
\item $\Sigma_n (S) = \Big\{ \sum_{i \in I} g_i \colon \ I \subset [1,\ell] \ \text{with} \ |I|=n \Big\}$ \ the {\it set of $n$-sums}  of $S$ for every $n \in \N$.
\end{itemize}
The sequence $S$ is said to be
\begin{itemize}
\item a {\it zero-sum sequence} \ if \ $\sigma (S) = 0$,

\item a {\it minimal zero-sum sequence} \ if it is a nontrivial zero-sum
      sequence and has no proper zero-sum subsequences.
\end{itemize}
The set of zero-sum sequences $\mathcal B (G_0) = \{S \in \mathcal F (G_0) \colon \sigma (S)=0\} \subset \mathcal F (G_0)$ is a submonoid, and the set of minimal zero-sum sequences is the set of atoms of $\mathcal B (G_0)$. Since the inclusion $\mathcal B (G_0) \hookrightarrow \mathcal F (G_0)$ is a divisor homomorphism, $\mathcal B (G_0)$ is a Krull monoid by Condition (c) in Subsection \ref{2.2}. If $|G| \le 2$, then $\mathcal B (G)$ is factorial and if $|G| \ge 3$, then the inclusion $\mathcal B (G) \hookrightarrow \mathcal F (G)$ is a divisor theory with class group isomorphic to $G$ and every class contains precisely one prime divisor.

Let $G_0 \subset G$ be a subset and let $S \in \mathcal F (G)$.
If $G'$ is an abelian group and $\varphi \colon G \to G'$ is a map, then $\varphi (S) = \varphi (g_1) \cdot \ldots \cdot \varphi (g_{\ell}) \in \mathcal F (G')$ and we set $-S = (-g_1) \cdot \ldots \cdot (-g_{\ell}) \in \mathcal F (G)$.  Moreover,
\begin{equation} \label{sequence-sub}
S_{G_0} = \prod_{g \in G_0} g^{\mathsf v_g (S)} \in \mathcal F (G_0)
\end{equation}
denotes the subsequence of $S$ consisting of all terms lying in $G_0$.
We denote by $\Aut (G)$ the automorphism group of $G$. Then
\[
\Stab (G_0) = \{ \varphi \in \Aut (G) \colon \varphi (G_0)=G_0 \} \subset \Aut (G)
\]
is the {\it stabilizer} of $G_0$ and
\[
\Stab (S) = \{ \varphi \in \Aut (G) \colon \varphi (S)=S \} \subset \Aut (G)
\]
is the {\it stabilizer} of $S$. Clearly, $\Aut (G)$ acts on $\mathcal B (G)$, on $\mathcal A (G)$, and on $\mathsf Z (G)$. For $S \in \mathcal B (G)$,
\[
\orb (S) = \{\varphi (S) \colon \varphi \in \Aut (G) \}
\]
denotes the {\it orbit} of $S$, and 
\begin{equation} \label{orbit-stability}
\mathsf L (S) = \mathsf L (S') \quad \text{ for every $S' \in \orb (S)$} \,.
\end{equation}
Furthermore, we have
\begin{equation} \label{shift-by-zeros}
\mathsf L (0^kS) = k + \mathsf L (S) \quad \text{ for every $S \in \mathcal B (G^{\bullet})$} \,.
\end{equation}
We use the above two properties without further mention.

For any arithmetical invariant $*(H)$ defined for a monoid $H$, we write $*(G_0)$ instead of $* \big(\mathcal B (G_0)\big)$. In particular, $\mathcal A (G_0) = \mathcal A \big(\mathcal B (G_0) \big)$ is the set of atoms of $\mathcal B (G_0)$, $\mathsf Z (G_0) = \mathsf Z \big( \mathcal B (G_0) \big)$ is the factorization monoid of $\mathcal B (G_0)$, $\mathcal L (G_0)=\mathcal L \big( \mathcal B(G_0) \big)$ is the system of  length sets of $\mathcal B (G_0)$, and so on. We denote by
\[
\mathsf D (G_0) = \sup \{ |S| \colon S \in \mathcal A (G_0) \} \in \N_0 \cup \{\infty\}
\]
the {\it Davenport constant} of $G_0$.  Suppose that $G$ is finite, say $G \cong C_{n_1} \oplus \ldots \oplus C_{n_r}$, with $r \in \N_0$, $1 < n_1 \t \ldots \t n_r$. Then $r = \mathsf r (G)$ is the rank of $G$, $n_r = \exp (G)$ is the exponent of $G$, and $\mathsf D (G)$ is the smallest integer such that, for any $S \in \mathcal F (G)$,
\[
|S| \ge \mathsf D (G) \quad \text{implies} \quad 0 \in \Sigma (S) \,.
\]
Similarly, the
\begin{itemize}
\item {\it Erd{\H{o}}s-{G}inzburg-{Z}iv constant} $\mathsf s (G)$ is the smallest integer such that, for any $S \in \mathcal F (G)$,
      \[
      |S| \ge \mathsf s (G) \quad \text{implies} \quad 0 \in \Sigma_{\exp (G)} (S) \,.
      \]

\item {\it Harborth constant} $\mathsf g (G)$ is the smallest integer such that, for any squarefree $S \in \mathcal F (G)$,
      \[
      |S| \ge \mathsf g (G) \quad \text{implies} \quad 0 \in \Sigma_{\exp (G)} (S) \,.
      \]
\end{itemize}
To recall some facts on these invariants, note that
\begin{equation} \label{inequality-Davenport}
1 + \sum_{i=1}^r (n_i-1) \le \mathsf D (G) \le |G| \,,
\end{equation}
with equality on the left side if $G$ is a $p$-group.  If $r \le 2$ (with $n_1=1$ if $r=1$), then
\begin{equation} \label{davenport-rank2}
\mathsf D (C_{n_1} \oplus C_{n_2}) = n_1+n_2-1 \quad \text{and} \quad \mathsf s (C_{n_1} \oplus C_{n_2}) = 2n_1+2n_2-3 \,
\end{equation}
 (\cite[Theorem 5.8.3]{Ge-HK06a}). Moreover, we have (for more on the Erd{\H{o}}s-{G}inzburg-{Z}iv constant, see \cite{G-H-S-T07,E-E-G-K-R07,Fa-Ga-Zh11a})
\begin{equation} \label{EGZ}
 \mathsf s (C_3^3) = 19 \quad \text{and} \quad \mathsf s (C_3^4) = 41 \,.
\end{equation}

\smallskip
\noindent
\subsection{Transfer Homomorphisms and Transfer Krull monoids.} \label{2.4}
A monoid homomorphism  $\theta \colon H \to B$ is called a {\it transfer homomorphism} if it has the following properties.

\begin{enumerate}
\item[{\bf (T\,1)\,}] $B = \theta(H) B^\times$  and  $\theta^{-1} (B^\times) = H^\times$.

\item[{\bf (T\,2)\,}] If $u \in H$, \ $b,\,c \in B$  and  $\theta (u) = bc$, then there exist \ $v,\,w \in H$ such that  $u = vw$, \  $\theta (v) \in bB^{\times}$, and  $\theta (w) \in c B^{\times}$.
\end{enumerate}
It is easy to check that, for all $a \in H$, we have
\[
\mathsf L_H (a) = \mathsf L_B \big( \theta (a) \big) \quad \text{and hence} \quad \mathcal L (H) = \mathcal L (B) \,.
\]
An atomic monoid $H$ is said to be {\it transfer Krull} (over $G_0$) if one of the following equivalent conditions holds.
\begin{enumerate}
\item[(a)] There exist a Krull monoid $B$ and a transfer homomorphism $\theta \colon H \to B$.

\item[(b)] There exist an abelian group $G$, a subset $G_0 \subset G$, and a transfer homomorphism $\theta \colon H \to \mathcal B (G_0)$.
\end{enumerate}

\smallskip
If $H$ is a Krull monoid with class group $G$ and prime divisors in all classes, then there is a transfer homomorphism $\theta \colon H \to \mathcal B (G)$. Furthermore, there are  transfer Krull monoids and domains that are not Krull. Transfer Krull monoids need not be cancellative and (in a more general setting which we do not introduce here) they also include classes of non-commutative rings. We  give some  examples and refer to extended lists of transfer Krull monoids that can be found in \cite[Chapter 1]{Ge-Gr-Zh26a}, \cite{Ge-Zh20a} and  \cite{Ba-Re22a}.
There is a characterization of orders in Dedekind domains with torsion class groups that are transfer Krull (\cite{Ra26a}; these orders are commutative noetherian domains that are not integrally closed whence they are not Krull). There are monoids of modules, which are unit-cancellative but not necessarily cancellative, and transfer Krull (\cite{Ba-Sm22a}).
Next, let $\mathcal O$ be a holomorphy domain in a global field $K$, $A$ a central simple algebra over $K$, and let $R$ be a classical maximal $\mathcal O$-order of $A$ (classical maximal orders are non-commutative Dedekind domains). If every stably free left $R$-ideal is free, then $R$ is transfer Krull over a ray class group of $\mathcal O$ (\cite[Theorem 1.1]{Sm13a}).  Ray class groups are finite abelian groups (\cite{HK20a, HK22a}), and computer algebra systems such as PARI/GP~\cite{pari:2.19.1} and Hecke~\cite{Fi-Ha-Ho-Jo17} can be used to determine them explicitly. Whenever the ray class group in question is isomorphic to $C_3 \oplus C_3 \oplus C_3$, Theorem~\ref{1.1} describes the system of length sets of $R$.

\smallskip
\noindent
\subsection{Affine Spaces and Caps} \label{2.5}
We need some concepts from finite geometry. Let $r, s \in \N$, let $p$ be a prime, and let $q$ be a  power of $p$. An elementary $p$-group of rank $r$ can be viewed as an $r$-dimensional vector space over the finite field $\F_p$.
We identify the affine space $\AG (r,q)$ with $\F_q^r$.
A subset $C \subset \F_q^r$ is called an {\it $m$-cap} if $|C|=m$ and any three points of $C$ are not collinear.  

We apply this machinery to $p=q=3$ and recall a lemma that reveals how zero-sum concepts are connected with the geometric concepts (for a proof and more details see \cite[Lemma 5.2]{E-E-G-K-R07}). In particular, we note that a cap in $\F_3^r$ is a subset where no three elements sum to zero. 

\begin{lemma} \label{2.6}
Let $G = \F_3^r$ with $r \in \N$.
\begin{enumerate}
\item The map $f$, defined by $f (T) = \supp (T)$, is a bijection
       from the set
       \[
       \{ T \in \mathcal F ( G ) \colon T \quad \text{is
       squarefree and has no zero-sum subsequence of length $3$}
       \}
       \]
       onto the set $\{ C \subset G \colon C \quad
       \text{is a cap} \}$.

\smallskip
\item $\mathsf g ( G ) - 1$ is the maximal size of a
       cap in $\AG (r,3)$.
\smallskip
\item $\mathsf s ( G ) = 2 \mathsf g ( G )-1$, and every
       sequence $S \in \mathcal  F ( G )$ of length
       $|S| = \mathsf s ( G ) - 1$ which has
       no zero-sum subsequence
       of length $3$ has
       the form $S = T^2$ where $\supp (T)$ is a maximal cap in
       $G$. Conversely, if $T \in \mathcal F ( G )$ is squarefree
       and has no zero-sum subsequence of
       length $3$, then $T^2$ has no zero-sum subsequence of length
       $3$.
\end{enumerate}
\end{lemma}

The {\it Cap Set Problem} asks for the largest possible size of a cap in $\F_3^r$.
Part 2 implies that the maximal size of a cap in $\F_3^r$ equals $(\mathsf s (G)-1)/2$. Thus, for $r=2$, $r=3$, resp. $r=4$ we get $4$ by \eqref{davenport-rank2}, and  $9$ resp. $20$ by \eqref{EGZ}. The precise values for $\mathsf s (G)$ are known for $r \le 6$ (see \cite{E-F-L-S02,Po08a}), and 
Ellenberg and Gijswijt \cite{El-Gi17a} proved (by the polynomial method) that there is some $c < 3$ such that the largest size of a cap is bounded by $c^r$ for  all $r \in \N$.

For $r \in [3,5]$ the maximal cap in $\F_3^r$ is unique up to affine transformations, that is, up to automorphisms and translations (\cite[Lemma 5.4]{E-E-G-K-R07}). The automorphism group $\Aut (G)$ of $G = \F_3^3$ is the general linear group $\GL (3, \F_3)$, and we have $|\GL (3, \F_3)|= 11232$. If $(e_1, e_2,e_3)$ is a basis of $G =\F_3^3$, then
\[
G_0 = \{0, e_1,e_2,e_3,  e_1+e_2, e_1+e_3, -e_1+e_2+e_3, -e_1-e_2+e_3, -e_1+e_2-e_3\}
\]
is a maximal cap. A cap $G_0$ is called \emph{symmetric} if
$G_0 = -G_0$.
If $G_0 = \{e_0,e_1,e_2,e_3\}$ with $e_0=e_1+e_2+e_3$,  then it is easy to check that
$G_0 \cup (-G_0)$ is a cap, and we call $G_0 \cup (-G_0)$ the
\emph{doubled four-cap}.  

\begin{lemma}\label{lem:symmetric-caps}
Let $G = \F_3^3$ and let $G_0 \subset G^{\bullet}$ be a symmetric cap.  Then
$|G_0| \le 8$, and there is a basis $(g_1,g_2,g_3)$ of $G$ such that
\[
 G_0 \subset \{\pm g_1, \pm g_2, \pm g_3, \pm (g_1+g_2+g_3)\} ,
\]
with equality if $|G_0| = 8$.  In particular, the symmetric caps with eight
elements are precisely the images of the doubled four-cap under the
automorphisms of $G$.
\end{lemma}

\begin{proof}
Since $0 \notin G_0 = -G_0$ and $|G_0| \le 9$, it follows  that $|G_0| = 2k$ with $k \in [1,4]$, say  $G_0=\{g_1, \ldots, g_k, -g_1, \ldots, -g_k\}$.  We claim that any three of them are
linearly independent. Indeed,  a nontrivial relation between three of them has all
its coefficients nonzero, because $g_i \ne \pm g_j$ for $i \ne j$, and hence
gives three distinct elements of $G_0$ with sum zero.  If $k \le 3$, we
extend $g_1, \ldots, g_k$ to a basis $(g_1,g_2,g_3)$ of $G$.  If $k = 4$, then
$(g_1,g_2,g_3)$ is a basis and $g_4 = a_1g_1+a_2g_2+a_3g_3$ with all
$a_i \ne 0$; replacing $g_i$ by $-g_i$ if necessary, we obtain
$g_4 = g_1+g_2+g_3$.  The last assertion follows because $\Aut (G)$ acts
transitively on the bases of $G$.
\end{proof}

\subsection{Graphs.} \label{sec:graphs}

We need some basic concepts from graph theory.
An \emph{undirected graph} is a pair $\Gamma = (V, E)$ consisting of a finite set $V$ of
\emph{vertices} and a set $E$ of \emph{edges}, where each edge is a subset
$e \subset V$ with $|e| \in \{1, 2\}$.  Two distinct vertices $v, w$ are
\emph{adjacent} if $\{v,w\} \in E$; an edge
$\{v\}$ is a \emph{loop}, and we say that $v$ \emph{has a loop}.
A subset $W \subset V$ is a \emph{clique} if any two distinct vertices of $W$
are adjacent; loops play no role in this definition.    A clique is \emph{maximal} if it is maximal with respect to
inclusion.  Since $V$ is finite, every clique is contained in a (not
necessarily unique) maximal clique.

If instead $E \subset V \times V$, then every edge $(v,w)$ has an
orientation from $v$ to $w$, whence $\Gamma = (V,E)$ is a \emph{directed graph}.
A \emph{directed path} of length $\ell \in \N_0$ from $v$ to $w$ is a sequence
of vertices $v = v_0, v_1, \ldots, v_\ell = w$ with $(v_{i-1}, v_i) \in E$ for
all $i \in [1,\ell]$; a \emph{cycle} is a directed path of positive length
from a vertex to itself.  A directed graph without cycles is a
\emph{directed acyclic graph} (DAG).  DAGs in particular are useful for
visualizing factorizations; see Figure~\ref{fig:factorization-dag} in
Subsection~\ref{sec:software:intro}.

\smallskip
\section{Basic  Properties of  Zero-Sum Sequences and Length Sets} \label{3}
\smallskip

Let $G=C_3^3$. 
In Subsection \ref{sec:basic} we gather some basic arithmetic properties of the monoid $\mathcal B (G)$.
In Subsection~\ref{sec:atoms} we determine the
structure of all minimal zero-sum sequences over $G$, and in
Subsection~\ref{sec:factorization-lemmas} we collect elementary lemmas on
factorizations  which are used in Sections~\ref{software}, \ref{5}, 
and~\ref{6}.

\subsection{Summary of Basic Properties} \label{sec:basic}
 We start with some further definitions. Let $G$ be a finite abelian group. Then
\[
\rho (G) = \max \{ \rho (L) \colon \{0\} \ne L \in \mathcal L (G)  \}
\]
is the elasticity of $G$. For $k \in \N$, the $k$-th elasticity $\rho_k (G)$ is the maximum of all $\ell \in \N$ for which there is an equation of the form $U_1 \cdot \ldots \cdot U_k = V_1 \cdot \ldots \cdot V_{\ell}$, where all $U_i$ and $V_j$ are atoms.

\begin{lemma} \label{lem:basic-properties}
Let $G = C_3^3$.
\begin{enumerate}
\item  $\mathsf D (G)=7$, $\rho (G)=7/2$ and  $\rho_k (G) = \lfloor 7 k /2 \rfloor$ for all $k \ge 2$.

\item $\Delta (G)  = [1,2]$ and $\mathsf c(G)=4$.
\end{enumerate}
\end{lemma}

\begin{proof}
The assertion on $\mathsf D (G)$ follows from \eqref{inequality-Davenport}. For proofs of the remaining statements, we refer to \cite[Corollaries 5.3.18 and 6.2.14]{Ge-Gr-Zh26a}.
\end{proof}

We need these invariants in the proof of Theorem \ref{1.1}. On the other hand, regardless of the proof of Theorem \ref{1.1}, knowledge of these invariants 
 helps us gain a better understanding of Theorem \ref{1.1}. For example, we observe that the  elasticity $\rho (G)=7/2$ is attained only by the length sets of $\mathcal L_2$. This is no surprise, as  length sets in most groups, whose elasticity equals $\rho (G)$, are arithmetic progressions with difference $1$ (for work in this direction see 
\cite{Ya26a}).

\subsection{The minimal zero-sum sequences}\label{sec:atoms}

Let $G=C_3^3$ and let $U \in \mathcal A (G)$. Then $|U| \in [1, \mathsf D (G)] = [1,7]$ by Lemma \ref{lem:basic-properties}. Obviously, $U = 0$ is the only minimal zero-sum sequence of length $1$ and $U = (-g)g$ for all $g \in G^{\bullet}$ are the minimal zero-sum sequences of length $2$. If $g, h \in G^{\bullet}$ with $h \ne -g$, then $U = g h (-g-h)$ is a minimal zero-sum sequence and every minimal zero-sum sequence of length three has this form; if $h = g$, then $U = g^3$, and otherwise $g$ and $h$ are linearly independent. Since every $g \in G^{\bullet}$ and every pair of linearly independent elements of $G$ extend to a basis of $G$, these are the forms listed in Part 5 of Theorem~\ref{3.3}. The structure of the minimal zero-sum sequences of length in $[4,7]$ is more involved; we start with two propositions.

\smallskip
\begin{proposition} \label{3.1}
If  $G = C_3 \oplus C_3$, then
\[
\begin{aligned}
\mathcal A (G)  = \ \{0\} & \sqcup  \{ (- e) e :\; e \in G^{\bullet} \} \cup \{ e^3 :\; e \in G^{\bullet} \} \\
                          & \sqcup \bigl\{ e_1^{b_1}  e_2^{b_2}  (-b_1e_1 - b_2 e_2)  :\;\mbox{ $(e_1, e_2)$ a basis for  $G$ with  $b_1, b_2 \in \{1,2\}$} \bigr\}.
\end{aligned}
\]
In particular, if $T$ is a zero-sum free sequence over $G$ of length $4$, then $\mathsf h(T)=2$.
\end{proposition}

\begin{proof}
See \cite[Corollary 4.5.10]{Ge-Gr-Zh26a}.
\end{proof}

\smallskip
\begin{proposition} \label{3.2}
Let $G = C_3 \oplus C_3 \oplus C_3$.  A sequence  $U \in \mathcal F (G)$ is a minimal zero-sum sequence of length $\mathsf D (G)=7$ if and only if
there is a basis $(e_1,e_2,e_3)$ of $G$ such that
\[
U=e_1^2 \prod_{i\in [1,2]}(x_ie_1+e_2)\prod_{j\in [1,3]}(x_{j+2}e_1 + y_je_2 +  e_3) ,
\]
where $x_1,\ldots,x_5,y_1,\ldots y_3 \in [0,2]$ with $\sum_{i=1}^5x_i\equiv
  \sum_{j=1}^3y_j \equiv 1 \mod{3}$.
\end{proposition}

\begin{proof}
See \cite[Theorem 4.6.2]{Ge-Gr-Zh26a}.
\end{proof}

\medskip
\begin{theorem} \label{3.3}
Let $G = C_3 \oplus C_3 \oplus C_3$ and let $U \in \mathcal A (G)$.
\begin{enumerate}
\item $|U|=7$  if and only if there is a basis $(e_1,e_2,e_3)$ of $G$ such that $U$ can be written as one of the following forms
	  \begin{enumerate}
	  \item $U=e_1^2e_2^2e_3^2(e_1+e_2+e_3)$;
			
	  \item $U=e_1^2e_2^2e_3(e_1+e_3)(e_2+e_3)$;
			
	  \item $U=e_1^2e_2e_3(e_1+e_2)(-e_1+e_3)(e_1+e_2+e_3)$.
	  \end{enumerate}
			
\item $|U|=6$  if and only if there is a basis $(e_1,e_2,e_3)$ of $G$ such that $U$ can be written as one of the following forms
	  \begin{enumerate}
	  \item $U=e_1^2e_2^2e_3(e_1+e_2-e_3)$;
			
	  \item $U=e_1^2e_2e_3(e_1+e_2)(e_2-e_3)$;
			
	  \item $U=e_1e_2e_3(-e_1+e_2)(e_1+e_3)(-e_1+e_2+e_3)$.
	  \end{enumerate}
		
\item $|U|=5$  if and only if there is a basis $(e_1,e_2,e_3)$ of $G$ such that $U$ can be written as one of the following forms
	  \begin{enumerate}
	  \item $U=e_1^2e_2^2(e_1+e_2)$;
			
	  \item $U=e_1^2e_2e_3(e_1-e_2-e_3)$;
			
      \item $U=e_1e_2e_3(e_1-e_2)(e_1-e_3)$.
	  \end{enumerate}
		
\item $|U|=4$  if and only if there is a basis $(e_1,e_2,e_3)$ of $G$ such that $U$ can be written as one of the following forms
	  \begin{enumerate}
	  \item $U=e_1^2e_2(e_1-e_2)$;
			
	  \item $U=e_1e_2e_3(-e_1-e_2-e_3)$.
	  \end{enumerate}

\item $|U| \le 3$ and $U \ne 0$  if and only if there is a basis $(e_1,e_2,e_3)$ of $G$ such that $U$ can be written as one of the following forms
	  \begin{enumerate}
	  \item $U=(-e_1)e_1$;

	  \item $U=e_1^3$;

	  \item $U=e_1e_2(-e_1-e_2)$.
	  \end{enumerate}
\end{enumerate}
\end{theorem}

\begin{proof}
Let $U\in \mathcal A(G)$ with $U\neq 0$ and let $g \in \supp (U)$. Then $\ord (g)=3$, whence there is a subgroup $H_g \subset G$ such that $G = \langle g \rangle \oplus H_g$. Let $\varphi_g:G\rightarrow G$ be the projection  onto $H_g$, whence $\varphi_g(xg+h)=h$ for $h\in H_g$ and $x\in \Z$. Then   $$\ker \varphi_g=\la g\ra\;\und\; \im \varphi_g=\varphi_g(G)=H_g\cong C_3\oplus C_3.$$

1. This follows from Proposition \ref{3.2}.
	
2. Suppose $|U|=6$. We distinguish several cases.
	
	Suppose $\mathsf h(U)=2$ and there exist $h_1,h_2\in G$ such that $h_1^2h_2^2\t U$. Then $U=h_1^2h_2^2h_3(h_1+h_2-h_3)$ and $(h_1,h_2,h_3)$ is a basis of $G$.

	Suppose  $\mathsf h(U)=2$ and there exists $h_1\in G$ such that $h_1^2\t U$ and $(h_1^2)^{-1}U$ is squarefree.
	 We see that $\varphi_{h_1}((h_1^2)^{-1}U)$ is a minimal zero-sum sequence over $H_{h_1}$. It follows from Proposition \ref{3.1}  that $\varphi_{h_1}((h_1^2)^{-1}U)= h_2^2h_3(h_2-h_3)$ for independent $h_2,h_3\in H_{h_1}$.
	Since $\sigma((h_1^2)^{-1}U)=h_1$, there exist $a_1,\ldots,a_4\in [0,2]$ with $\sum_{j=1}^4a_j\equiv 1\mod 3$ such that
	\[
	U=h_1^2\left(\prod_{i=1}^2(a_ih_1+h_2)\right)(a_3h_1 + h_3) (a_4h_1+h_2-h_3) \,.
	\]
	Let $e_1=h_1$, $e_2=a_1h_1+h_2$, and $e_3=a_3h_1+h_3$. Then $h_2=e_2-a_1e_1$ and $h_3=e_3-a_3e_1$. Set $x=a_2-a_1$ and $y=a_4-a_1+a_3$. Then $x+y\equiv 1\mod 3$ and
	\[
	U=e_1^2 e_2 e_3 (xe_1+e_2)(ye_1+e_2-e_3) \,.
	\]
	
	If $x\equiv 0 \mod 3$, then we obtain a contradiction to the fact that $(h_1^2)^{-1}U$	is squarefree.
	If $x\equiv 1 \mod 3$, then $y\equiv 0 \mod 3$ and $U=e_1^2 e_2 e_3 (e_1+e_2)(e_2-e_3)$.
		If $x\equiv -1\mod 3$, then $y\equiv -1 \mod 3$ and $U=e_1^2 e_2 e_3 (-e_1+e_2)(-e_1+e_2-e_3)$.
	By letting $(f_1,f_2,f_3)=(e_1, -e_1+e_2, e_3)$, we obtain $U=f_1^2 f_2 f_3 (f_1+f_2)  (f_2-f_3)$.

	Suppose $\mathsf h(U)=1$. Let $h_1\t U$ and let $W=h_1^{-1}U$.
	Suppose $\varphi_{h_1}(W)$ is a minimal zero-sum sequence over $H_{h_1}$.
	Set $U^*=h_1^{-1}(-h_1)^2U$. Then $U^*$ is a minimal zero-sum sequence of length $7$, whence it has the form 1(c), that is,  there exists a basis $(e_1,e_2,e_3)$ such that $U^*=e_1^2e_2e_3(e_1+e_2)(-e_1+e_3)(e_1+e_2+e_3)$, whence
	$U=(-e_1)e_2e_3(e_1+e_2)(-e_1+e_3)(e_1+e_2+e_3)$.

	Suppose $\varphi_{h_1}(W)$ is not a minimal zero-sum sequence over $H_{h_1}$.
	Then $\varphi_{h_1}(W)=g(-g)g_1g_2(-g_1-g_2)$ for some $g,g_1,g_2\in H_{h_1}$ and there exist $a_1,\ldots,a_5\in [0,2]$ with $a_1+a_2\equiv a_3+a_4+a_5\equiv 1\mod 3$ such that
	\[
	U=h_1(a_1h_1+g)(a_2h_1-g)(a_3h_1+g_1)(a_4h_1+g_2) (a_5h_1-g_1-g_2) \,.
	\]
	Let $h_2=a_1h_1+g$. There exists $i\in [1,2]$ such that $g_i\not\in \langle h_1,h_2 \rangle=\langle h_1,g\rangle$, say $i=1$, whence $H_{h_1}=\langle g, g_1\rangle$. Let $h_3=a_3h_1+g_1$ and $b,c\in [0,2]$ such that $g_2=bg+cg_1=bh_2+ch_3-(a_1b+a_3c)h_1$, and let $d\equiv a_4-a_1b-a_3c \mod 3$. Then
	\[
	U=h_1h_2(h_1-h_2)h_3(dh_1+bh_2+ch_3) ((1-d)h_1-bh_2+(-1-c)h_3) \,.
	\]	
	If $c\neq 1$, then either $c\equiv 0\mod 3$, or $-1-c\equiv 0 \mod 3$, whence there exists a zero-sum free subsequence $T$ of $U$ such that $|T|=4$, $T\in \mathcal F(\langle h_1, h_2\rangle)$, and $\mathsf h(T)=1$. But this is impossible from Proposition \ref{3.1}.
	If $c=1$, then \[
	U=h_1h_2(h_1-h_2)h_3(dh_1+bh_2+h_3) ((1-d)h_1-bh_2+h_3) \,.
	\]
	Since $\varphi_{h_2}(h_2^{-1}U)=h_1^2h_3(dh_1+h_3)((1-d)h_1+h_3)$ is a minimal zero-sum sequence over $H_{h_2}=\langle h_1,h_3\rangle$, we go back to a case above.

3. Suppose $|U|=5$. If  $\langle \supp(U)\rangle\cong C_3^2$, then $U=e_1^2e_2^2(e_1+e_2)$, where $(e_1,e_2,e_3)$ is a basis of $G$.
	If $\langle \supp(U)\rangle\cong C_3^3$ and $\mathsf h(U)=2$, then $U=e_1^2e_2e_3(e_1-e_2-e_3)$, where $(e_1,e_2,e_3)$ is a basis of $G$.
	If $\langle \supp(U)\rangle\cong C_3^3$ and $\mathsf h(U)=1$, then there exists a basis $(f_1,f_2,f_3)$ such that $U=f_1f_2f_3(x_1f_1+x_2f_2+x_3f_3)(x_1'f_1+x_2'f_2+x_3'f_3)$, where $x_1,x_2,x_3,x_1',x_2',x_3'\in [-1,1]$ such that $x_1+x_1'\equiv x_2+x_2'\equiv x_3+x_3'\equiv -1\mod 3$. Note that $x_i+x_i'=2$ implies that $x_i=x_i'=1$ for each $i\in [1,3]$. Thus $\mathsf h(U)=1$ implies that  $$|\{i\in[1,3]\colon x_i+x_i'=2\}|\le 2\,.$$
	
	Suppose $|\{i\in[1,3]\colon x_i+x_i'=2\}|= 2$. Then, by symmetry, we may assume that $x_1=x_1'=x_2=x_2'=1$, $x_3=0$, and $x_3'=-1$. By letting $(e_1,e_2,e_3)=(f_1+f_2, f_1+f_2-f_3, f_2)$, we obtain that $U=e_1e_2e_3(e_1-e_2)(e_1-e_3)$.
	Suppose $|\{i\in[1,3]\colon x_i+x_i'=2\}|= 1$. Then, by symmetry, we may assume that $x_1=x_1'=1$, $x_2=0$, and $x_2'=-1$. If $x_3=0$, then $(x_1f_1+x_2f_2+x_3f_3)=f_1$, a contradiction. Thus $x_3=-1$, $x_3'=0$, and  $U=f_1f_2f_3(f_1-f_3)(f_1-f_2)$. Suppose $|\{i\in[1,3]\colon x_i+x_i'=2\}|= 0$. Then $|\{i\in[1,3]\colon x_i+x_i'=-1\}|= 3$. Since $0\nmid U$, by symmetry, we may assume that $x_1=x_2=x_3'=0$ and $x_1'=x_2'=x_3=-1$, whence $-f_3\mid U$, a contradiction.

4. Suppose $|U|=4$. If  $\langle \supp(U)\rangle\cong C_3^3$, then there exists  a basis $(e_1,e_2,e_3)$ of $G$ such that  $U=e_1e_2e_3(-e_1-e_2-e_3)$.	If $\langle \supp(U)\rangle\cong C_3^2$, then  $\mathsf h(U)=2$ and there exists a basis $(e_1,e_2,e_3)$ of $G$ such that $U=e_1^2e_2(e_1-e_2)$.

5. See the discussion at the beginning of this subsection.
\end{proof}

The next lemma characterizes those atoms of $\mathcal A (G)$ whose support lies in a  doubled four-cap (see Subsection~\ref{2.5}).

\begin{lemma}\label{lem:cap-atoms}
Let $(e_1,e_2,e_3)$ be a basis of $G$, let $e_0 = e_1+e_2+e_3$, let
$e_4 = -e_0$, and let $G_0 = \{e_0,e_1,e_2,e_3\}$.
\begin{enumerate}
\item If $U \in \mathcal A (G_0)$, then $|U| \in \{3,5,7\}$.

\item Let $U \in \mathcal A \big(  G_0 \cup (-G_0) \big)$ and suppose that $U$ is
      neither of the form $g(-g)$ nor of the form $g^3$ for any $g \in \supp (U)$.
      \begin{enumerate}
      \item Then
            \[
            U = \prod_{i=1}^{4}(\varepsilon_i e_i)^{n_i}
            \qquad \text{with} \qquad \varepsilon_i \in \{1,-1\}, \ n_i \in [1,2] ,
            \]
            $\varepsilon_i = \varepsilon_j$ if $n_i = n_j$, and $\varepsilon_i = -\varepsilon_j$ if $n_i \ne n_j$ for all $i,j \in [1,4]$.  In particular $|U| \in [4,7]$.
      
      \item If $|U| = 7$, then there are $j \in [1,4]$ and
            $\varepsilon \in \{1,-1\}$ such that
            \[
            U =  (-\varepsilon e_j)  \prod_{i \in [1,4] \setminus \{ j\}}(\varepsilon e_i)^2 ,
            \]
            and conversely each of these eight sequences is an atom.
            
      \item If $|U| = 4$, then $U$ has one or three terms in $-G_0$, and if $|U| = 5$, then $U$ has zero, two, three or five terms in $-G_0$.
      \end{enumerate}      
\end{enumerate}
\end{lemma}

\begin{proof}
1. Let
$U = e_0^{n_0}e_1^{n_1}e_2^{n_2}e_3^{n_3} \in \mathcal A (G_0)$ with $n_0, \ldots, n_3 \in \N_0$. If there is some $i \in [0,3]$ with $n_i \ge 3$, then $n_i=3$ and $U =e_i^3$. Suppose that $n_0, \ldots, n_3 \in [0,2]$.
Since $0 = \sigma (U) = \sum_{i=1}^3 (n_i+n_0)e_i$, it follows that $n_i \equiv -n_0 \mod 3$ for
$i \in [1,3]$.  Thus $n_0 = 1$ and $n_i = 2$, or
$n_0 = 2$ and $n_i = 1$ for $i \in [1,3]$, whence $|U| \in \{7,5\}$.

2.(a)  Since $U$ has none of the indicated forms, we may write
$U = \prod_{i=1}^{4}(\varepsilon_i e_i)^{n_i}$ with $\varepsilon_i \in \{1,-1\}$
and $n_i \in [0,2]$ for all $i \in [1,4]$.   Since 
$0 = \sigma (U) = \sum_{i=1}^3(\varepsilon_i n_i - \varepsilon_4 n_4)e_i$, it follows that 
$\varepsilon_i n_i \equiv \varepsilon_4 n_4 \mod 3$ for all $i \in [1,3]$.  If one of the  $n_i$ would be
$0$, then all $n_i$ would be $0$,  whence $n_i \in [1,2]$ for all $i \in [1,4]$, and all assertions follow.

2.(b)  With the notation of 2.(a), we get  $n_1+n_2+n_3+n_4=|U| = 7$, whence there is some $j \in [1,4]$ with $n_j=1$ and $n_i=2$ for all $i \in [1,4]\setminus \{j\}$.  This implies that  $\varepsilon_i = -\varepsilon_j$ for
all $i \ne j$, which is the displayed form with $\varepsilon = \varepsilon_i$.
Conversely, every sequence of the given form is an atom because $(e_1,e_2,e_3)$ is a basis.

2.(c)  We use the notation of 2.(a).  If $|U| = 4$, then 
$n_i = 1$ for all $i \in [1,4]$, all signs agree, and $U$ has one
resp.\ three terms in $-G_0$.  Suppose $|U| = 5$. Then
 $U = (-\varepsilon e_j)^2 \prod_{i \ne j} (\varepsilon e_i)$ for some $j \in [1,4]$.  If $j = 4$, then
 $U=e_0^2e_1e_2e_3$ or $U=(-e_0)^2(-e_1)(-e_2)(-e_3)$, with zero or five
terms in $-G_0$.  If $j \in [1,3]$ and $\varepsilon = 1$, then $U$ has the three terms
 $-e_j, -e_j, e_4$ in $-G_0$, and if $\varepsilon = -1$, then two terms of $U$ are in $-G_0$.
\end{proof}

\subsection{Elementary Lemmas on Factorizations}\label{sec:factorization-lemmas}

By \eqref{factorization-hom} and the abbreviations introduced in Subsection \ref{2.3}, 
\[
\mathsf Z (G) = \mathcal F \big(\mathcal A(G) \big)
 \qquad \text{and} \qquad
 \pi\colon \mathsf Z (G)\longrightarrow\mathcal B(G)
\]
denote the  factorization monoid and the  factorization homomorphism of $\mathcal B (G)$. 
Since $\Delta (G) = \{1,2\}$ by Lemma \ref{lem:basic-properties}, a length set  $L  \in \mathcal L (G)$ is an
arithmetic progression with difference two (possibly a singleton), or
$1 \in \Delta ( L )$.

\begin{lemma} \label{lem:hereditary-gap}
Let $B \in \mathcal B (G)$, let $z \in \mathsf Z(B)$, and let $y \in \mathsf Z (G)$ with $y\mid z$. 
\begin{enumerate}
\item If $|z|+1 \notin \mathsf L (B)$, then $|y|+1\notin \mathsf L \big(\pi(y)   \big)$. If in addition $|z| = \min \mathsf L (B)$, then also
     $\min\mathsf L \big(\pi(y) \big)=|y|$.  In particular, if $z=U_1\cdot\ldots\cdot U_k$
     and if $i,j,h \in [1,k]$ are distinct, then
     \[
     3\notin\mathsf L(U_iU_j), \qquad 4\notin\mathsf L(U_iU_jU_h) \,.
     \]
     
\item If $|z| = \max \mathsf L (B)$ and $|z|-1 \notin \mathsf L (B)$,  then
      $|y| = \max \mathsf L \big(\pi (y) \big)$ and $|y|-1 \notin \mathsf L \big(\pi (y) \big)$.
      In particular, if $z = U_1 \cdot \ldots \cdot U_k$ and $i,j,h \in [1,k]$ are distinct, then
      \[
      \mathsf L (U_iU_j)=\{2\}, \qquad \mathsf L (U_iU_jU_h)=\{3\} .
      \]
\end{enumerate}     
\end{lemma}

\begin{proof}
Let $x \in \mathsf Z (G)$ with $z = xy = U_1 \cdot \ldots \cdot U_k$ with $U_1, \ldots, U_k \in \mathcal A (G)$. Then $B = \pi (z) = \pi (x)\pi (y)$ and $k = |z|=|x|+|y| \in \mathsf L \big( \pi (x) \big) + \mathsf L \big( \pi (y) \big) \subset \mathsf L (B)$.
Thus both claims follow.
\end{proof}

\begin{lemma} \label{lem:upper-hole-affine-cap}
Let $B \in \mathcal B (G)$ with 
$\max\mathsf L(B)-1\notin\mathsf L(B)$. If  $B = U_1 \cdot \ldots \cdot U_k$ with $k=\max \mathsf L (B)$ , then $|\{i \in [1,k] \colon |U_i|=2\}| \le 8$. Moreover, if $|\{i \in [1,k] \colon |U_i|=2\}| = 8$, then there exists a basis $(e_1,e_2,e_3)$ of $G$ such that
\[
\prod_{i\in [1,k] \text{ with }|U_i|=2}U_i \quad = \quad \big(e_1(-e_1) \big)^2 \big(e_2(-e_2)\big)^2 \big(e_3(-e_3)\big)^2 \big((e_1+e_2+e_3)(-e_1-e_2-e_3)\big)^2\,.
\] 
\end{lemma}

\begin{proof}
Let $G_0 = \{ g \in G^{\bullet} \colon (-g)g = U_i \ \text{for some $i \in [1,k]$} \}$. Note that $-G_0 = G_0 \subset G^{\bullet}$, whence $|G_0|$ is even.  Since  $\big( (-g)g \big)^3 = g^3 (-g)^3$ and $\max\mathsf L(B)-1\notin\mathsf L(B)$,  $(-g)g$ does not occur three times in the factorization for any $g \in G_0$. We assert that  $G_0$ is a cap in $\F_3^3$. Indeed, if $g_1, g_2, g_3 \in G_0$ have sum zero, then 
\[
\big( (-g_1)g_1 \big) \big( (-g_2)g_2 \big) \big( (-g_3)g_3 \big) = \big(g_1g_2g_3 \big) \big((-g_1)(-g_2)(-g_3) \big)
\]
implies that $\max \mathsf L (B)-1 \in \mathsf L (B)$, a contradiction. By Subsection~\ref{2.5}, a cap in $\F_3^3$
has at most nine elements, whence $|G_0|\le 8$. Since each pair $(-g)g$ occurs at most twice, we infer that $|\{i \in [1,k] \colon |U_i|=2\}| \le 8$. If $|\{i \in [1,k] \colon |U_i|=2\}| = 8$, then each pair $(-g)g$ occurs twice, whence
$|G_0| = 8$.  Thus the assertion follows by  Lemma~\ref{lem:symmetric-caps}.
\end{proof}

\begin{lemma} \label{lem:parity-formula}
Let $B \in \mathcal B (G)$ and let
$z = W_1 \cdot \ldots \cdot W_{\ell} \in \mathsf Z (B)$ with
$W_1, \ldots, W_{\ell} \in \mathcal A (G)$.
\begin{enumerate}
\item $\ell \equiv |B| + \bigl|\{ i \in [1,\ell] \colon |W_i| \ \text{is even}\}\bigr| \mod 2$.
\item $1 \notin \Delta \big(\mathsf L (B) \big)$ if and only if any two factorizations
      of $B$ have the same number of atoms of even length modulo two.
\item Let $G_0 \subset G^{\bullet}$ be such that every atom of
      $\mathcal B (G_0)$ has odd length.  Then
      $\mathsf L (S) \subset |S| + 2\Z$ for every $S \in \mathcal B (G_0)$.
      In particular, $1 \notin \Delta \big(\mathsf L (S) \big)$.
\end{enumerate}
\end{lemma}

\begin{proof}
1. We have $|B| = |W_1| + \ldots + |W_{\ell}|$, and $1 \equiv |W_i| \mod 2$ if
$|W_i|$ is odd, $1 \equiv |W_i| + 1 \mod 2$ if $|W_i|$ is even.  Summing
over $i \in [1,\ell]$ gives the assertion.

2. Part 1 shows that the lengths of two factorizations are congruent modulo $2$ if and
only if the number of atoms of even length in the respective factorizations agree modulo two. This implies the assertion.

3. Let $S \in \mathcal B (G_0)$ and let $z = U_1 \cdot \ldots \cdot U_k \in \mathsf Z (S)$.  Then $|U_i| \equiv 1 \mod 2$ for every $i \in [1,k]$.  Thus Part 1 implies that 
$k=|z| \equiv |S| \mod 2$, whence $\mathsf L (S) \subset |S| + 2\Z$.
\end{proof}

The following criterion refines Part 3 of Lemma~\ref{lem:parity-formula}. If every atom over a subset $G_0 \subset G^{\bullet}$ has odd
length, then the parity of the length of a factorization $y$ depends only on those
atoms dividing $y$ that are not in  $\mathcal A (G_0)$.

\begin{lemma} \label{lem:parity-criterion}
Let $G_0 \subset G^{\bullet}$ be a subset such that every atom of $\mathcal B (G_0)$ has
odd length, let $A \in \mathcal B (G)$, let
$C_1, \ldots, C_r \in \mathcal A (G_0)$ with
$\supp (C_1 \cdot \ldots \cdot C_r) = G_0$, and put
$B_m = A C_1^{m_1} \cdot \ldots \cdot C_r^{m_r}$ for
$m = (m_1, \ldots, m_r) \in \N_0^r$.  Let $\mathsf Y$ be the set of all
$y \in \mathsf Z (G)$ such that every atom dividing $y$ has a term outside $G_0$,
and $\pi(y)$ and $A$ have the same terms outside $G_0$ (with the same
multiplicities), that is, $\pi(y)_{G\setminus G_0} = A_{G\setminus G_0}$.
Then $\mathsf Y$ is finite, and
for every $m \in \N_0^r$ we have $1 \in \Delta \big(\mathsf L (B_m)\big)$ if and
only if there are $y, y' \in \mathsf Y$ with
$|y| + |\pi(y)| \not\equiv |y'| + |\pi(y')| \mod 2$ and
\begin{equation}
 \sum_{i=1}^{r} m_i \mathsf v_g (C_i)
 \ \ge\ \max \{ \mathsf v_g (\pi(y)), \mathsf v_g (\pi(y')) \} - \mathsf v_g (A)
\qquad \text{for all } g \in G_0 .\label{eq:cone-feasibility}
\end{equation}
In particular, if $|y| + |\pi(y)|$ has the same parity for all
$y \in \mathsf Y$, then $1 \notin \Delta \big(\mathsf L (B_m)\big)$ for every
$m \in \N_0^r$.
\end{lemma}

\begin{proof}
Every atom occurring in an element $y \in \mathsf Y$ has a term outside $G_0$, so
$|y| \le |A_{G \setminus G_0}|$, and $\mathsf Y$ is finite.  Since
$C_1, \ldots, C_r \in \mathcal F (G_0)$, we have
$(B_m)_{G \setminus G_0} = A_{G \setminus G_0}$ and
$\mathsf v_g (B_m) = \mathsf v_g (A) + \sum_{i=1}^r m_i \mathsf v_g (C_i)$ for
$g \in G_0$.  Hence, for $y \in \mathsf Y$, we have $\pi(y) \mid B_m$ if and only
if $\sum_{i=1}^r m_i \mathsf v_g (C_i) \ge \mathsf v_g (\pi(y)) - \mathsf v_g (A)$
for all $g \in G_0$.

If $z \in \mathsf Z (B_m)$, write $z = yw$, where no atom dividing $y$ lies
in $\mathcal A (G_0)$ and every atom dividing $w$ does; then $y \in \mathsf Y$
and $\pi(y) \mid B_m$;
conversely, if $y \in \mathsf Y$ with $\pi(y) \mid B_m$, then
$B_m\pi(y)^{-1} \in \mathcal B (G_0)$, and $y$ times any factorization of
$B_m\pi(y)^{-1}$ is a factorization of $B_m$.  In both cases, $y^{-1}z$ is a
factorization of $B_m\pi(y)^{-1} \in \mathcal B (G_0)$, so its length is
congruent to $|B_m| - |\pi(y)|$ modulo $2$ by 
Lemma~\ref{lem:parity-formula}.3, whence
$|z| \equiv |B_m| + |y| + |\pi(y)| \mod 2$.  Since $\Delta (G) = \{1,2\}$ by
Lemma~\ref{lem:basic-properties}, we have $1 \in \Delta \big(\mathsf L (B_m)\big)$
if and only if $\mathsf L (B_m)$ contains integers of both parities, and the
assertion follows.
\end{proof}

Let  $w \colon G \to (\Q, +)$ be  a (weight) map. Then $w$ lifts to a homomorphism $w \colon \mathcal F (G) \to \mathcal F (\Q)$ (see Subsection \ref{2.3}), and for $S = g_1 \cdot \ldots \cdot g_\ell$ we have $\sigma \big( w (S) \big) = w(g_1) + \ldots + w(g_\ell)$.

\begin{lemma} \label{lem:weight-bound}
Let $B \in \mathcal B (G)$ and let $w \colon G \to \Q$ be a map.
\begin{enumerate}
\item If $\sigma \big(w(W) \big) \ge 1$ for every $W \in \mathcal A (G)$ with $W \mid B$, then
      $\max \mathsf L (B) \le \sigma \big(w(B) \big)$.
\item If $\sigma \big( w(W) \big) \le 1$ for every $W \in \mathcal A (G)$ with $W \mid B$, then
      $\min \mathsf L (B) \ge \sigma \big(w(B) \big)$.
\end{enumerate}
\end{lemma}

\begin{proof}
If $B = W_1 \cdot \ldots \cdot W_\ell$ with $W_1, \ldots, W_\ell \in \mathcal A (G)$,
then $\sigma \big( w(B) \big)  = \sigma \big( w(W_1) \big)  + \ldots +  \sigma \big(w(W_\ell) \big)$, which is at least $\ell$ under the
hypothesis of Part 1 and at most $\ell$ under the hypothesis of Part 2.
\end{proof}

\begin{lemma}\label{A:short-atom}
	Let $B \in \mathcal B (G)$ and let $V \in \mathcal A (G)$ with $V \mid B$ and
	$|V| \le 3$.  If $\max\mathsf L(B)-1\notin\mathsf L(B)$, then there is some $z \in \mathsf Z (B)$ with $|z|=\max \mathsf L (B)$ such that $V \mid z$, 
	and $\max\mathsf L(B)=1+\max\mathsf L(BV^{-1})$.
\end{lemma}

\begin{proof}
Let $B = U_1 \cdot \ldots \cdot U_k$ with $k =\max \mathsf L (B)$ and $U_1, \ldots, U_k \in \mathcal A (G)$. After renumbering if necessary there is $j \in [1,3]$ such that $V \mid U_1 \cdot \ldots \cdot U_j$, say $W = V^{-1}U_1 \cdot \ldots \cdot U_j$. If $j=1$, then $V=U_1$ and we are done. Thus $j \in [2,3]$ and $W \ne 1$. Since $k= \max \mathsf L (B)$, we infer that $1+\max \mathsf L (W) \le \max \mathsf L (U_1 \cdot \ldots \cdot U_j) = j$. If $1 \le \max \mathsf L (W) \le j-2\le 1$, then $j=3$, $W \in \mathcal A (G)$, and $B = VWU_4 \cdot \ldots \cdot U_k$ has a factorization of length $k-1$, a contradiction. Thus, $\max \mathsf L (W)=j-1$, say $W = W_1 \cdot \ldots \cdot W_{j-1}$ with $W_1, \ldots, W_{j-1} \in \mathcal A (G)$. Therefore, we get $z=VW_1 \cdot \ldots \cdot W_{j-1}U_{j+1} \cdot \ldots \cdot U_k \in \mathsf Z (B)$, $V \mid z$, and $\max \mathsf L ( W_1 \cdot \ldots \cdot W_{j-1}U_{j+1} \cdot \ldots \cdot U_k) = k-1$.
\end{proof}

\smallskip
\section{Software for Studying Length Sets and Computational Results} \label{software}
\smallskip

Some steps in the proof of Theorem~\ref{1.1} reduce to finite calculations
over $C_3 \oplus C_3 \oplus C_3$, for example the characterization of
the length sets of all products of two atoms. In this section
we describe how we carried out these calculations using the newly developed
Python software package \code{zero-sum-sequences}.
We want to emphasize that the package is a general purpose tool
for computations with zero-sum sequences; it is not tailored to
our particular investigation here. A guided tour through the capabilities of our
package can be found in Subsection~\ref{sec:software:intro}.

The structural results obtained based on our computations are then
stated and proved in Subsections~\ref{sec:software:results}
to~\ref{sec:software:length-two-atoms}. The databases produced by our
calculations, together with the code that generates them and the
additional verification programs are distributed in
the Zenodo record accompanying this article~\cite{Ge-Ha-Zh27b}.
Details on what this record contains and how the databases are verified
are provided in Subsection~\ref{sec:software:structural-certificates}.

\subsection{Introduction to our Python package}\label{sec:software:intro}

While working on this project, the number of ad-hoc computations for various
special cases grew rather quickly; finding a suitable way to share our
results with others that allowed reproducing them in a simple and
straightforward way became a priority. The consolidation of our individual
computations, together with the desire to have a convenient and efficient
way to generate, inspect, manipulate, factorize, ... -- and store
zero-sum sequences over generic finite abelian groups was the main motivation
behind building the \code{zero-sum-sequences} Python package.

For the sake of convenience there is an optional
interface with the free and open source computer mathematics system
SageMath~\cite{sagemath:10.9} that allows using SageMath's finite abelian
groups directly.
Our package is free software (published
under the MIT license), available from the Python Package Index
and on GitHub at \url{https://github.com/behackl/zero-sum-sequences}.

The computations made throughout this article have been carried out
using version \code{0.3.1} of the package, which is installed by
\begin{center}
	\code{pip install zero-sum-sequences==0.3.1}
\end{center}
and requires Python 3.12 or newer.  The documentation at
\url{https://zero-sum-sequences.readthedocs.io} contains an extensive
introduction as well as a detailed reference manual of all its classes
and functions.

In particular, our package has an \emph{interactive tutorial}
at \url{https://zero-sum-sequences.readthedocs.io/en/stable/tutorial/edit/} that
runs directly in a web browser.
Let us now walk through some illustrative examples
of using the \code{zero-sum-sequences} package; a broader overview of the
features of the package follows afterwards.

\noindent
\textbf{Quickstart.}
Consider the group $G = C_2 \oplus C_4$. As is common
practice, the elements of $G$ are written as tuples,
and sequences are written multiplicatively.

\begin{lstlisting}[style=pythonconsole, caption={The atoms of $G$.}]
>>> import zero_sum_sequences as zss
>>> G = zss.FiniteAdditiveGroup.cyclic_product(2, 4)
>>> S = zss.AdditiveSequenceSpace(G)
>>> S.davenport_bound
PYOUT:5
>>> atoms = S.enumerate_atom_catalogue()
>>> len(atoms), atoms[0], atoms[-1]
PYOUT:(38, (0, 1) · (0, 3), (0, 3) · (1, 1)^3 · (1, 2))
\end{lstlisting}

The object \code{S} (the programmatic equivalent of the sequence monoid
$\mathcal{F}(G)$)
constructs sequences over $G$: \code{S(g, h, ...)} is the sequence
$g\cdot h \cdot \ldots$, and sequences are multiplied with \code{*}
and raised to powers with \code{**} (as is usual in Python). Atoms
have length at most $\mathsf{D}(C_2 \oplus C_4) = 5$. Inspecting the
output of \code{list(atoms)} yields a complete list of
38 atoms. Its first and last entries are
$(0, 1) \cdot (0, 3)$ and $(0, 3) \cdot (1, 1)^3 \cdot (1, 2)$, respectively.
Let us next compute some factorizations.

\begin{lstlisting}[style=pythonconsole, caption={A length set and its factorizations.}]
>>> x, xbar, y, z = S((0, 1)), S((0, 3)), S((1, 0)), S((1, 2))
>>> B = x**2 * xbar**4 * y**3 * z**3
>>> B.is_zero_sum(), len(B)
PYOUT:(True, 12)
>>> sorted(B.length_set())
PYOUT:[3, 4, 5]
>>> B.factorization_witnesses().get(3)
PYOUT:((0, 1)^2 · (1, 0) · (1, 2), (0, 3)^2 · (1, 0) · (1, 2), (0, 3)^2 · (1, 0) · (1, 2))
>>> len(list(B.factorizations()))
PYOUT:3
\end{lstlisting}

Mathematically, this snippet shows that for $x = (0,1)$, $\bar{x} = (0,3)$,
$y = (1,0)$, $z = (1,2)$, the zero-sum sequence $B = x^2 \bar{x}^4 y^3 z^3$
has $\mathsf{L}(B) = \{3, 4, 5\}$ with factorizations
\[
 B = \bigl(x^2yz\bigr)\bigl(\bar x^2yz\bigr)^2
   = \bigl(x^2yz\bigr)\bigl(\bar x^4\bigr)\bigl(y^2\bigr)\bigl(z^2\bigr)
   = \bigl(x\bar x\bigr)^2\bigl(\bar x^2yz\bigr)\bigl(y^2\bigr)\bigl(z^2\bigr).
\]

Proper support for termwise actions of automorphisms of $G$ on sequences
is another core feature of our package. In the snippet below we show
that two given length-5 atoms $U$ and $V$ lie in the same orbit
and determine an automorphism by providing the images of the generators.

\begin{lstlisting}[style=pythonconsole, caption={Working with automorphism orbits.}]
>>> U = x**3 * y * S((1, 1))
>>> V = x * y * S((1, 1))**3
>>> U.is_atom(), V.is_atom(), len(U.orbit())
PYOUT:(True, True, 8)
>>> U.orbit_witness(V).show()
PYOUT:(1, 0) ↦ (1, 0)
PYOUT:(0, 1) ↦ (1, 1)
>>> sorted((U * -U).length_set())
PYOUT:[2, 4, 5]
\end{lstlisting}

\noindent
\textbf{Package Overview.}
The documentation at \url{https://zero-sum-sequences.readthedocs.io}
describes the full feature set of the package in detail. For this article,
the following parts of the package have been used.

\begin{itemize}
\item \emph{Groups and sequences.} The class \code{FiniteAdditiveGroup}
  is the central interface for setting up the finite additive group $G$.
  The class method \code{cyclic_product}
  allows easily constructing products of cyclic groups.
  As mentioned above, the
  monoid $\mathcal{F}(G)$ corresponds to an instance of
  \code{AdditiveSequenceSpace}. 
  Each ``sequence space'' knows an upper bound for $\mathsf{D}(G)$.
\item \emph{Atom catalogues.}  \code{AtomCatalogue} is the list of all atoms
  of $G$ (or of the atoms with support in a given subset), with fixed
  indices, its orbits under the automorphism group, and the atoms dividing a
  given sequence.  Catalogues can be written to and read from files.
\item \emph{Factorizations.}  \code{FactorizationSolver} answers all
  questions about the factorizations of one zero-sum sequence: its length
  set, its minimal and maximal length, whether a given length occurs, a
  factorization of a given length (possibly with a prescribed number of
  atoms of a given kind), and all factorizations.  The methods of a
  sequence used above call it.  \code{FactorizationCache} stores the answers
  for many sequences, optionally only for one representative of each orbit.
\item \emph{Symmetry.}  \code{AutomorphismGroup} computes orbits, canonical
  representatives, stabilizers and transporters of sequences.
\item \emph{Zero-sum constants.}  The \code{zero_sum_constant} function
  computes the Davenport constant and related constants by exhaustive search.
\end{itemize}

\noindent
\textbf{Factorization algorithm.}
Let us now outline the theoretical considerations behind the implementation
of our efficient factorization algorithm.

Let $B \in \mathcal B (G)$ with $0 \nmid B$; sequences with $0$-occurrences
are reduced by applying~\eqref{shift-by-zeros} beforehand.
We then first determine all atoms
dividing $B$, either by checking against a given atom catalogue, or otherwise
by exhaustively checking all zero-sum subsequences of $B$ whose length is
at most the known bound for $\mathsf{D}(G)$.

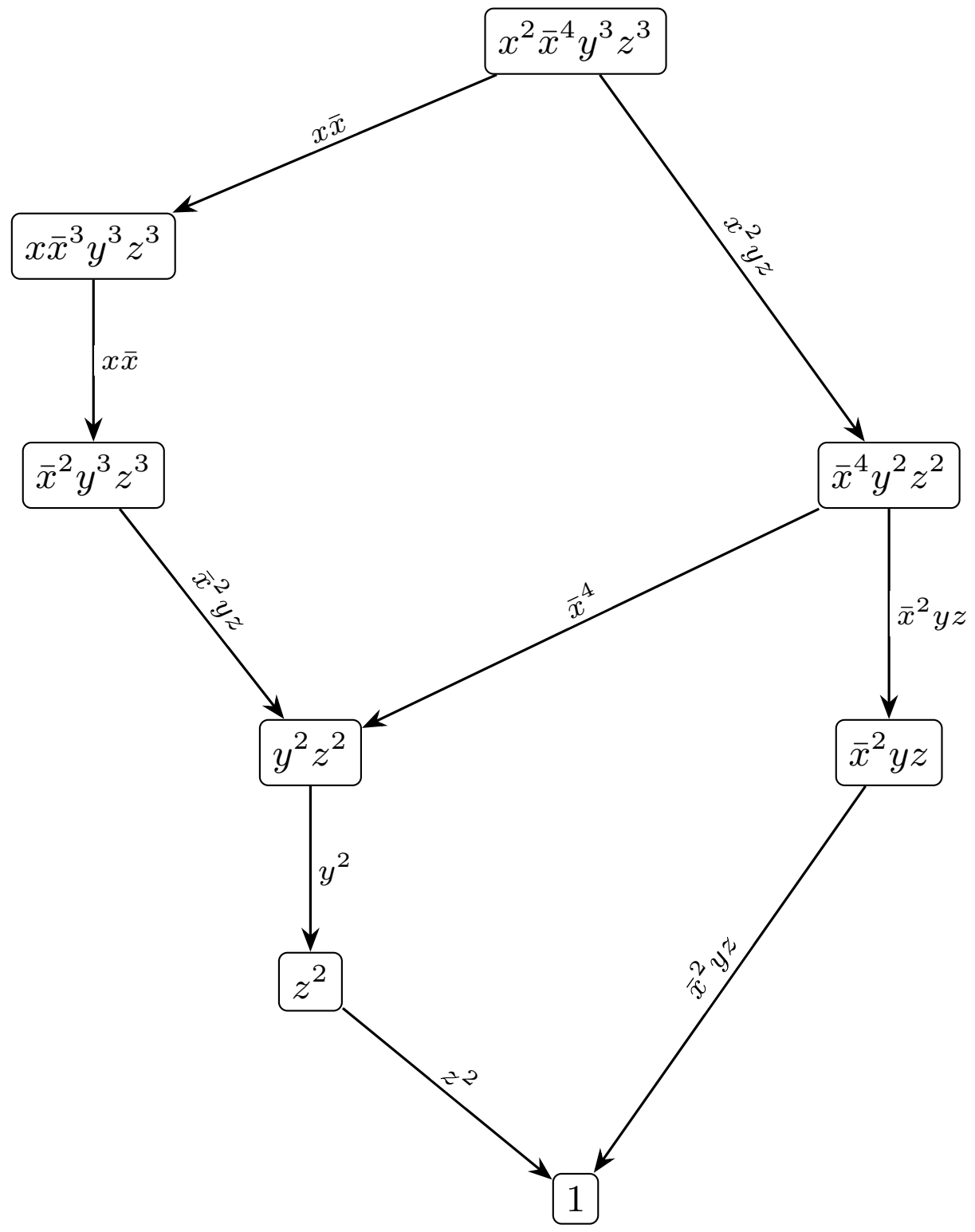
\begin{figure}
   \centering
   \begin{tikzpicture}[
     >=Stealth,
     state/.style={
       draw,
       rounded corners=2pt,
       inner sep=3pt,
       font=\small
     },
     atom/.style={
       fill=white,
       inner sep=1.5pt,
       font=\scriptsize
     },
     every edge/.style={draw, ->, semithick}
   ]
     \node[state] (S)  at ( 0,   0.0) {$x^2\bar{x}^4y^3z^3$};
     \node[state] (R1) at (-4,  -1.7) {$x\bar{x}^3y^3z^3$};
     \node[state] (R2) at (-4,  -3.6) {$\bar{x}^2y^3z^3$};
     \node[state] (T)  at ( 2.6,-3.6) {$\bar{x}^4y^2z^2$};
     \node[state] (U)  at (-2.2,-5.9) {$y^2z^2$};
     \node[state] (Q)  at ( 2.6,-5.9) {$\bar{x}^2yz$};
     \node[state] (V)  at (-2.2,-7.8) {$z^2$};
     \node[state] (E)  at ( 0,  -9.6) {$1$};

     \path
       (S)  edge node[atom, above, sloped] {$x\bar{x}$} (R1)
       (S)  edge node[atom, above, sloped] {$x^2yz$} (T)

       (R1) edge node[atom, right] {$x\bar{x}$} (R2)
       (R2) edge node[atom, above, sloped] {$\bar{x}^2yz$} (U)

       (T)  edge node[atom, above, sloped] {$\bar{x}^4$} (U)
       (T)  edge node[atom, right] {$\bar{x}^2yz$} (Q)

       (U)  edge node[atom, right] {$y^2$} (V)
       (V)  edge node[atom, above, sloped] {$z^2$} (E)
       (Q)  edge node[atom, above, sloped] {$\bar{x}^2yz$} (E);
   \end{tikzpicture}
   \caption{The factorization DAG of
     $B=x^2\bar{x}^4y^3z^3$ over $C_2\oplus C_4$,
	 where $x = (0, 1)$, $\bar{x} = (0, 3)$, $y = (1, 0)$, and $z = (1, 2)$.
     Edges are labelled by the removed atoms. The three directed paths
     from $B$ to $1$ have lengths $3$, $4$, and $5$, which corresponds to $\mathsf{L}(B) = \{3, 4, 5\}$.}
	 \label{fig:factorization-dag}
\end{figure}

Then, with a fixed order of the
terms of $B$, we consider the \emph{remainders} (i.e., the divisors $R$ of $B$
obtained by removing atoms). For a remainder $R \neq 1$ of $B$, let $g_R$ denote
the first term of $R$. Then every factorization of $R$ contains an atom $A$ with
$g_R \mid A$, and the remaining atoms form a factorization of $A^{-1}R$. Therefore,
\[
 \mathsf L (1) = \{0\}
 \qquad \text{and} \qquad
 \mathsf L (R) = \bigcup_{\substack{A \in \mathcal A (G),\ A \mid R \\ g_R \mid A}} \bigl(1 + \mathsf L (A^{-1}R)\bigr)
 \quad \text{for } R \ne 1 .
\]
The solver builds the directed graph whose vertices are the remainders and
whose edges join $R$ to $A^{-1}R$ for the atoms $A$ in this union, labelled
by $A$ (see Subsection~\ref{sec:graphs}).  The graph is acyclic, since
every edge lowers the length, and $\mathsf L (R)$ is the set of lengths of
the directed paths from $R$ to $1$. See also Figure~\ref{fig:factorization-dag}
for an example.  The solver computes $\mathsf L (R)$ for all vertices,
from $1$ upwards, by the recursion above.

Two features make this fast:
First, a remainder reached along different paths is a single vertex, so the
work is proportional to the number of remainders and edges, not to the
number of factorizations.  Second, requiring $g_R \mid A$ restricts the
edges leaving $R$ to the atoms through one term.

Observe that the graph does not correspond to the factorizations one-to-one.
A factorization of $R$ with several atoms containing $g_R$ may be reached
along several paths.  This does not matter for length sets.  To enumerate
the factorizations, the solver therefore uses a separate depth-first search,
which removes atoms in a fixed order and never returns to an
earlier atom, so that every factorization is produced exactly once.

For other calculations and queries (e.g., finding the minimal factorization length, or for checking whether a length $\ell$ factorization exists) adequately
modified versions of this procedure are used such that superfluous computations
are avoided. The data structure used to represent remainders internally
has been chosen with care too: essentially, we use vectors of
multiplicities, packed into a single integer with one spare bit per term,
so that whether an atom divides a remainder, and the next remainder, are
decided by one subtraction.

To give an idea of the sizes, let $U = e_1^2e_2^2e_3^2(e_1+e_2+e_3)$ be the
atom of length seven of form 1(a) of Theorem~\ref{3.3}.  The sequence
$\bigl(U(-U)\bigr)^3$ of length $42$ has the length set $[6,21]$, which is
computed in a hundredth of a second from a graph with $3093$ remainders and
$16236$ edges; the same sequence has $29536$ factorizations.

\subsection{Computational results}\label{sec:software:results}

From here on, let $G = C_3 \oplus C_3 \oplus C_3$.  All computations below
use the complete catalogue of the $27911$ atoms of $\mathcal B (G^\bullet)$.  The databases
they produce, their contents and the way they are checked are described in
Subsection~\ref{sec:software:structural-certificates}.

We begin with the following lemma that characterizes all length sets of
zero-sum sequences that admit a factorization into two atoms.

\begin{lemma}\label{lem:comp:length-sets-of-pairs}
	Let $H$ be a transfer Krull monoid over $G$ and let $L \in \mathcal{L}(H)$
	be a length set. Then we have $\min L = 2$ if and only if $L$ is one of the
	following sets:
	\begin{gather*}
	\{2\},\{2,3\},\{2,4\},
 	\{2,3,4\}, \{2,3,5\}, \{2,4,5\},
 	\{2,3,4,5\}, \{2,3,4,6\}, \\
 	\{2,3,4,5,6\}, \{2,3,4,5,7\}, \{2,3,4,5,6,7\}.
	\end{gather*}
	In particular, Theorem~\ref{1.1} holds for all length sets $L$ with $\min L = 2$.
\end{lemma}

\begin{proof}
	By the transfer property (see Subsection \ref{2.4}) it suffices to consider $\mathcal B(G)$.  Let $B \in \mathcal B (G)$ with $\min \mathsf L (B)=2$, say
	$B = UV$ where $U,V\in\mathcal A(G)$ with $|U| \ge |V|$.  If  $V=0$, then $\mathsf L(B)=\{2\}$. 
By	Theorem~\ref{3.3}, there is a basis of $G$ such that $U$ has one of the given normal forms.  The generator of
	\code{length-sets-with-minimum-2.jsonl} iterates
	over all $143832$ such pairs $(U, V)$, computes $\mathsf L(UV)$ and stores
	a factorization witness for every length; precisely the eleven sets of
	the statement occur.  Explicit representatives, with a factorization of
	every length, are as follows.

	\paragraph*{$\mathsf L(B)=\{2\}$.}
	If $U= (-e_1)e_1$ and $V= (-e_3){e_3}$, then 
	\[
	B=UV =  (-e_1){e_1} (-e_3){e_3}.
	\]

	\paragraph*{$\mathsf L(B)=\{2,3\}$.}
	If $U=e_1^3$, then
	\[
	B= (-U) U =(-e_1)^3 {e_1}^3.
	\]

	\paragraph*{$\mathsf L(B)=\{2,4\}$.}
	If $U=e_1e_2e_3(-e_1-e_2-e_3)$, then
	\[
	B= (-U) U
	= \big( (-e_1){e_1} \big)  \big( (-e_2){e_2} \big)  \big( (-e_3){e_3} \big)  \big( (-e_1-e_2-e_3) (e_1+e_2+e_3) \big).
	\]

	\paragraph*{$\mathsf L(B)=\{2,3,4\}$.}
	If $U=e_1^2e_2(e_1-e_2)$, then
	\begin{align*}
	B & = (-U)U \\
	& = \bigl(e_2(e_1-e_2)(-e_1)\bigr)
	\,  \bigl((-e_2)e_1(-e_1+e_2)\bigr)\, \big( (-e_1){e_1} \big),\\
	&= \big( (-e_2){e_2} \big)  \big( (-e_1){e_1}\big)^2  \big( (-e_1+e_2)({e_1-e_2}) \big) .
	\end{align*}

	\paragraph*{$\mathsf L(B)=\{2,3,5\}$.}
	If $U=e_1^2e_2e_3(e_1-e_2-e_3)$, then
	\begin{align*}
	B & = (-U)U \\
	  &=\bigl(e_2e_3(e_1-e_2-e_3)(-e_1)\bigr)
	\, \bigl((-e_2)(-e_3)e_1(-e_1+e_2+e_3)\bigr)
	\, \big( (-e_1){e_1} \big) ,\\
	&= \big( (-e_1){e_1}\big)^2  \big( (-e_2){e_2}\big)  \big( (-e_3){e_3}\big)   \big( (-e_1+e_2+e_3)({e_1-e_2-e_3}) \big) .
	\end{align*}

	\paragraph*{$\mathsf L(B)=\{2,4,5\}$.}
	If $a=e_1+e_2+e_3$,  $U=e_1^2e_2^2e_3^2a$, and  $V=e_2^2e_3^2(-e_1-e_2)(-e_1-e_3)(-a)$, then
	\begin{align*}
	B & = UV \\
	&=e_2^3  e_3^3
	 \bigl(e_1e_2e_3(-a)\bigr)
	 \bigl(e_1a(-e_1-e_2)(-e_1-e_3)\bigr),\\
	&=e_2^3  e_3^3
	 \bigl(e_1e_2(-e_1-e_2)\bigr)
	 \bigl(e_1e_3(-e_1-e_3)\bigr)  \big( (-a)a \big).
	\end{align*}

	\paragraph*{$\mathsf L(B)=\{2,3,4,5\}$.}
	If $a=e_1+e_2$ and  $U=e_1^2e_2^2a$, then
	\begin{align*}
	B & = (-U)U \\
	& = \big( (-e_2)e_2 \big)  \bigl(e_2(-e_1)^2(-a)\bigr)
	\, \bigl((-e_2)e_1^2a\bigr),\\
	&= \big( (-e_2){e_2}\big) \bigl(e_1e_2(-a)\bigr)
	\, \bigl((-e_1)(-e_2)a\bigr)\,  \big( (-e_1){e_1} \big),\\
	&= \big( (-e_1){e_1}\big)^2\, \big( (-e_2){e_2}\big)^2\, \big( (-a){a} \big).
	\end{align*}

	\paragraph*{$\mathsf L(B)=\{2,3,4,6\}$.}
	If $a=e_1+e_2-e_3$ and  $U=e_1^2e_2^2e_3a$, then
	\begin{align*}
	B & = (-U)U \\
	&=\bigl(e_2^2e_3a(-e_1)\bigr)
	\, \bigl((-e_2)^2(-e_3)e_1(-a)\bigr) \, \big( (-e_1){e_1} \big),\\
	&=\bigl((-e_1)(-e_2)e_3a\bigr)
	\,\bigl(e_1e_2(-e_3)(-a)\bigr)\, \big( (-e_1){e_1} \big) \, \big( (-e_2){e_2} \big),\\
	&= \big( (-e_1){e_1}\big)^2\, \big( (-e_2){e_2}\big)^2\, \big( (-e_3){e_3}\big)\, \big( (-a){a}\big).
	\end{align*}

	\paragraph*{$\mathsf L(B)=\{2,3,4,5,6\}$.}
	If $a=e_2-e_3$,  $b=e_1+e_2$, and 
	$U=e_1^2e_2e_3ab$, then
	\begin{align*}
	B & = (-U)U \\
	&=\bigl((-e_2)e_3a\bigr)
	\,\bigl(e_1^2(-e_3)(-a)b\bigr)
	\,\bigl((-e_1)^2e_2(-b)\bigr),\\
	&= \big( (-e_3){e_3}\big)\,\bigl((-e_1)^2e_2(-b)\bigr)
	\, \big( (-a) {a}\big) \,\bigl(e_1^2(-e_2)b\bigr),\\
	&= \big( (-e_3){e_3}\big) \,\bigl(e_1e_2(-b)\bigr)
	\, \big( (-a){a}\big) \,\bigl((-e_1)(-e_2)b\bigr)\, \big( (-e_1){e_1}\big) ,\\
	&= \big( (-e_3){e_3}\big) \, \big( (-e_2){e_2}\big) \, \big( (-a){a}\big) \, \big( (-e_1){e_1}\big)^2\, \big( (-b){b}\big).
	\end{align*}

	\paragraph*{$\mathsf L(B)=\{2,3,4,5,7\}$.}
	If $a=e_1+e_2+e_3$ and   $U=e_1^2e_2^2e_3^2a$, then
	\begin{align*}
	B & = (-U)U \\
	&= \big( (-e_3){e_3}\big) \,\bigl((-e_1)^2(-e_2)^2e_3(-a)\bigr)
	\,\bigl(e_1^2e_2^2(-e_3)a\bigr),\\
	&= \big( (-e_3){e_3}\big) \,\bigl((-e_1)^2e_2e_3(-a)\bigr)
	\,\bigl(e_1^2(-e_2)(-e_3)a\bigr)\, \big( (-e_2){e_2}\big),\\
	&= \big( (-e_3){e_3}\big) \,\bigl(e_1e_2e_3(-a)\bigr)
	\,\bigl((-e_1)(-e_2)(-e_3)a\bigr)\, \big( (-e_2){e_2}\big) \, \big( (-e_1){e_1} \big),\\
	&= \big( (-e_3){e_3}\big)^2\, \big( (-e_2){e_2}\big)^2\, \big( (-e_1){e_1}\big)^2\, \big( (-a){a}\big).
	\end{align*}

	\paragraph*{$\mathsf L(B)=\{2,3,4,5,6,7\}$.}
	If $a=e_2+e_3$,  $b=e_1+e_3$, and 
	$U=e_1^2e_2^2e_3ab$, then
	\begin{align*}
	B & = (-U)U \\
	&= \big( (-e_3){e_3}\big) \,\bigl((-e_1)^2e_2^2a(-b)\bigr) \,\bigl(e_1^2(-e_2)^2(-a)b\bigr),\\
	&= \big( (-e_3){e_3}\big) \, \big( (-e_2){e_2}\big)
	\,\bigl(e_1^2e_2(-a)b\bigr)
	\,\bigl((-e_1)^2(-e_2)a(-b)\bigr),\\
	&= \big( (-e_3){e_3}\big) \, \big( (-e_2){e_2}\big)
	\,\bigl((-e_1)e_2(-a)b\bigr)
	\,\bigl(e_1(-e_2)a(-b)\bigr)\, \big( (-e_1){e_1}\big),\\
	&=\bigl(e_2e_3(-a)\bigr)\,\bigl((-e_2)(-e_3)a\bigr)
	\, \big( (-e_2){e_2}\big)\, \big( (-e_1){e_1}\big)^2\, \big( (-b) {b}\big),\\
	&= \big( (-e_3){e_3}\big) \, \big( (-e_2){e_2}\big)^2\, \big( (-a){a}\big) \, \big( (-e_1){e_1}\big)^2\, \big( (-b){b}\big). \qedhere
	\end{align*}
	\end{proof}

\begin{lemma} \label{lem:comp:mixed-triples}
Let $B=UVW$ with $U,V,W \in \mathcal A (G)$.  Suppose that
$4 \in \mathsf L (B)$ and that $1 \notin \Delta \big(\mathsf L (E) \big)$ for every
proper divisor $E$ of $B$.  Then
there is a basis $(e_1,e_2,e_3)$ of $G$ such that
\[
 A = e_1e_2e_3(-e_1-e_2-e_3),  \qquad B = A^3, \qquad \text{and} \qquad \mathsf L (B) = \{3,4\}.
\]
\end{lemma}

\begin{proof}
If one of the atoms, say $U$, is the atom $0$, then $4 \in \mathsf L (B)$
gives $\{2,3\} \subset \mathsf L (VW)$, and the proper divisor $VW$ has two
consecutive lengths.  Thus $U, V, W \in \mathcal A (G^{\bullet})$, and we may
assume that $|U| \ge |V| \ge |W|$ and that $U$ is one of the fourteen normal forms
of Theorem~\ref{3.3}.  Over the fourteen normal forms $U$, and up to the
automorphisms of $G$ fixing $U$, there are $201471$ such triples for which none
of $UV$, $UW$ and $VW$ has two consecutive lengths.  Of these, $135789$ satisfy $4 \in \mathsf L (B)$.
For all but one of them the generator exhibits a factorization of $B$ of
length at least three containing two atoms whose product $E$ satisfies
$1 \in \Delta \big(\mathsf L (E) \big)$; such an $E$ is a proper divisor of $B$.  The
remaining product is the displayed $A^3$, and a direct enumeration of all
its proper zero-sum divisors confirms that none of them has two consecutive
lengths.

The database \code{consecutive-triples.jsonl} stores, for each of the
$201471$ triples, a representative $(U,V,W)$ of its orbit together with the
size of the orbit and, if $4 \in \mathsf L(B)$ and $B \ne A^3$, a
factorization of $B$ and the two atoms exhibiting $E$.
\end{proof}

\begin{lemma} \label{lem:comp:six-long-atoms}
Let $B \in \mathcal B (G)$, let $z \in \mathsf Z (B)$ with
$|z| = \min \mathsf L (B)$,  and suppose that $|z|+1 \notin \mathsf L (B)$. Let
$U_1, \ldots, U_6 \in \mathcal A (G)$ with $U_1 \cdot \ldots \cdot U_6 \mid z$ and
$|U_i| \ge 4$ for all $i \in [1,6]$.  Then there are a subset $I \subset [1,6]$, 
$k \in \N$, and  $V_1, \ldots, V_k \in \mathcal A (G)$ such that
\[
 \prod_{i \in I} U_i = V_1 \cdot \ldots \cdot V_k , 
 \qquad |V_1| \le 3 , \qquad \text{and such that}
\]
\[
 k = |I|
 \qquad \text{or} \qquad
 |V_j| \le 3 \ \text{ for all } j \in [1,k] .
\]
\end{lemma}

\begin{proof}
By Lemma~\ref{lem:hereditary-gap}.1, every subfamily $(U_i)_{i \in I}$ with
$|I| \ge 2$ satisfies
\begin{equation}
 \min \mathsf L \Bigl( \prod_{i \in I} U_i \Bigr) = |I|
 \qquad \text{and} \qquad
 |I|+1 \notin \mathsf L \Bigl( \prod_{i \in I} U_i \Bigr) .\label{eq:six-long-conditions}
\end{equation}
Suppose, to the contrary, that no subfamily of $U_1, \ldots, U_6$ satisfies
the conclusion of the lemma.  Then $U_1, \ldots, U_6$ is a family of atoms of
length at least four such that every subfamily with at least two members
satisfies \eqref{eq:six-long-conditions} and no subfamily satisfies the
conclusion of the lemma; clearly, every subfamily of such a family has these
two properties as well.
All families with these two properties are enumerated by adding one atom at a
time, starting from a member of maximal length, which by
Theorem~\ref{3.3} may be assumed to be one of the eleven normal forms of
length at least four.  The database \code{length-sets-with-minimum-2.jsonl} of
Lemma~\ref{lem:comp:length-sets-of-pairs} restricts the atoms that may be
added, and for every subfamily containing the new member the conditions
\eqref{eq:six-long-conditions} and the conclusion of the lemma are tested.
Here the alternative $k=|I|$ has to be tested only when
$\sum_{i \in I}(7-|U_i|) \ge 4$: since
$\sum_{j=1}^{k}|V_j| = \sum_{i \in I}|U_i|$, $|V_1| \le 3$ and $|V_j| \le 7$
for $j \in [2,k]$, the alternative $k=|I|$ gives
$\sum_{i \in I}|U_i| \le 3 + 7(|I|-1) = 7|I|-4$.  The enumeration tests
$476\,820$ subfamilies, of which $340\,075$ produce a factorization as in
the lemma.  It finds $1\,270$ families with the two properties: $992$ with
two members, $158$ with three, $84$ with four, $36$ with five, and none with
six.  This contradiction proves the lemma.
The database \code{six-long-atoms.jsonl} contains the $1\,270$ families;
each record stores the family and its complete length set.
\end{proof}

\subsection{The doubled four-cap}\label{sec:software:doubled-cap}

Let $(e_1,e_2,e_3)$ be a basis of $G$,  let  $e_0=e_1+e_2+e_3$, and write
\[
 G_0=\{e_0,e_1,e_2,e_3\}\qquad\text{and}\qquad
 D=\prod_{i=0}^{3}\bigl((-e_i)e_i\bigr)^2 .
\]
Thus $G_0 \cup (-G_0) = \{ \pm e_0, \pm e_1, \pm e_2, \pm e_3 \}$ is the
doubled four-cap of Subsection~\ref{2.5}.
For $a \in [5,8]$ we write $D_a$ for the subproduct of $D$ consisting of $a$
of its eight inverse pairs, at least one for each index, that is
\[
 \begin{aligned}
 D_5&=\bigl((-e_0)e_0\bigr)\bigl((-e_1)e_1\bigr)\bigl((-e_2)e_2\bigr)\bigl((-e_3)e_3\bigr)^2, &
 D_6&=\bigl((-e_0)e_0\bigr)\bigl((-e_1)e_1\bigr)\bigl((-e_2)e_2\bigr)^2\bigl((-e_3)e_3\bigr)^2, \\
 D_7&=\bigl((-e_0)e_0\bigr)\bigl((-e_1)e_1\bigr)^2\bigl((-e_2)e_2\bigr)^2\bigl((-e_3)e_3\bigr)^2, &
 D_8&=D .
 \end{aligned}
\]

\begin{lemma} \label{lem:comp:doubled-cap-length-three}
Let $T\in\mathcal A(G)$ with $|T|=3$.  If
$8\notin\mathsf L(DT)$, 
then $T=(\varepsilon e_i)^3$ for some $i\in[0,3]$ and some
$\varepsilon\in\{1,-1\}$.
\end{lemma}

\begin{proof}
The stabilizer of the cap $G_0 \cup (-G_0)$ has order $48$, it fixes $D$,
and it has ten orbits on the $130$ atoms of length three, of sizes $8$,
$12$, $24$, $12$, $24$, $6$, $12$, $12$, $12$ and $8$.  The orbit of size
$8$ listed first is $\{(\varepsilon e_i)^3:i\in[0,3], \varepsilon\in\{1,-1\}\}$.
Since the condition $8 \in \mathsf L (DT)$ is invariant under the
stabilizer, it suffices to exhibit, for one representative $T$ of each of
the other nine orbits, a factorization of $DT$ of length eight.
The ten orbits, with one representative $T$ each, a length-8 
factorization of $DT$ for the nine orbits other than the cubes, and the
complete length sets $\mathsf{L}(DT)$, are stored in
the database file \code{doubled-cap-length-three.jsonl}.
A direct comparison completes the proof.
\end{proof}

We shall also use the following fixed products.  They are included in
\code{fixed-length-sets.jsonl}; the database stores a witness for every displayed length.

\begin{lemma} \label{lem:comp:doubled-cap-seeds}
With the notation above, we obtain the following length sets.
\begin{align*}
 \mathsf L(D_5)&=[2,3]\cup\{5\},\\
 \mathsf L(D_6)&=[2,4]\cup\{6\},\\
 \mathsf L(D_7)&=[2,5]\cup\{7\},\\
 \mathsf L(D_8)=\mathsf L(D)&=[3,6]\cup\{8\},\\
 \mathsf L(D e_1^3e_2^3)&=[4,8]\cup\{10\},\\
 \mathsf L(D e_1^3e_2^3e_3^3)&=[4,9]\cup\{11\},\\
 \mathsf L(D e_0^3(-e_1)^3e_2^3e_3^3)&=[4,10]\cup\{12\},\\
 \mathsf L\bigl(D (e_0^3(-e_1)^3e_2^3e_3^3)^2\bigr)
   &=[6,14]\cup\{16\}.
\end{align*}
\end{lemma}

\smallskip
In the following lemma we prove that length sets of zero-sum sequences having a factorization of minimal length
that is not divisible by $A^3$ for any atom $A$ of length 7,
can only attain one particular form. Moreover,  up to automorphisms of $G$
there is only one such zero-sum sequence.

\begin{lemma} \label{lem:comp:doubled-cap-cubes}
Let $B \in \mathcal B (G)$ and  $k \in \N$  with
\[
 \min \mathsf L (B) = 3k,
 \qquad
 \max \mathsf L (B) = 7k+2,
 \qquad \text{and} \qquad
 7k+1 \notin \mathsf L (B) .
\]
Suppose $x \in \mathsf Z (B)$ with $|x| = 3k$ such that $A^3\nmid x$ for any $A\in \mathcal A(G)$ of length $7$.  Then $k=2$, and there is a basis
$(e_1,e_2,e_3)$ of $G$ such that, with $e_0=e_1+e_2+e_3$ and
$D=\prod_{i=0}^{3}\bigl((-e_i)e_i\bigr)^2$ as in
Subsection~\ref{sec:software:doubled-cap},
\[
 B = D\bigl(e_0^3(-e_1)^3e_2^3e_3^3\bigr)^2 
 \qquad \text{and} \qquad
 \mathsf L (B) = [6,14]\cup\{16\} .
\]
\end{lemma}

\begin{proof}
Throughout, let $z \in \mathsf Z(B)$ with $|z| = 7k+2 = \max \mathsf L (B)$.
By  Lemma~\ref{lem:hereditary-gap}.2, every divisor $y$ of $z$ satisfies
\begin{equation}
 |y| = \max \mathsf L \big(\pi (y) \big)
 \qquad \text{and} \qquad
 |y|-1 \notin \mathsf L \big(\pi (y) \big) .\label{eq:doubled-cap-cubes-local}
\end{equation}

\smallskip
\noindent
\emph{Step 1: the atoms of length two.}  The atoms of length two occurring
in $z$ have the form $(-g)g$ with $g \in G^{\bullet}$.  As in the proof of
Lemma~\ref{lem:upper-hole-affine-cap}, none of them occurs three times, and
the union $C_0$ of the sets $\{g,-g\}$ with $(-g)g \mid z$ is
a symmetric cap.  By Lemma~\ref{lem:symmetric-caps}, there is a basis
$(e_1,e_2,e_3)$ of $G$ such that, with $e_0 = e_1+e_2+e_3$ and
$G_0 = \{e_0,e_1,e_2,e_3\}$, we have $C_0 \subset G_0 \cup (-G_0)$; thus the
product of the atoms of length two occurring in $z$, counted with
multiplicity, is the divisor
\[
 y_{\lambda} = \prod_{i=0}^{3}\bigl((-e_i)e_i\bigr)^{\lambda_i} \in \mathsf Z (G)
 \qquad \text{of $z$, with} \qquad \lambda_0, \ldots, \lambda_3 \in [0,2] ,
\]
and their number is $a = \lambda_0+\ldots+\lambda_3$.  Every permutation of
the four sets $\{\pm e_0\}, \ldots, \{\pm e_3\}$ is induced by an automorphism
of $G$ fixing $G_0 \cup (-G_0)$: with $e_4 = -e_0$ we have
$e_1+e_2+e_3+e_4 = 0$, so for every permutation $\tau$ of $[1,4]$ the
automorphism with $e_i \mapsto e_{\tau(i)}$ for $i \in [1,3]$ maps $e_4$ to
$e_{\tau(4)}$, and $\{\pm e_4\} = \{\pm e_0\}$.

\smallskip
\noindent
\emph{Step 2: counting terms.}  Put
$\delta = \sum_{W}\bigl(|W|-3\bigr)$, the sum being taken over the atoms
$W$ of length at least three occurring in $z$, counted with multiplicity.  Counting terms in $z$ gives
$|B| = 2a+3\bigl(|z|-a\bigr)+\delta = 21k+6-a+\delta$, while counting terms
in $x$ gives $|B| \le 7|x| = 21k$ by $\mathsf D (G) = 7$.  Hence
\begin{equation}
 7|x|-|B| = a-6-\delta \ \ge \ 0 ,\label{eq:doubled-cap-cubes-defect}
\end{equation}
so that $a \ge 6$ and $\delta \le a-6$; and $a \le 8$ by
Lemma~\ref{lem:upper-hole-affine-cap}.  By Step 1 we may therefore assume
that
\[
 (\lambda_1,\lambda_2,\lambda_3,\lambda_0) \in
 \bigl\{ (2,2,2,0),\ (2,2,1,1),\ (2,2,2,1),\ (2,2,2,2) \bigr\} ,
\]
these being all vectors in $[0,2]^4$ with sum in $[6,8]$ up to the order of
their entries.

\smallskip
\noindent
\emph{Step 3: the atoms of length at least three.}  Fix one of the four
vectors $\lambda$.  Let $\mathcal T_{\lambda}$ be the set of atoms $T$ of
length three such that $y_{\lambda}T \in \mathsf Z (G)$ satisfies
\eqref{eq:doubled-cap-cubes-local} in place of $y$, and let $\Gamma_{\lambda}$ be the graph with vertex set
$\mathcal T_{\lambda}$ in which two distinct vertices $T, T'$ are adjacent if
$y_{\lambda}TT'$ satisfies \eqref{eq:doubled-cap-cubes-local}.  By
\eqref{eq:doubled-cap-cubes-local}, applied to the divisors $y_\lambda T$ and
$y_\lambda TT'$ of $z$, the distinct atoms of length three occurring in $z$ form a
clique of $\Gamma_{\lambda}$.  The computation shows that every maximal
clique of $\Gamma_\lambda$ contains at least one cube $g^3$ and at most one
atom of length three that is not of the form $g^3$ (a \emph{non-cube}); a
non-cube $T = g_1g_2g_3$ occurs at most twice in $z$,
since otherwise replacing $T^3$ by $g_1^3g_2^3g_3^3$ gives a factorization of $B$ of
length $|z|$, to which \eqref{eq:doubled-cap-cubes-local} applies as well, whereas the
computation shows that for each of the (two) non-cubes in the sets
$\mathcal T_\lambda$ one of the cubes $g_1^3, g_2^3, g_3^3$ does not lie in
$\mathcal T_\lambda$.

If $a = 6$, then $\delta = 0$ by \eqref{eq:doubled-cap-cubes-defect}.  If $a \in
\{7,8\}$, then $\delta \le 2$, so every atom occurring in $z$ has length at most five;
the computation shows that no atom $W$ of length four (if $a = 7$) and no
atom of length four or five (if $a = 8$) is such that $y_\lambda W$ satisfies
\eqref{eq:doubled-cap-cubes-local}, whence $\delta = 0$ as well.  Thus every atom occurring in
$z$ has length two or three, and \eqref{eq:doubled-cap-cubes-defect} becomes
$7|x|-|B| = a-6$: all atoms occurring in $x$ have length seven, except for one of
length six if $a = 7$, and except for one of length five or two of length six
if $a = 8$.

\smallskip
\noindent
\emph{Step 4: the enumeration.}  By Step 3, every atom occurring in $z$ has length
two or three, and the atoms of length three occurring in $z$ lie in a maximal clique
$\mathcal K$ of $\Gamma_\lambda$.  Let $G_{\mathcal K}$ be the set of all
$g \in G$ with $g^3 \in \mathcal K$, and let $T$ be the non-cube in
$\mathcal K$ (if there is one).  Then
\[
 B = \pi(y_\lambda)\, T^j \prod_{g \in G_{\mathcal K}} g^{3n_g}
 \qquad \text{with } j \in [0,2] \text{ and } n_g \in \N_0 ,
\]
so every term of $B$ outside $G_{\mathcal K}$ is a term of
$\pi(y_\lambda)T^j$.  Consider the atoms occurring in $x$.  Those with a term outside
$G_{\mathcal K}$ together, counted with multiplicity, contain exactly the
terms of $\pi(y_\lambda)T^j$ outside $G_{\mathcal K}$, which leaves only
finitely many possibilities for them.  Every other atom of length seven
occurring in $x$ has support in $G_{\mathcal K}$, and by hypothesis it occurs at most twice in $x$.  For
every $\lambda$, every maximal clique $\mathcal K$ of $\Gamma_\lambda$ and
every $j$, the generator enumerates all factorizations $x$ of this kind,
with the atoms of length less than seven allowed by Step 3.  It
determines the exponents $n_g$ from the multiplicities of the terms in
$G_{\mathcal K}$, keeps those with $|x| = 3k$ and $|z| = 7k+2$, and computes
$\mathsf L (B)$ for each resulting sequence $B$.
The results are as follows; the last two columns give the number of
sequences $B$ enumerated and the number of those satisfying the three
hypotheses on $\mathsf L (B)$.
\[
\begin{array}{c|c|c|c|c}
(\lambda_1,\lambda_2,\lambda_3,\lambda_0) & |\mathcal T_\lambda| &
 \text{maximal cliques of } \Gamma_\lambda & \text{sequences } B & \text{hypotheses hold}\\ \hline
(2,2,2,0)&14&64\text{ of size }4&8&0\\
(2,2,1,1)&10&8\text{ of size }4,\
             8\text{ of size }5&2&0\\
(2,2,2,1)&8&16\text{ of size }4&2&0\\
(2,2,2,2)&8&16\text{ of size }4&8&6
\end{array}
\]
The six sequences in the last column have
$(\lambda_1,\lambda_2,\lambda_3,\lambda_0)=(2,2,2,2)$ and $k=2$, and they
form a single orbit of size six under the automorphisms of $G$ fixing
$G_0 \cup (-G_0)$; one of them is $D\bigl(e_0^3(-e_1)^3e_2^3e_3^3\bigr)^2$,
whose length set is $[6,14]\cup\{16\}$ by
Lemma~\ref{lem:comp:doubled-cap-seeds}.

The database \code{doubled-cap-cubes.jsonl} contains the four rows of
the table and the twenty sequences $B$, each with its length set and
factorizations of minimal and maximal length.
\end{proof}

\subsection{On 2-4-sequences and their extensions}\label{sec:software:carriers}

The finite calculations of this subsection are the ones used in the proof of
the structural results in Section~\ref{5}.  A zero-sum sequence
$D \in \mathcal B (G)$ is called a \emph{$2$-$4$-sequence} if $\{2,4\} \subset \mathsf L (D)$ and 
$3 \notin \mathsf L (D)$.  
We consider pairs $(D,x)$, for $2$-$4$-sequences $D$ and for factorizations $x\in \mathsf Z(G)$,  satisfying the following properties:
\begin{equation}\label{eq:carrier-extension-state}
 \begin{cases}
 \begin{aligned}\hspace{1ex}
  &\text{(P1)} && 1 \in \Delta \bigl(\mathsf L (D\pi(x))\bigr) ,\\
  &\text{(P2)} && 1 \notin \Delta \bigl(\mathsf L (D\pi(x'))\bigr) \ \text{ for every proper divisor $x'$ of $x$, and}\\
  &\text{(P3)} && |x|+3 \notin \mathsf L \big(D\pi(x) \big) .
 \end{aligned}
\end{cases}
\end{equation}
We call such a pair a \emph{critical pair}, and for a fixed $2$-$4$-sequence
$D$ we say that $x$ is \emph{critical} if $(D,x)$ is a critical pair.  The
name refers to the point at which two consecutive lengths first appear:
multiplying $D$ by $\pi(x)$ creates two consecutive lengths by (P1), while
multiplying $D$ by $\pi(x')$ for a proper divisor $x'$ of $x$ does not by (P2).
Since $yx$ is a factorization of $D\pi(x)$ of length $|x|+2$ for every
$y \in \mathsf Z_2 (D)$, condition (P3) says that the length directly above
it is missing.  Recall that $\Aut (G)$ acts on $\mathcal B (G)$ and on
$\mathsf Z (G)$; we identify two critical pairs if they differ by an
automorphism of $G$.  For a fixed $D$, we call $x \in \mathsf Z (G)$
\emph{pure} if $1 \notin \Delta \big( \mathsf L ( D\pi(x) )\big)$.  Thus $x$ is
critical if and only if $x$ is not pure, every proper divisor of $x$ is pure,
and $|x|+3 \notin \mathsf L \big( D\pi(x) \big)$.  We note:
\begin{itemize}
\item Every divisor $x'$ of a pure element $x$ is pure: multiplying
      factorizations of $D\pi(x')$ of lengths $\ell$ and $\ell+1$ by
      $xx'^{-1}$ gives factorizations of $D\pi(x)$ of lengths
      $\ell + |xx'^{-1}|$ and $\ell+1+|xx'^{-1}|$.
\item A pure $x$ satisfies $|x|+3 \notin \mathsf L \big( D\pi(x) \big)$, since
      $|x|+2 \in \mathsf L \big( D\pi(x) \big)$ and all elements of
      $\mathsf L \big( D\pi(x) \big)$ have the same parity.
\item No critical $x$ is divisible by the atom $0$: if $x = 0x''$, then
      $\mathsf L \big( D\pi(x) \big) = 1 + \mathsf L \big( D\pi(x'') \big)$ by
      \eqref{shift-by-zeros}, so that $x''$ would not be pure.  Hence no
      element of $\mathsf Z (G)$ considered below is divisible by the atom $0$.
\end{itemize}
The critical pairs are obtained by an explicit enumeration up to $|x| = 5$ and by a
parity argument beyond.

\begin{lemma} \label{lem:comp:carrier-antichain}
Up to automorphisms of $G$, there are exactly $52$ different $2$-$4$-sequences $D$, one
with $\mathsf L (D) = \{2,4,5\}$ and $51$ with $\mathsf L (D) = \{2,4\}$,
and there are exactly $652$ critical pairs $(D,x)$ with $|x| \le 14$:
$(D,1)$ for the $2$-$4$-sequence with $\mathsf L (D) = \{2,4,5\}$, and $651$
with $|x| \in [1,4]$.
\end{lemma}

\begin{proof}
A $2$-$4$-sequence is a product of two atoms.  Hence the $52$
$2$-$4$-sequences, their length sets and their factorizations of length two
are read off the atom catalogue and the database of
Lemma~\ref{lem:comp:length-sets-of-pairs},
\code{length-sets-with-minimum-2.jsonl}; their length sets are $\{2,4\}$ and
$\{2,4,5\}$.  If $\mathsf L (D) = \{2,4,5\}$, then $(D,1)$ is
a critical pair, and it is the only one with this $D$, since $1$ is a proper divisor of
every $x \ne 1$.  Now fix one of the $51$ $2$-$4$-sequences $D$ with
$\mathsf L (D) = \{2,4\}$; then every critical $x$ satisfies $x \ne 1$, and all
enumerations below are carried out up to the stabilizer of $D$ in $\Aut (G)$,
which maps critical elements to critical elements and pure elements to pure
elements.

Step~1 determines all critical $x$ with $|x| \le 5$.  Steps~2 to~4 show that
there is no critical $x$ with $|x| \in [6, 14]$: Step~2 shows that such an $x$ is built from a
divisor of length three and the atoms of a maximal clique of a finite graph,
Step~3 excludes all but $37$ of these cliques by the parity criterion of
Lemma~\ref{lem:parity-criterion}, and Step~4 excludes the remaining $37$ by
enumeration.  The counts are collected in the table at the end of the proof.

\emph{Step 1: the critical $x$ with $|x| \le 5$.}  For $k \in [0,5]$, let
$P_k$ be the set of all pure $x \in \mathsf Z (G)$ with $|x| = k$; thus
$P_0 = \{1\}$, since $\mathsf L (D) = \{2,4\}$.  Every divisor of an element of
$P_k$ is pure, and so is every proper divisor of a critical $x$.  Hence, for
$k \in [1,5]$, all elements of $P_k$ and all critical $x$ with $|x| = k$ are
among the $x$ of length $k$ all of whose divisors of length $k-1$ lie in
$P_{k-1}$.  For $k \in [1, 5]$ in turn, we form these $x$ and compute
$\mathsf L \big( D\pi(x) \big)$ for one representative of each orbit.  This determines
$P_1, \ldots, P_5$ and all critical $x$ with $|x| \le 5$: there are $651$ of
them, all with $|x| \le 4$.

\emph{Step 2: prefixes and cliques.}  Let $x$ be critical with
$6 \le |x| \le 14$, and let $p \mid x$ with $|p| = 3$.  Since all proper divisors of $x$ are
pure, $p \in P_3$, $pC \in P_4$ for every atom $C \mid xp^{-1}$, and
$pCC' \in P_5$ for all atoms $C, C'$ with $CC' \mid xp^{-1}$ (where $C = C'$ is
allowed).  After applying an element of the stabilizer of $D$, we may assume
that $p$ is one of a fixed set of representatives of the orbits of $P_3$.
Consider the graph whose vertices are the atoms $C$ with $pC \in P_4$, in
which two distinct vertices $C$ and $C'$ are adjacent if $pCC' \in P_5$, and
in which $C$ has a loop if $pC^2 \in P_5$.  Then $\supp (xp^{-1})$ is
contained in a maximal clique $\mathcal K$ of this graph, and
\begin{enumerate}
\item[(a)] $\mathsf v_C (xp^{-1}) \le 1$ for every $C \in \mathcal K$ without a
      loop, and
\item[(b)] $\mathsf v_C (x) \le 2$ for every $C \in \mathcal K$ with
      $1 \in \Delta \big(\mathsf L (C^3)\big)$: if $C^3 \mid x$, then, since
      $|x| \ge 6$, there is an atom $C'$ with $C' \mid x$ and $C^3 \mid xC'^{-1}$, and
      the proper divisor $xC'^{-1}$ of $x$ would not be pure, since its
      divisor $C^3$ is not.
\end{enumerate}
(Whether $1 \in \Delta \big(\mathsf L (C^3)\big)$ is read off the length sets
$\mathsf L (W^3)$ of the fourteen normal forms $W$ of Theorem~\ref{3.3}.)  Let
$\mathcal K_{\infty}$ be the set of atoms of $\mathcal K$ to which neither (a)
nor (b) applies, that is, the atoms with a loop and
$1 \notin \Delta \big(\mathsf L (C^3)\big)$, let
$\mathcal K_b = \mathcal K \setminus \mathcal K_{\infty}$, and, for
$C \in \mathcal K_b$, let $c(C)$ be the smallest of the bounds $1$ (if (a)
applies) and $2 - \mathsf v_C (p)$ (if (b) applies) on
$\mathsf v_C (xp^{-1})$.  Hence
\[
 x = p \prod_{C \in \mathcal K_{b}} C^{n_C} \prod_{C \in \mathcal K_{\infty}} C^{m_C}
 \qquad \text{with} \quad 0 \le n_C \le c(C) \quad \text{and} \quad m_C \in \N_0 ,
\]
and it remains to show that no element of this form with
$6 \le |x| \le 14$ is critical.

\emph{Step 3: the parity argument.}  Let $G_0$ be the set of terms of the atoms
of $\mathcal K_{\infty}$, and suppose that every atom of $\mathcal B (G_0)$ has
odd length; this holds for all but $37$ of the maximal cliques.  Fix the
exponents $(n_C)_{C \in \mathcal K_b}$ and apply
Lemma~\ref{lem:parity-criterion} with
$A = D\pi(p)\prod_{C \in \mathcal K_{b}} C^{n_C}$ and with
$C_1, \ldots, C_r$ the atoms of $\mathcal K_{\infty}$, so that
$B_m = D\pi(x)$ for
$x = p\prod_{C \in \mathcal K_b} C^{n_C}\, C_1^{m_1} \cdot \ldots \cdot C_r^{m_r}$.
Let $M$ be the set of all $m \in \N_0^r$ with
$1 \in \Delta \big(\mathsf L (B_m)\big)$.  For $y, y' \in \mathsf Y$ with
$|y| + |\pi(y)| \not\equiv |y'| + |\pi(y')| \mod 2$, let $M_{y,y'}$ be the
set of all $m \in \N_0^r$ satisfying \eqref{eq:cone-feasibility}.  By
Lemma~\ref{lem:parity-criterion}, $M$ is the union of the finitely many sets
$M_{y,y'}$; in particular, $M = \emptyset$ if all residues
$|y| + |\pi(y)| \bmod 2$, $y \in \mathsf Y$, agree.

The left-hand side of \eqref{eq:cone-feasibility} increases with $m$, so each
$M_{y,y'}$ is closed upwards with respect to the componentwise order.
Moreover, a minimal element of $M_{y,y'}$ has all its coordinates at most the
largest right-hand side in \eqref{eq:cone-feasibility}, so $M_{y,y'}$ has only
finitely many minimal elements.  Hence the minimal elements of $M$ are the
minimal ones among the minimal elements of the sets $M_{y,y'}$, and there are
finitely many.

If $x$ is critical with exponent vector $m$, then $m \in M$ by (P1).  Every
$m' \ne m$ with $m' \le m$ componentwise gives a proper divisor of $x$, which
is pure, so that $m' \notin M$.  Thus $m$ is a minimal element of $M$.  For
every minimal element $m$ of $M$, we compute $\mathsf L \big( D\pi(x) \big)$ for the
corresponding $x$.  It has two consecutive elements, as predicted, and either
$|x|+3 \in \mathsf L \big( D\pi(x) \big)$, or there is an atom $C \mid x$ such that
$xC^{-1}$ is not pure.  In both cases $x$ is not critical.

\emph{Step 4: the remaining cliques.}  For the $37$ maximal cliques not
covered by Step 3 we enumerate all
$x = p \prod_{C \in \mathcal K} C^{n_C}$ with $n_C \le c(C)$ for
$C \in \mathcal K_b$ and $|x| \le 14$, by increasing length, adding one atom
at a time and extending only pure elements; this is exhaustive, because all
proper divisors of a critical $x$ are pure.  For every element with $|x| \ge 6$ which
is not pure and has $|x|+3 \notin \mathsf L \big( D\pi(x) \big)$, there is an atom
$C \mid x$ such that $xC^{-1}$ is not pure; hence none of them is critical.

The counts, summed over the $51$ $2$-$4$-sequences, are as follows.
\[
\begin{array}{l|l}
 \text{Step 1} & 651 \text{ critical pairs } (150, 398, 100, 3 \text{ with } |x| = 1, 2, 3, 4)\\
 \text{Step 2} & 17796 \text{ prefixes } p,\ 154710 \text{ maximal cliques}\\
 \text{Step 3} & 154673 \text{ cliques},\ 1143508 \text{ choices of } (n_C),\ 849955 \text{ of them with } M = \emptyset,\\
               & 724335 \text{ minimal elements of } M, \text{ giving } 229097 \text{ elements } x \text{ with } |x| \le 10:\\
               & 209423 \text{ with } |x|+3 \in \mathsf L \big( D\pi(x) \big),\ 19674 \text{ with } xC^{-1} \text{ not pure}\\
 \text{Step 4} & 37 \text{ cliques},\ 91275 \text{ elements},\ 12540 \text{ of them with } |x| \ge 6 \text{ not pure}
\end{array}
\]
The database \code{critical-pair-layers.jsonl} contains the $2$-$4$-sequences, the
layers $P_1, \ldots, P_5$ and the critical pairs of Step~1 with their length sets
$\mathsf L \big( D\pi(x) \big)$, and
\code{critical-pair-cliques.jsonl} contains, for every prefix $p$, the graph and its
maximal cliques, and, for every clique, the data of Step~3 or Step~4.
\end{proof}

\begin{lemma} \label{lem:comp:carrier-minimum-side}
Let $(D,x)$ be one of the critical pairs in Lemma~\ref{lem:comp:carrier-antichain}.
Suppose that
\[
 |x|+2 = \min \mathsf L \big( D \pi (x) \big)
\]
and that there are $U,V \in \mathcal A (G)$ such that $D=UV$, 
\[
 1 \notin \Delta \big( \mathsf L (\pi (Ux)) \big),
 \quad \text{and} \quad
 1 \notin \Delta \big(\mathsf L (\pi (Vx)) \big) .
\]
Then one of the following holds.
\begin{enumerate}
\item $|x| \le 2$.  Up to automorphisms of $G$, this occurs for exactly
      one critical pair with $|x|=0$, namely $(D,1)$ with
      $\mathsf L (D) = \{2,4,5\}$, for exactly $127$ with $|x|=1$,
      and for exactly $65$ with $|x|=2$.
\item $|x|=3$, and there are both a factorization $w$ of $D\pi(x)$ with $|w|=5$
      and some $T \in \mathcal A (G)$ with  $T \mid w$,  $|T| \le 3$ such that
      $1 \in \Delta \big( \mathsf L (\pi (T^{-1}w)) \big)$.  Up to automorphisms of
      $G$, this occurs for exactly nine critical pairs.
\end{enumerate}
In particular, $|x| \le 3$.
\end{lemma}

\begin{proof}
The critical pair $(D,1)$ with $\mathsf L (D) = \{2,4,5\}$ satisfies the
hypotheses, since $\min \mathsf L (D) = 2$ and
$\mathsf L (U) = \mathsf L (V) = \{1\}$.  For the other $651$ critical pairs,
$\mathsf L (D) = \{2,4\}$, and their length sets $\mathsf L \big( D\pi(x) \big)$ are
known from Lemma~\ref{lem:comp:carrier-antichain}.  For $216$ of them,
$\min \mathsf L \big( D\pi(x) \big) \ne |x|+2$.  For each of the remaining $435$, we
compute $\mathsf L \big( \pi(Ux) \big)$ and $\mathsf L \big( \pi(Vx) \big)$ for every
$UV \in \mathsf Z_2 (D)$ and find that for $234$ of them no $UV$ satisfies the
hypotheses. This leaves $127$ critical pairs with $|x| = 1$, $65$ with
$|x| = 2$, nine with $|x| = 3$, and none with $|x| = 4$, which proves Part~1
and the final assertion.  For each of the nine with $|x| = 3$, enumerating the
factorizations of $D\pi(x)$ of length five yields $w$ and $T$ as in Part~2.

The database \code{bounded-endpoint.jsonl} contains, for each of these $202$
critical pairs, the factorizations $UV \in \mathsf Z_2 (D)$ with the length
sets $\mathsf L \big( \pi(Ux) \big)$ and $\mathsf L \big( \pi(Vx) \big)$, and, for the nine with
$|x| = 3$, the factorization $w$ and the atom $T$.
\end{proof}

The next calculation extends the $193$ critical pairs
of Lemma~\ref{lem:comp:carrier-minimum-side}.1 by further atoms of length at
least four.

\begin{lemma} \label{lem:comp:bounded-endpoint}
Let $(D,x)$ be one of the $193$ critical pairs of
Lemma~\ref{lem:comp:carrier-minimum-side}.1, let $t \in [0,5]$, let
$C_1, \ldots, C_t \in \mathcal A (G)$ with $|C_i| \ge 4$ for all
$i \in [1,t]$, and put $B = D\pi(x)C_1 \cdot \ldots \cdot C_t$.  If
\[
 \min \mathsf L (B) = |x|+t+2
 \qquad \text{and} \qquad
 |x|+t+3 \notin \mathsf L (B) ,
\]
then $\max \mathsf L (B) - 1 \in \mathsf L (B)$.
\end{lemma}

\begin{proof}
Put $F = D\pi(x)$.  For $t \in [0,5]$, let $\mathcal Q_t$ be the set of all
$c = C_1 \cdot \ldots \cdot C_t \in \mathsf Z (G)$ with $|C_i| \ge 4$ for all
$i \in [1,t]$ such that
\[
 \min \mathsf L \big( F\pi(c) \big) = |x|+t+2
 \qquad \text{and} \qquad
 |x|+t+3 \notin \mathsf L \big( F\pi(c) \big) .
\]
Then $\mathcal Q_0 = \{1\}$ by Lemma~\ref{lem:comp:carrier-minimum-side} and
(P3), and we have to show that
$\max \mathsf L \big( F\pi(c) \big) - 1 \in \mathsf L \big( F\pi(c) \big)$ for every
$c \in \mathcal Q_t$.  Let $\mathcal E_D$ be the set of all atoms $C$ with
\[
 |C| \ge 4, \qquad \min\mathsf L(DC)=3, \qquad 4\notin\mathsf L(DC).
\]

Let $t \ge 1$, $c = C_1 \cdot \ldots \cdot C_t \in \mathcal Q_t$,
and $y \in \mathsf Z_2 (D)$.  Then $yxc$ is a factorization of $F\pi(c)$ of
minimal length $|x|+t+2$, and $|x|+t+3$ is not a length, whence
Lemma~\ref{lem:hereditary-gap}.1 applies to every divisor of $yxc$.  For the
divisors $yC_i$ it gives $C_i \in \mathcal E_D$, and for the divisors
$yxcC_j^{-1}$ it gives $cC_j^{-1} \in \mathcal Q_{t-1}$.  Hence every element
of $\mathcal Q_t$ is a product of $t$ atoms of $\mathcal E_D$ all of whose
divisors of length $t-1$ lie in $\mathcal Q_{t-1}$.

The automorphisms of $G$ fixing $D$ and $x$ map
$\mathcal E_D$ and each $\mathcal Q_t$ to themselves.  For $t \in [1,  5]$
in turn, we form all products $c$ of $t$ atoms of $\mathcal E_D$ all of whose
divisors of length $t-1$ lie in $\mathcal Q_{t-1}$, and we decide for one
representative of each orbit whether it lies in $\mathcal Q_t$.  This
determines $\mathcal Q_1, \ldots, \mathcal Q_5$.  For one representative $c$ of
each orbit of $\mathcal Q_t$, $t \in [0,5]$, we compute
$\mathsf L \big( F\pi(c) \big)$ and find that
$\max \mathsf L \big( F\pi(c) \big) - 1 \in \mathsf L \big( F\pi(c) \big)$, which proves the
lemma.

The $34$ $2$-$4$-sequences occurring in the $193$ critical pairs have $689$
atoms in their sets $\mathcal E_D$ in total.  Summed over the $193$ critical
pairs, the numbers of products $c$ formed in the enumeration are as follows.
\[
\begin{array}{l|ccccc}
 t & 1 & 2 & 3 & 4 & 5\\ \hline
 \text{elements of } \mathcal Q_t & 1540 & 5128 & 11328 & 20263 & 32001\\
 \text{other products formed} & 3702 & 4850 & 57 & 21 & 2
\end{array}
\]
The elements of $\mathcal Q_1, \ldots, \mathcal Q_5$ form $57050$ orbits under
the automorphisms of $G$ fixing $D$ and $x$, and each of the $8632$ other
products $c$ is excluded by a factorization of $F\pi(c)$ of length $|x|+t+3$.
Over all $193$ critical pairs and $t \in [0,5]$, the sequences $F\pi(c)$ with
$c \in \mathcal Q_t$ are $1354$ distinct sequences up to automorphisms of $G$.

The database \code{bounded-endpoint.jsonl} contains the sets $\mathcal E_D$
and, for each critical pair and each $t$, one representative of every orbit of
the products formed, with the size of its orbit and either a factorization
showing that it does not lie in $\mathcal Q_t$, or the length set
$\mathsf L \big( F\pi(c) \big)$ with factorizations of lengths
$\max \mathsf L \big( F\pi(c) \big)$ and $\max \mathsf L \big( F\pi(c) \big) - 1$.
\end{proof}

The last calculation of this subsection concerns the products
$D\pi(x) C_1 \cdot \ldots \cdot C_e$ for one of the $652$
critical pairs $(D,x)$ of Lemma~\ref{lem:comp:carrier-antichain}, where the
adjoined atoms $C_1, \ldots, C_e$ are arbitrary, subject only to the
condition that the gap in the length set  $\mathsf L \big( D\pi(x) \big)$ persists.  It states that the
missing length then stays next to the endpoint on which it started.

\begin{lemma} \label{lem:comp:extensions}
Let $(D,x)$ be one of the $652$ critical pairs of
Lemma~\ref{lem:comp:carrier-antichain}.  Then
$\max \mathsf L \big( D\pi(x) \big) = |x|+4$ or $\min \mathsf L \big( D\pi(x) \big) = |x|+2$.
Moreover, let $e \in \N_0$ and let $C_1, \ldots, C_e \in \mathcal A (G)$ be
such that
\[
 |x|+e+3 \notin \mathsf L (D\pi(x) C_1 \cdot \ldots \cdot C_e) .
\]
\begin{enumerate}
\item If $\max \mathsf L \big( D\pi(x) \big) = |x|+4$, then
      $\max \mathsf L \big(D\pi(x) C_1 \cdot \ldots \cdot C_e \big) = |x|+e+4$.
      
\item If $\min \mathsf L \big( D\pi(x) \big) = |x|+2$, then
      $\min \mathsf L \big( D\pi(x) C_1 \cdot \ldots \cdot C_e \big) = |x|+e+2$.
\end{enumerate}
\end{lemma}

\begin{proof}
The length sets $\mathsf L \big( D\pi(x) \big)$ of the $651$ critical pairs with
$x \ne 1$ are known from Lemma~\ref{lem:comp:carrier-antichain}, and
$\mathsf L (D) = \{2,4,5\}$ for the remaining one.  Exactly $216$ of the $652$
critical pairs satisfy $\max \mathsf L \big( D\pi(x) \big)  = |x|+4$ (case 1), and the
other $436$ satisfy $\min \mathsf L \big( D\pi(x) \big) = |x|+2$ (case 2).  This proves
the first assertion.

Now put $B = D\pi(x) C_1 \cdot \ldots \cdot C_e$.  The atom $0$ divides neither
$D$ (since $\mathsf L (0V) = \{2\}$ for every $V \in \mathcal A (G)$) nor $x$
(see the remarks before Lemma~\ref{lem:comp:carrier-antichain}), and by
\eqref{shift-by-zeros}, removing those $C_i$ that are equal to $0$ changes
neither the hypothesis nor the conclusions.  We may therefore assume that
$0 \nmid B$.  Let $y \in \mathsf Z_2 (D)$.  Then
$z = yxC_1 \cdot \ldots \cdot C_e$ is a factorization of $B$ of length
$|x|+e+2$, and multiplying a factorization of $D\pi(x)$ of maximal length by
$C_1 \cdot \ldots \cdot C_e$ shows
$\max \mathsf L (B) \ge \max \mathsf L \big( D\pi(x) \big) + e$.  It therefore suffices
to prove $\max \mathsf L (B) \le |x|+e+4$ in case 1 and
$\min \mathsf L (B) \ge |x|+e+2$ in case 2.

\emph{A graph construction.}  Since $|z|+1 = |x|+e+3 \notin \mathsf L (B)$, 
Lemma~\ref{lem:hereditary-gap}.1 applies to every divisor of $z$.  For the
divisors $C_iW$, $C_iWW'$, $yxC_i$ and $C_iC_j$, $C_iC_jW$, $yxC_iC_j$, where
$i \ne j$ and $W, W' \in \mathcal A (G)$ with $WW' \mid yx$, it gives
\begin{equation}
 \begin{aligned}
 &3 \notin \mathsf L \big( \pi (C_iW) \big), &\qquad
 &3 \notin \mathsf L (C_iC_j),\\
 &4 \notin \mathsf L \big( \pi (C_iWW') \big), &\qquad
 &4 \notin \mathsf L \big( \pi (C_iC_jW) \big),\\
 &|x|+4 \notin \mathsf L \big( D\pi(x)C_i \big), &\qquad
 &|x|+5 \notin \mathsf L \big( D\pi(x)C_iC_j \big),
 \end{aligned}\label{eq:extension-conditions}
\end{equation}
and the right-hand conditions also with $i = j$ whenever $C_i$ occurs at least
twice among $C_1, \ldots, C_e$.  Let $\mathcal V$ be the set of all nonzero
atoms $C$ satisfying the left-hand conditions of
\eqref{eq:extension-conditions} in place of $C_i$, and consider the graph with
vertex set $\mathcal V$ in which $C, C' \in \mathcal V$ are joined by an edge (a
loop if $C = C'$) if they satisfy the right-hand conditions in place of
$C_i, C_j$.  Then
$C_1, \ldots, C_e$ lie in a maximal clique $\mathcal K$ of this graph, and
every atom of $\mathcal K$ without a loop occurs at most once among
$C_1, \ldots, C_e$.  There are $622$ maximal cliques over the $216$ critical
pairs of case 1, and $3336$ over the $436$ of case 2.

Now let $\mathcal R \subset \mathcal K$
be the set of atoms with a loop, $\mathcal F = \mathcal K \setminus \mathcal R$,
and $P = D\pi(x) \prod_{C \in \mathcal F} C$.  For each maximal clique
$\mathcal K$ there is a map $w \colon G \to \Q$ such that 
\begin{enumerate}
\item[(i)] $\sigma \big( w(W) \big) \ge 1$ in case 1, and $\sigma \big( w(W) \big) \le 1$ in case 2, for every
      $W \in \mathcal A (G)$ dividing $P \bigl(\prod_{C \in \mathcal R} C\bigr)^N$
      for some $N \in \N$,
\item[(ii)] $\sigma \big( w(C) \big) = 1$ for every $C \in \mathcal R$, and
\item[(iii)] $\sigma \big(w(P) \big) < |x|+5+|\mathcal F|$ in case 1, and
      $\sigma \big( w(P) \big) > |x|+1+|\mathcal F|$ in case 2.
\end{enumerate}
Only the values of $w$ on the terms of $P$ and of the atoms of $\mathcal R$
matter, and (i) concerns finitely many atoms, so (i) to (iii) is a finite
system of linear inequalities.  We find a solution by linear programming for
each of the $3958$ maximal cliques and check (i) to (iii) in exact rational
arithmetic.

Since the atoms of $\mathcal F$ occur at most once among
$C_1, \ldots, C_e$, the sequence $B$ divides
$P \bigl(\prod_{C \in \mathcal R} C\bigr)^N$ for large $N$; hence every atom
$W \mid B$ satisfies the hypothesis of (i), and so does every
$C \in \mathcal F$.  By (ii),
\[
 \sigma \big(w(B) \big)  - \sigma \big( w(P) \big)
 = \sum_{i=1}^{e} \sigma \big( w(C_i) \big) - \sum_{C \in \mathcal F} \sigma \big( w(C) \big)
 = e - |\mathcal F| + \sum_{\substack{C \in \mathcal F\\ C \notin \{C_1, \ldots, C_e\}}} \bigl(1 - \sigma \big( w(C) \big)\bigr) .
\]
In case 1, each summand is at most $0$, so
$\sigma \big( w(B) \big)  \le \sigma \big( w(P) \big) + e - |\mathcal F| < |x|+e+5$ by (iii), and 
Lemma~\ref{lem:weight-bound}.1 gives $\max \mathsf L (B) \le \sigma \big( w(B) \big) < |x|+e+5$.
In case 2, each summand is at least $0$, so
$\sigma \big( w(B) \big) \ge \sigma \big( w(P) \big) + e - |\mathcal F| > |x|+e+1$ by (iii), and 
Lemma~\ref{lem:weight-bound}.2 gives $\min \mathsf L (B) \ge \sigma \big( w(B) \big) > |x|+e+1$.

The database \code{extensions.jsonl} contains, for each of the $652$ critical
pairs, its case, the graph, and, for each maximal clique, the map $w$.  No
bound on $e$ and no bound on the multiplicities of the $C_i$ is used.
\end{proof}

\subsection{The maximal number of length-two atoms in a factorization}\label{sec:software:length-two-atoms}

For $B \in \mathcal B (G)$ we define
\begin{equation}
 p(B)=\sum_{\{g,-g\} \subset G^{\bullet}}
       \min\{\mathsf v_g(B),\mathsf v_{-g}(B)\} ,       \label{eq:p-of-B}
\end{equation}
where the sum runs over the $13$ sets $\{g,-g\}$ with $g \in G^{\bullet}$.
Since the atoms of length two are the sequences $(-g)g$ with
$g \in G^{\bullet}$, $p(B)$ is the maximal number of atoms of length two in a
factorization of $B$.  Removing one term from $B$ lowers at most one of the
summands, and by at most one, whence
\begin{equation} \label{eq:removing-one-term}
p(B) \le p(BA^{-1}) + |A| \quad \text{ for every $A \in \mathcal B (G)$ with $A \mid B$} \,.
\end{equation}
For every $k \in \N$, we define a map
\[
\lex_k \colon \mathsf Z_k (G) \to \N^k, \quad V_1 \cdot \ldots \cdot V_k \longmapsto (|V_{\tau (1)}|, \ldots, |V_{\tau (k)}|) \,,
\]
where $\tau \colon [1,k] \to [1,k]$ is a permutation such that
$|V_{\tau (1)}| \ge \ldots \ge |V_{\tau (k)}|$, and we compare elements of
$\N^k$ lexicographically.  Then, for all $j \in \N$, $y, y' \in \mathsf Z_j (G)$
and $x \in \mathsf Z (G)$,
\begin{equation}
 \lex_j (y) < \lex_j (y')
 \quad \text{implies} \quad
 \lex_{j+|x|} (yx) < \lex_{j+|x|} (y'x) . \label{eq:lex-monotone}
\end{equation}

A subset $\mathbb H \subset G^{\bullet}$ is called a \emph{half-basis} if
$G^{\bullet} = \mathbb H \sqcup (- \mathbb H)$.  Let $\mathbb H$ be a
half-basis.  For $S \in \mathcal F (G)$, the subsequence $S_{\mathbb H}$ of $S$ is
defined by \eqref{sequence-sub}, so that
$|S_{\mathbb H}| = \sum_{g \in \mathbb H} \mathsf v_g (S)$ and $|(SS')_{\mathbb H}| = |S_{\mathbb H}| + |S'_{\mathbb H}|$ for
all $S, S' \in \mathcal F (G)$.  For an atom $A \in \mathcal A (G)$ we put
\begin{equation}
 e_{\mathbb H} (A)=|A_{\mathbb H}|-(7-|A|) ,\label{eq:orientation-excess}
\end{equation}
and we call $e_{\mathbb H} (A)$ the \emph{excess} of $A$ under $\mathbb H$.  By
\eqref{eq:p-of-B},
\begin{equation}
 p(B) \le |B_{\mathbb H}| \qquad \text{for every half-basis $\mathbb H$ and every $B \in \mathcal B (G)$} , \label{eq:budget-by-choice}
\end{equation}
because $\mathbb H$ contains exactly one element of every set $\{g,-g\}$,
and the smaller of $\mathsf v_g (B)$ and $\mathsf v_{-g} (B)$ is at most
$\mathsf v_g (B)$.

\begin{lemma}\label{lem:comp:length-two-atoms}
Let $k \in \N$ and $V_1,\ldots, V_k\in \mathcal A(G)$ with
$|V_1|\ge \ldots \ge |V_k|\ge 4$ and $|V_1|\ge 5$, and put
$B = V_1\cdot\ldots\cdot V_k$ and $z = V_1\cdot\ldots\cdot V_k \in \mathsf Z_k (B)$.
Suppose that
\[
 k=\min \mathsf L(B), \qquad k+1\not\in \mathsf L(B), \qquad \text{and} \qquad
 \lex_k(z)\ge \lex_k(z') \ \text{ for all } z'\in \mathsf Z_k(B) .
\]
Then
\[
 p(B) \le 7k - |B| + 1 .
\]
\end{lemma}

\begin{proof}
For atoms $V, W \in \mathcal A (G)$, not necessarily distinct, we consider
the conditions
\begin{enumerate}
\item[(a)] $3 \notin \mathsf L (VW)$, and
\item[(b)] $\lex_2 (VW) \ge \lex_2 (y)$ for all $y \in \mathsf Z_2 (VW)$,
      that is, no factorization of $VW$ of length two contains an atom of length
      greater than $\max \{|V|,|W|\}$.
\end{enumerate}
Only the following consequence of the hypotheses will be used.

\smallskip
\noindent
We first check the conditions (a) and (b): for all distinct $i, j \in [1,k]$,
the atoms $V_i$ and $V_j$ satisfy (a) and (b).  Indeed, (a) holds by 
Lemma~\ref{lem:hereditary-gap}.1, applied to the divisor $V_iV_j$ of $z$.  If some
$y \in \mathsf Z_2 (V_iV_j)$ had $\lex_2 (y) > \lex_2 (V_iV_j)$, then
$z' = y\,(V_iV_j)^{-1}z \in \mathsf Z_k (B)$ would satisfy
$\lex_k (z') > \lex_k (z)$ by \eqref{eq:lex-monotone}, contradicting the
choice of $z$.  Thus (b) holds.

The quantity $p$, lengths, length sets and the
conditions (a) and (b) are invariant under automorphisms of $G$.  Thus, by
Theorem~\ref{3.3}, we may assume that $U = V_1$ is one of the nine normal
forms $1(a)$, $1(b)$, $1(c)$, $2(a)$, $2(b)$, $2(c)$, $3(a)$, $3(b)$, $3(c)$
of length at least five.  Let $\Gamma_U$ be the graph whose vertices are the
atoms $V$ with $|V| \in [4,|U|]$ such that $U$ and $V$ satisfy (a) and (b), in
which two vertices $V$ and $W$ are joined by an edge (a loop if $V = W$) if
$V$ and $W$ satisfy (a) and (b).  By the argument above, $V_2, \ldots, V_k$ lie in a
maximal clique $\mathcal K$ of $\Gamma_U$ (recall from
Subsection~\ref{sec:graphs} that loops play no role for cliques), and every
atom of $\mathcal K$ without a loop occurs at most once among
$V_2, \ldots, V_k$.  Let $\mathcal R \subset \mathcal K$ be the set of atoms
with a loop, and $\mathcal F = \mathcal K \setminus \mathcal R$.

Next we show that for every maximal
clique $\mathcal K$ of $\Gamma_U$ there is a half-basis $\mathbb H$ with
\begin{equation}
 e_{\mathbb H} (V) \le 0 \ \text{ for all } V \in \mathcal R
 \qquad \text{and} \qquad
 e_{\mathbb H} (U) + \sum_{V \in \mathcal F} \max \{ e_{\mathbb H} (V), 0 \} \le 1 . \label{eq:clique-choice}
\end{equation}
Given this, \eqref{eq:budget-by-choice} and \eqref{eq:orientation-excess}
yield
\[
 p(B) \le |B_{\mathbb H}| = \sum_{i=1}^k |(V_i)_{\mathbb H}| = \sum_{i=1}^k \bigl( 7 - |V_i| + e_{\mathbb H} (V_i) \bigr)
 = 7k - |B| + e_{\mathbb H} (U) + \sum_{i=2}^k e_{\mathbb H} (V_i) \le 7k - |B| + 1 ,
\]
since the atoms of $\mathcal R$ contribute at most $0$ and every atom of
$\mathcal F$ occurs at most once among $V_2, \ldots, V_k$.

It remains to verify \eqref{eq:clique-choice} for every maximal clique of
each of the nine graphs $\Gamma_U$.  To determine $\Gamma_U$, we read off
whether two atoms satisfy (a) from the length sets of
Lemma~\ref{lem:comp:length-sets-of-pairs}, and we decide whether they satisfy
(b) by computing the factorizations of length two of their product.  The
graphs $\Gamma_U$ have the following numbers of vertices, edges, loops and
maximal cliques.
\[
\begin{array}{c|ccccccccc}
 U & 1(a) & 1(b) & 1(c) & 2(a) & 2(b) & 2(c) & 3(a) & 3(b) & 3(c)\\ \hline
 \text{vertices} & 116 & 24 & 10 & 115 & 21 & 14 & 54 & 49 & 16\\
 \text{edges} & 1164 & 46 & 13 & 1085 & 46 & 22 & 288 & 132 & 38\\
 \text{loops} & 31 & 9 & 5 & 21 & 9 & 8 & 18 & 16 & 8\\
 \text{maximal cliques} & 1060 & 43 & 10 & 1362 & 33 & 17 & 486 & 106 & 21
\end{array}
\]
Searching the $2^{13}$ half-bases, we find for each of the $3138$ maximal
cliques a half-basis $\mathbb H$ satisfying \eqref{eq:clique-choice}.  Only
$13$ distinct half-bases are needed: three for $U$ of the form $1(a)$, two
each for $1(b)$ and $2(a)$, and a single one, serving all maximal cliques,
for each of the six remaining forms.

The database \code{length-two-atoms.jsonl} contains the graphs $\Gamma_U$,
their maximal cliques and the half-bases, together with, for every pair of
atoms violating (b), a factorization of length two showing this.
\end{proof}

\subsection{Databases and verification}
\label{sec:software:structural-certificates}

Table~\ref{tab:databases} lists the databases used in the proofs above.  They
are contained in the Zenodo record~\cite{Ge-Ha-Zh27b}, together with
\begin{itemize}
\item one generating script per database, named after it (for instance
      \code{critical_pair_layers.py} writes \code{critical-pair-layers.jsonl}),
\item \code{DATABASES.md}, a description of every file and its fields,
\item the verification notebook \code{computational-verification.ipynb}, and
\item the verifier \code{verify_databases.py}.
\end{itemize}

\begin{table}
\centering
\footnotesize
\setlength{\tabcolsep}{3pt}
\begin{tabular}{lll}
 database & content & used in\\ \hline
 \code{atom-catalogue.jsonl} & the $27911$ atoms, with their normal forms & all computations\\
 \code{length-sets-with-minimum-2.jsonl} & $\mathsf L (UV)$ for normal forms $U$ and atoms $V$, $|V| \le |U|$ & Lemma~\ref{lem:comp:length-sets-of-pairs}\\
 \code{consecutive-triples.jsonl} & the triples $U, V, W$, with factorizations exhibiting $E$ & Lemma~\ref{lem:comp:mixed-triples}\\
 \code{six-long-atoms.jsonl} & the families of atoms of length at least four & Lemma~\ref{lem:comp:six-long-atoms}\\
 \code{doubled-cap-length-three.jsonl} & orbits of atoms $T$ of length three, factorizations of $DT$ & Lemma~\ref{lem:comp:doubled-cap-length-three}\\
 \code{fixed-length-sets.jsonl} & length sets of fixed products, with factorizations & Lemmas~\ref{lem:comp:doubled-cap-seeds}, \ref{le3.5}, Section~\ref{6}\\
 \code{doubled-cap-cubes.jsonl} & the sequences $B$ of Step 4, with their length sets & Lemma~\ref{lem:comp:doubled-cap-cubes}\\
 \code{critical-pair-layers.jsonl} & the $2$-$4$-sequences, the layers $P_k$ and the critical pairs & Lemma~\ref{lem:comp:carrier-antichain}\\
 \code{critical-pair-cliques.jsonl} & the graphs and maximal cliques of Steps 2 to 4 & Lemma~\ref{lem:comp:carrier-antichain}\\
 \code{bounded-endpoint.jsonl} & the critical pairs of Lemma~\ref{lem:comp:carrier-minimum-side}, the sets $\mathcal E_D$ and $\mathcal Q_t$ & Lemmas~\ref{lem:comp:carrier-minimum-side}, \ref{lem:comp:bounded-endpoint}\\
 \code{extensions.jsonl} & the graphs, maximal cliques and maps $w$ & Lemma~\ref{lem:comp:extensions}\\
 \code{length-two-atoms.jsonl} & the graphs $\Gamma_U$, maximal cliques and half-bases & Lemma~\ref{lem:comp:length-two-atoms}
\end{tabular}
\caption{The databases of the Zenodo record~\cite{Ge-Ha-Zh27b}.}
\label{tab:databases}
\end{table}

All files refer to the atoms of $G$ by their index in
\code{atom-catalogue.jsonl}.  Every file begins with a metadata record
carrying its schema and the digests of its input files; for instance, the
databases of Lemmas~\ref{lem:comp:carrier-minimum-side}
to~\ref{lem:comp:extensions} are tied to the database of
Lemma~\ref{lem:comp:carrier-antichain} they were computed from.  Enumerations
carried out up to automorphisms store one representative of each orbit
together with the size of the orbit.

Apart from Jupyter, the verification notebook needs only the Python standard
library, and it runs in about a minute.  It checks every stored atom, multiplies out every stored
factorization, and compares the stored length sets and counts with the
statements of this section.  It takes the stored data as given: that an
enumeration is exhaustive, or that a length does not occur, is not re-derived
by the notebook, and for the databases of
Lemmas~\ref{lem:comp:length-sets-of-pairs} to~\ref{lem:comp:doubled-cap-cubes}
it rests on the generators.  The databases of
Lemmas~\ref{lem:comp:carrier-antichain} to~\ref{lem:comp:length-two-atoms}
are in addition re-derived from their inputs by the verifier, which does not
use the code of the generators and has one section per database.  It recomputes the layers and
critical pairs of Step 1 of Lemma~\ref{lem:comp:carrier-antichain}, the
graphs, maximal cliques, orbits, bounds and parity arguments of
Lemmas~\ref{lem:comp:carrier-antichain} to~\ref{lem:comp:length-two-atoms},
and all stored length sets with the factorization solver; it multiplies out
every stored factorization, and it checks the maps $w$ of
Lemma~\ref{lem:comp:extensions}, which are found by linear programming (with
the HiGHS solver of SciPy), over the rationals and exact arithmetic.
With $14$ worker
processes on an $8$-core machine, regenerating all databases takes
about 15 minutes, and running the verification script takes about $10$ minutes.
In all cases, a stored factorization shows that a
length occurs. The fact that a length does \emph{not} occur,
and that the computed length sets are complete, rest on
an exhaustive computation or on Lemma~\ref{lem:weight-bound}.

\section{Precise Structural Properties of General Length Sets} \label{5}

\smallskip

\centerline{\it Throughout this section, let $G = C_3 \oplus C_3 \oplus C_3$.}

\smallskip

\subsection{On small divisors $E$ of zero-sums $B$ with $\Delta \big( \mathsf L (E) \big)$ reflecting $\Delta \big( \mathsf L (B) \big)$} \label{sec:structural-lemmas}

\begin{lemma} \label{lem:first-jump-carrier}
Let $B \in \mathcal B (G)$ and let $k \in \N$ such that $k-1,k+1\in\mathsf L(B)$ and
$k\notin\mathsf L(B)$.  Then there are $z \in \mathsf Z (B)$ with $|z| = k-1$
and $U, V \in \mathcal A (G)$ with $UV \mid z$ such that
\[
 4 \in \mathsf L (UV)
 \qquad \text{and} \qquad
 3 \notin \mathsf L (UV) \,.
\]
In particular, $\pi (UV)$ is a $2$-$4$-sequence in the sense of
Subsection~\ref{sec:software:carriers}.
\end{lemma}

\begin{proof}
Since $\mathsf c (G) = 4$  by Lemma \ref{lem:basic-properties}, each two factorizations of $B$ of lengths $k-1$ and
$k+1$ can be concatenated by a $4$-chain. Since $k \notin \mathsf L (B)$, there are thus $z, z' \in \mathsf Z (B)$ with
$|z| = k-1$,  $|z'| = k+1$, and $\mathsf d (z,z') \le 4$.  If $y = z\gcd (z,z')^{-1}$ and
$y' = z'\gcd (z,z')^{-1}$, then $|y'|-|y| = 2$ and
$\max \{|y|,|y'|\} = \mathsf d (z,z') \le 4$.  Moreover $|y| \ge 2$, since
$\pi(y) = \pi(y')$ and $y, y'$ have no common nontrivial divisor: $|y|=0$ would force
$y'=1$, and $|y|=1$ would make the atom $\pi(y)$ a product of $|y'| \ge 2$
atoms.  It follows that $|y|=2$ and $|y'|=4$, whence $y = UV$ with $U, V \in \mathcal A (G)$ and
$\{2,4\} \subset \mathsf L (UV)$.  If $3 \in \mathsf L (UV)$, then it follows from $k-1=|y|+|\gcd (z,z')|=2+|\gcd (z,z')|$ that $k=3+|\gcd (z,z')|\in \mathsf L(B)$, a contradiction. Thus $3\not\in \mathsf L (UV)$ and the assertion follows.
\end{proof}

\begin{lemma} \label{lem:bounded-mixed-divisor}
Let $B \in \mathcal B (G)$ with $1 \in \Delta \big(\mathsf L (B) \big)$.  Then there is
a divisor $E$ of $B$ with $1 \in \Delta \big(\mathsf L (E) \big)$ and 
$1 \notin \Delta \big(\mathsf L (E') \big)$ for every proper divisor $E'$ of $E$, 
such that
\begin{equation}
\{k,k+1\}\subset\mathsf L(E) \quad \text{for some }k\in\{2,3\}
 \quad \text{and} \quad  |E| \le14.                                   \label{eq:mixed-core-bound}
\end{equation}
If $k=3$, then there is a basis $(e_1,e_2,e_3)$ of $G$ such that
\[
A=e_1e_2e_3(-e_1-e_2-e_3) , \qquad E=A^3, \quad \text{and} \quad \mathsf L(E)=\{3,4\}.
\]
\end{lemma}

\begin{proof}
Choose a divisor $E$ of $B$ with $1 \in \Delta \big(\mathsf L (E) \big)$, which is
minimal with respect to divisibility. Then $1 \notin \Delta \big(\mathsf L (E') \big)$
for every proper divisor $E'$ of $E$.  Since $\mathsf c (G) = 4$ by Lemma \ref{lem:basic-properties}, two
factorizations of $E$ whose lengths have different parity can be concatenated
by a $4$-chain, and some step of this chain joins two factorizations $z$ and
$z'$ of $E$ whose lengths have different parity.  Set $y = z\gcd (z,z')^{-1}$
and $y' = z'\gcd (z,z')^{-1}$. Then $|y|$ and $|y'|$ have different parity,
$\max \{|y|,|y'|\} \le 4$ by $\mathsf d (z,z') \le 4$, and $|y|,|y'| \ge 2$ by $\pi(y)=\pi(y')$ with $|y|\neq |y'|$.  Therefore
$\{|y|,|y'|\} = \{2,3\}$ or $\{|y|,|y'|\} = \{3,4\}$, so that
$\pi (y) = \pi (y')$ is a divisor of $E$ with $1\in \Delta \big(\mathsf L ( \pi(y) ) \big)$.
 By the minimality of $E$, we obtain
$\pi (y) = E$ and  $\{k,k+1\} \subset \mathsf L (E)$ with $k\in \{2,3\}$.

For $k=2$, the sequence $E$ is a product of two atoms, whence
$|E| \le 2\mathsf D(G)=14$.  For $k=3$, the assertion follows from Lemma~\ref{lem:comp:mixed-triples}.
\end{proof}

\begin{lemma} \label{lem:bounded-context}
Let $D$ be a $2$-$4$-sequence, and let $x \in \mathsf Z(G)$ be such that
$1\in \Delta \big(\mathsf L(D\pi(x)) \big)$ and $1\notin \Delta \big(\mathsf L(D\pi(x')) \big)$ for every proper divisor $x'$ of $x$.  Then $|x| \le 14$.
In particular, the critical pairs are precisely the $652$ listed in
Lemma~\ref{lem:comp:carrier-antichain}.
\end{lemma}

\begin{proof}
Since $1 \in \Delta \big(\mathsf L (D\pi(x)) \big)$,
Lemma~\ref{lem:bounded-mixed-divisor} provides a divisor
$E$ of $D\pi(x)$ with $1 \in \Delta \big(\mathsf L (E) \big)$ and $|E| \le 14$.  Let
$y = UV \in \mathsf Z (D)$ of length $2$ and let $x = X_1 \cdot \ldots \cdot X_{|x|}$. If there exists $i\in [1, |x|]$, say $i=1$, such that $E$ divides $UVX_2\cdot \ldots \cdot X_{|x|}$, then 
$1 \in \Delta (\mathsf L \big(D \pi (x X_1^{-1})) \big)$, contradicting the minimality
of $x$.  Thus $E\nmid D\pi(xX_i^{-1})$ for all $i \in [1,|x|]$.  Now
fix $E$ as a subsequence of $UVX_1 \cdot \ldots \cdot X_{|x|}$ and assign each
of its terms to an atom containing it.  If $|x| > |E|$, then some $X_i$
receives no term, whence $E \mid D\pi(xX_i^{-1})$, a contradiction.  Thus
$|x|\le |E|\le 14$.
\end{proof}

\subsection{On the structure  of length sets $L$ with $1 \in \Delta (L)$}

\begin{theorem} \label{thm:mixed-shape}
Let $L \in \mathcal L (G)$ with $1 \in \Delta ( L )$.
 Then, either
\begin{equation}
L  =[\min L, \max L] \quad\text{or} \quad
L =[\min L ,\max  L -2]\cup\{\max  L \}
 \quad\text{or}\quad
L =\{\min  L \}\cup[\min  L +2,\max  L].                     \label{eq:mixed-shapes}
\end{equation}
\end{theorem}

\begin{proof}
	
In order to prove
\eqref{eq:mixed-shapes}, we proceed in a series of assertions.

\begin{claim}\label{A:protected-complement}
Let $B \in \mathcal B (G)$ be minimal with respect to divisibility among all
sequences $B_0 \in \mathcal B (G)$ with $1 \in \Delta \big(\mathsf L (B_0) \big)$ and
\[
 \min\mathsf L(B_0)+1,\ \max\mathsf L(B_0)-1\notin\mathsf L(B_0).
\]
Among all  factorizations of $B$ with minimal length and all their divisors, choose a
divisor $z_F$ with minimal length such that 
$1 \in \Delta \big(\mathsf L (\pi (z_F)) \big)$ and let $z$ be a factorization of $B$ with minimal length for which the divisor
$z_F$ was chosen.  Then $zz_F^{-1}$ contains no atom of length at most three, and $|zz_F^{-1}|\le 5$.
\end{claim}

\begin{proof}[Proof of \ref{A:protected-complement}]
Let $m = \min \mathsf L (B)$, $M = \max \mathsf L (B)$,  and $F = \pi  (z_F)$. Suppose $zz_F^{-1}=U_1\cdot \ldots \cdot U_s$, where $s=|zz_F^{-1}|$. It suffices to show that $s\le 5$ and $|U_i|\ge 4$ for all $i\in [1,s]$. 

We first prove the stronger assertion that for every factorization $z^*=W_1\cdot\ldots\cdot W_s$ of $\pi(zz_F^{-1})$ with minimum length $s$, we have $|W_i|\ge 4$ for all $i\in [1,s]$. 
Assume to the contrary that there exists $i\in [1,s]$, say $i=s$, such that $|W_s|\le 3$. By Lemma \ref{A:short-atom}, 
there exists a factorization $z'=V_1\cdot\ldots\cdot V_{M}$ of $B$  with maximal length such that $V_1=W_s$.  Note that $M-1=\max\mathsf L(BW_s^{-1})$ and $m-1=\min\mathsf L(BW_s^{-1})$ by Lemma~\ref{lem:hereditary-gap}.1. 
It follows from $F\t BW_s^{-1}$, $m+1\notin \mathsf L(B)$, and $M-1\not\in \mathsf L(B)$ that $1\in \Delta \big(\mathsf L(BW_s^{-1}) \big)$, $m\not\in \mathsf L(BW_s^{-1})$, and $M-2\notin \mathsf L(BW_s^{-1})$, a contradiction to the minimality of $B$. Thus  $|W_i|\ge 4$ for all $i\in [1,s]$. 

We next show $s\le 5$.
Assume to the contrary that $s\ge 6$. Note that $s=\min \mathsf L \big(\pi(zz_F^{-1}) \big)$ and $s+1\not\in \mathsf L \big(\pi(zz_F^{-1}) \big)$. Then
Lemma~\ref{lem:comp:six-long-atoms}, applied to $\pi(zz_F^{-1})$, its factorization
$zz_F^{-1}$ of minimal length $s$, and the atoms $U_1, \ldots, U_6$, which satisfy
$U_1 \cdot \ldots \cdot U_6 \mid zz_F^{-1}$ and $|U_i| \ge 4$, provides a
subset $I \subset [1,6]$ and atoms $V_1, \ldots, V_k$ with
$\prod_{i \in I}U_i = V_1 \cdot \ldots \cdot V_k$ and $|V_1| \le 3$.

If $k = |I|$, then replacing in $zz_F^{-1}$ the atoms $U_i$ with $i \in I$ by
$V_1, \ldots, V_k$ gives a factorization $z^*$ of $\pi(zz_F^{-1})$ with $|z^*| = s$. But $|V_1|\le 3$,  contradicting  the stronger assertion.

If $|V_j| \le 3$ for all $j \in [1,k]$, then  by repeatedly applying Lemma \ref{A:short-atom},  we obtain that $B$ has a factorization $V_1\cdot\ldots\cdot V_kV_{k+1}\cdot\ldots\cdot V_M$ with maximum length. Note that $M-k=\max\mathsf L(B(V_1\cdot\ldots\cdot V_k)^{-1})$ and $m-|I|=\min\mathsf L(B(\prod_{i\in I}U_i)^{-1})$ by Lemma~\ref{lem:hereditary-gap}.1. 
It follows from $F\t B(V_1\cdot\ldots\cdot V_k)^{-1}$, $m+1\notin \mathsf L(B)$, and $M-1\not\in \mathsf L(B)$ that $1\in \Delta \big(\mathsf L(B(V_1\cdot\ldots\cdot V_k)^{-1}) \big)$, $m-|I|+1\not\in \mathsf L(B(\prod_{i\in I}U_i)^{-1})$, and $M-k-1\notin \mathsf L(B(V_1\cdot\ldots\cdot V_k)^{-1})$. Since $B(V_1\cdot\ldots\cdot V_k)^{-1}=B(\prod_{i\in I}U_i)^{-1}$, we obtain a contradiction to the minimality of $B$. 
\qedhere\ (\ref{A:protected-complement})
\end{proof}

\begin{claim}\label{A:no-double-gap}
Let $B\in\mathcal B(G)$ with $1 \in \Delta \big(\mathsf L (B) \big)$.  Then either
$\min\mathsf L(B)+1\in\mathsf L(B)$, or $\max\mathsf L(B)-1\in\mathsf L(B)$.
\end{claim}

\begin{proof}[Proof of \ref{A:no-double-gap}]
	Let
	$m=\min\mathsf L(B)$ and $M=\max\mathsf L(B)$.
Assume to the contrary that both $m+1$ and $M-1$ are missing, and choose $B$
minimal with respect to divisibility with this property.  Among all 
factorizations of $B$ with minimal length and all their divisors, choose a divisor $z_F$
with minimal length such that $1 \in \Delta \big(\mathsf L (\pi (z_F)) \big)$, and let $z$ be the  factorization of $B$ with minimal length for which the divisor $z_F$ was chosen.

Set $F = \pi (z_F)$.  By  Lemma~\ref{lem:hereditary-gap}.1, we have
$|z_F| = \min \mathsf L (F)$ and $|z_F|+1 \notin \mathsf L (F)$, and no
proper divisor $z'$ of $z_F$ satisfies
$1 \in \Delta \big(\mathsf L (\pi (z')) \big)$.

Since $1 \in \Delta (\mathsf L (F) )$, we have $|\mathsf L (F)|\ge 2$.  Since
$\min \mathsf L (F) = |z_F|$ and $|z_F| + 1 \notin \mathsf L (F)$ while
$\max\Delta (G)=2$, we obtain $|z_F| + 2\in\mathsf L (F)$,
and may apply Lemma~\ref{lem:first-jump-carrier} (with $k=|z_F| + 1$).  After
replacing $z_F$ by a factorization of $F$ of the same length (and $z$
accordingly), there exist a $2$-$4$-sequence $D$ and $x\in \mathsf Z(G)$ such that $F=D\pi(x)$ and such that $(D,x)$ satisfies
\eqref{eq:carrier-extension-state} with $|z_F| = |x|+2$.
  By
Lemma~\ref{lem:bounded-context}, $(D,x)$ is one of the critical pairs of
Lemma~\ref{lem:comp:carrier-antichain}, and the choice of $z_F$ shows that
the hypotheses of Lemma~\ref{lem:comp:carrier-minimum-side} are satisfied.
Hence $|x| \le 3$, that is, $|z_F| \le 5$.  If $|z_F| = 5$, then
Lemma~\ref{lem:comp:carrier-minimum-side}.2 provides a factorization $w$
of $F$ with $|w| = 5$ and some $T\in \mathcal A(G)$ with $T \mid w$ with $|T| \le 3$ such that
$1 \in \Delta \big( \mathsf L (\pi (T^{-1}w)) \big)$. Since $T^{-1}w$ is a divisor of length four of a  factorization of $B$ with minimal length, this contradicts the
choice of $z_F$.  Thus $|z_F| \le 4$.  

By Claim~\ref{A:protected-complement}, we obtain  $|zz_F^{-1}|\le 5$  and all atoms in the factorization $zz_F^{-1}$ have length at least $4$, whence 
$m=|z_F|+|zz_F^{-1}| \le 4+5 = 9$, and the  factorization of $B$ with minimal length has the form
\[
 z = y W_1 \cdot \ldots \cdot W_s C_1\cdot \ldots \cdot C_t, \qquad\text{ with }
 \qquad s\le2,
 \qquad t \in [0, 5],
\]
where $y \in \mathsf Z_2 (D)$, where $x = W_1 \cdot \ldots \cdot W_s$ is such
that $(D,x)$ is one of the $193$ critical pairs of
Lemma~\ref{lem:comp:carrier-minimum-side}.1 (by the choice of $z_F$, the
displayed $y$ satisfies the hypotheses of Lemma~\ref{lem:comp:carrier-minimum-side}), and where $|C_i| \ge 4$ for all
$i \in [1,t]$.  Since $z$ is a factorization of $B$ with minimal length and
$m+1 \notin \mathsf L (B)$, Lemma~\ref{lem:comp:bounded-endpoint} gives
$\max \mathsf L (B) - 1 \in \mathsf L (B)$, contradicting
$M-1 \notin \mathsf L (B)$.
\qedhere\ (\ref{A:no-double-gap})
\end{proof}

\begin{claim}\label{A:endpoint-hole}
Let $B \in \mathcal B (G)$ with $1 \in \Delta \big(\mathsf L (B) \big)$, and let $h\in [\min\mathsf L(B), \max\mathsf L(B)]$ be such that  $h\not\in\mathsf L (B)$.  Then $h = \min \mathsf L (B) +1$ or $h = \max \mathsf L (B)-1$.
\end{claim}

\begin{proof}[Proof of \ref{A:endpoint-hole}]
Since $\Delta (G) = \{1,2\}$, we have $h-1, h+1 \in \mathsf L (B)$.  By
Lemma~\ref{lem:first-jump-carrier} there are $z_0 \in \mathsf Z (B)$ with
$|z_0| = h-1$ and $y \in \mathsf Z (G)$ with $y \mid z_0$ and $|y| = 2$ such
that $D = \pi (y)$ is a $2$-$4$-sequence.  Let $x = z_0y^{-1}$, so that
$B = D\pi(x)$ and $|x| = h-3$. 

Among the divisors $x^*$ of $x$ with $1\in \Delta \big(\mathsf L(D\pi(x^*)) \big)$  choose one, say $x'$, which is
minimal with respect to divisibility, and set $x'' = x x'^{-1}$,
and $e = |x''|$.  Then
\[
 B = D\pi(x') \, \pi (x''),
 \qquad
 h = |x|+3 = |x'|+e+3 ,
\]
and $|x'|+e+3 \notin \mathsf L (B)$.  We claim that $(D,x')$ is one of the
critical pairs to which Lemma~\ref{lem:comp:extensions} refers.  Indeed, $1\in \Delta \big(\mathsf L(D\pi(x')) \big)$ and $1\notin \Delta \big(\mathsf L(D\pi(x^*)) \big)$ for every proper divisor $x^*$ of $x'$ by the choice of $x'$, and $|x'|+3\not\in \mathsf L \big(D\pi(x') \big)$,
because $|x'|+3 \in \mathsf L \big(D\pi(x') \big)$ would give
$|x'|+3+e = h \in \mathsf L (B)$.
Thus $(D,x')$ satisfies \eqref{eq:carrier-extension-state}, and by
Lemma~\ref{lem:bounded-context}, it is one of the $652$ critical pairs of
Lemma~\ref{lem:comp:carrier-antichain}. Suppose $C_1, \ldots, C_e\in \mathcal A(G)$ such that $x''=C_1\cdot \ldots\cdot C_e$.  By
Lemma~\ref{lem:comp:extensions}, $\max \mathsf L \big(D\pi(x') \big) = |x'|+4$ or
$\min \mathsf L \big(D\pi(x') \big) = |x'|+2$.  In the first case,
Lemma~\ref{lem:comp:extensions}.1 gives
$\max \mathsf L (B) = \max \mathsf L (D\pi(x')C_1\cdot \ldots \cdot  C_e) = |x'|+e+4 = h+1$.  In the
second case,  Lemma~\ref{lem:comp:extensions}.2 gives
$\min \mathsf L (B) = \min \mathsf L (D\pi(x')C_1\cdot \ldots \cdot  C_e) = |x'|+e+2 = h-1$.\qedhere\
(\ref{A:endpoint-hole})
\end{proof}

Let $B \in \mathcal B (G)$ with $\mathsf L (B)=L$, and set $m=\min \mathsf L (B)$ and $M = \max \mathsf L (B)$.
Assume that $L\neq [m,M]$. Let $h\in [m,M]$ be such that $h\notin L$. 
By
Claim~\ref{A:endpoint-hole}, we obtain $h\in \{m+1, M-1\}$, and by Claim~\ref{A:no-double-gap}, either $m+1\in L$ or $M-1\in L$.
This gives \eqref{eq:mixed-shapes}.  
\end{proof}

\smallskip
\subsection{On the elasticity of length sets $L$ with $\min L + 1 \notin L$}

We now prove an estimate for $\max L$ when $\min L+1\not\in L$.  We use the quantity $p(B)$ from
\eqref{eq:p-of-B}.
Among all  factorizations $z$ with minimal length, we choose one for which
$\lex_{|z|} (z)$ is maximal, with $\lex_{|z|}$ as defined in
Subsection~\ref{sec:software:length-two-atoms}.  This choice is
inherited by the divisors of $z$:  if $y \mid z$ and if
$y' \in \mathsf Z_{|y|} \big(\pi (y) \big)$ satisfies $\lex_{|y|} (y') > \lex_{|y|} (y)$,
then $y'(y^{-1}z)$ is a factorization of $\pi (z)$ of length $|z|$ with
$\lex_{|z|} \big(y'(y^{-1}z) \big) > \lex_{|z|} (z)$ by \eqref{eq:lex-monotone}.  Together with Lemma~\ref{lem:hereditary-gap}.1 this
shows:  if $|z| = \min \mathsf L \big(\pi (z) \big)$, $|z|+1 \notin \mathsf L \big(\pi (z) \big)$
and $\lex_{|z|} (z)$ is maximal, then every divisor $y$ of $z$ satisfies
$|y| = \min \mathsf L \big(\pi (y) \big)$, $|y|+1 \notin \mathsf L \big(\pi (y) \big)$, and
$\lex_{|y|} (y) \ge \lex_{|y|} (y')$ for all $y' \in \mathsf Z_{|y|} \big(\pi (y) \big)$.

\begin{theorem} \label{thm:lower-gap-endpoint}
Let $L \in\mathcal L (G)$ with $\min L  +1\notin  L$.  Then
\begin{equation}
\max  L  \le \frac{7}{3}\min L  + \frac{1}{3}                                            \label{eq:lower-gap-bound}
\end{equation}
Moreover, if $\Delta (L) = \{2\}$, then $\rho (L) \le 7/3 < 7/2 = \rho (G)$.
\end{theorem}

\begin{proof}
Let $B \in \mathcal B (G)$ with $\mathsf L (B) = L$ and  let $m=\min\mathsf L(B)$,  $M=\max\mathsf L(B)$. By Lemma \ref{lem:basic-properties}, we have $\rho (G)=7/2$.  By \eqref{shift-by-zeros}, we may suppose that  $0\nmid B$. Choose a  factorization $z$ of $B$ with $|z|=\min \mathsf L (B)$ for
which $\lex_{|z|} (z)$ is maximal. Recall the quantity $p(B)$ from \eqref{eq:p-of-B}.  A 
factorization of $B$ with maximal length has at most $p(B)$ atoms of length two, whence
\begin{equation}
 3M\le |B|+p(B).                             \label{eq:maximum-term-bound}
\end{equation}
It suffices to show 
\[
 p(B) \le 7m-|B|+1 .                    \tag{$\ast$}
\]
Let $z'$ be the product of the atoms occurring in $z$ of length at least four,
and set $B' = \pi (z')$.  Every atom $U$ occurring in $z z'^{-1}$ satisfies
$|U| \le 3$, hence $|U| \le 7-|U|$, and summing over
these atoms gives
$|B|-|B'| \le 7|z z'^{-1}|-|B(B')^{-1}| = (7|z|-|B|) -(7|z'|-|B'|)$.
Since $p(B) \le p(B') + |B|-|B'|$ by 
\eqref{eq:removing-one-term}, it suffices to prove
$$p(B') \le 7|z'|-|B'|+1\,.$$  If $z' = 1$, this is clear.  If every
atom occurring in $z'$ has length four, then
\[
 p(B')\le |B'|/2=2|z'|\le3|z'|=7|z'|-|B'|.
\]
Otherwise suppose $z' = V_1 \cdot \ldots \cdot V_n$ with
$|V_1| \ge \ldots \ge |V_n| \ge 4$, so that $|V_1| \ge 5$.  By the discussion
preceding the theorem, applied with $|z| = m$ and
$|z|+1 = m+1 \notin \mathsf L (B)$, the divisor $z'$ of $z$ satisfies
$n = \min \mathsf L (B')$, $n+1 \notin \mathsf L (B')$, and
$\lex_n (z') \ge \lex_n (y')$ for all $y' \in \mathsf Z_n (B')$.  Hence
Lemma~\ref{lem:comp:length-two-atoms} gives
\[
 p(B') \le 7n - |B'| + 1 = 7|z'|-|B'|+1 .
\]
Thus $(\ast)$ holds, and \eqref{eq:maximum-term-bound} gives
$3M \le |B| + p(B) \le 7m+1$.

Finally suppose that $\Delta (L) = \{2\}$. We have $3 \max L   \le 7 \min L+1$.  Since $\max L - \min L$ is even, it follows that $3 \max L - 7 \min L$ is even, whence  $3 \max L \le 7 \min L$.
\end{proof}

The bound \eqref{eq:lower-gap-bound} is attained by the following explicit families.

\begin{lemma}\label{le3.5}
Let $(e_1,e_2,e_3)$ be a basis of $G$ and let
\[
U=e_1^2e_2^2e_3^2(e_1+e_2+e_3), \ V=e_1^2e_3^2(-e_2-e_3)(-e_1-e_2)(-e_1-e_2-e_3), \ \text{and} \ W=e_1^2e_2(e_1-e_2+e_3)(-e_3) \,.
\]
Then, for every $k\in \N_0$, we have
\begin{enumerate}
\item  $\mathsf L(U^{3k+1}V)=\{3k+2\}\cup [3k+4, 7k+5]$;

\item $\mathsf L(U^{3k+2}V)=\{3k+3\}\cup [3k+5, 7k+7]$;

\item $\mathsf L(U^{3k+3}V)=\{3k+4\}\cup [3k+6, 7k+9]$;

\item  $\mathsf L(U^{3k+1}W)=\{3k+2\}\cup [3k+4, 7k+4]$.
\end{enumerate}
\end{lemma}

\begin{proof}
Let $G_0=\{e_1,e_2,e_3,e_1+e_2+e_3\}$.  By  Lemma~\ref{lem:cap-atoms}.1,
every atom of $\mathcal B (G_0)$ has odd length, whence
$\mathsf L (S) \subset |S| + 2\Z$ for every $S \in \mathcal B (G_0)$ by
Lemma~\ref{lem:parity-formula}.3.  Since $U^{3k+1}\in \mathcal B(G_0)$ and $|U^{3k+1}| = 21k+7$, this gives
$3k+2\not\in \mathsf L(U^{3k+1})$ for every $k\in \N_0$. 
	Let $1 \ne T_1$ be a subsequence of $(-e_2-e_3)(-e_1-e_2)(-e_1-e_2-e_3)$ and let $T_2 \in \mathcal F (G_0)$ such that $T_1T_2$ is an atom.
	Then it is easy to see that 
	\begin{equation}\label{neq2}
		|T_1|\neq 2\,.
	\end{equation}
	Note that $\mathsf L (U^3)=\{3,5,7\}$,
	$\mathsf L (UV)=\{2\}\cup[4,5]$, $\mathsf L (U^2V)=\{3\}\cup[5,7]$, $\mathsf L(UW)=\{2,4\}$, and $\mathsf L(U^2W)=\{3\}\cup[5,6]$
	(these length sets are the records 
	\code{u1-cubed}, \code{u1-v}, \code{u1-squared-v}, \code{u1-w}, and \code{u1-squared-w} of
	\code{fixed-length-sets.jsonl}, each with a factorization witness for
	every attained length).
	
\smallskip
	1. We have $\big(\{2\}\cup [4,5]\big)+k\{3,5,7\}=\{3k+2\}\cup [3k+4, 7k+5]\subset \mathsf L(UV)+\mathsf L(U^{3k})\subset \mathsf L(U^{3k+1}V)$.
	
	Since $\mathsf D(G)=7$ and $(-e_1-e_2-e_3)(e_1+e_2+e_3)$ is the only  atom of length-$2$  dividing $U^{3k+1}V$, we obtain $\mathsf L(U^{3k+1}V)\subset[3k+2, 7k+5]$.
	It suffices to show $3k+3\not\in \mathsf L(U^{3k+1}V)$.
	
	Assume to  the contrary that there exist $W_1,\ldots, W_{3k+3}\in \mathcal A(G)$  such that $U^{3k+1}V=W_1\cdot\ldots\cdot W_{3k+3}$.
	If there exists $i\in [1,3k+3]$, say $i=1$, such that $(-e_2-e_3)(-e_1-e_2)(-e_1-e_2-e_3)$ divides $W_1$, then $W_1=V$ and hence  $3k+2\in \mathsf L(U^{3k+1})$, a contradiction. In view of \eqref{neq2}, every atom $W_i$ has at most one term from $(-e_2-e_3)(-e_1-e_2)(-e_1-e_2-e_3)$. 	
	After renumbering if necessary, we may assume that $(-e_2-e_3)\t W_1$, $(-e_1-e_2)\t W_2$, and $(-e_1-e_2-e_3)\t W_3$. Thus $|W_t|\le 4$ for each $t\in [1,3]$, whence
	\[
	|U^{3k+1}V|=21k+14=|W_1W_2W_3|+|W_4\cdot\ldots\cdot W_{3k+3}|\le 12+21k\,,
	\]
	a contradiction.

2.  We have $\big(\{3\}\cup [5,7]\big)+k\{3,5,7\}=\{3k+3\}\cup [3k+5, 7k+7]\subset \mathsf L(U^2V)+\mathsf L(U^{3k})\subset \mathsf L(U^{3k+2}V)$.
	
	Since $\mathsf D(G)=7$ and $(-e_1-e_2-e_3)(e_1+e_2+e_3)$ is the only  atom of  length-$2$  dividing $U^{3k+2}V$, we obtain $\mathsf L(U^{3k+2}V)\subset[3k+3, 7k+7]$.
	Assume to  the contrary that $3k+4\in \mathsf L(U^{3k+2}V)$. Then
	$3k+6=3(k+1)+3\in \mathsf L(U^{3k+4}V)$, a contradiction to Part 1. Thus $\mathsf L(U^{3k+2}V)=\{3k+3\}\cup[3k+5, 7k+7]$.

3.  Since $\{4\}\cup [6,9]=\{2,4\}+(\{2\}\cup [4,5])\subset\mathsf L(U^2)+\mathsf L(UV)\subset \mathsf L(U^3V)$,
	we have $\big(\{4\}\cup [6,9]\big)+k\{3,5,7\}=\{3k+4\}\cup [3k+6, 7k+9]\subset \mathsf L(U^3V)+\mathsf L(U^{3k}) \subset\mathsf L(U^{3k+3}V)$.
	
	Since $\mathsf D(G)=7$ and $(-e_1-e_2-e_3)(e_1+e_2+e_3)$ is the only  atom  of length-$2$  dividing $U^{3k+3}V$, we obtain $\mathsf L(U^{3k+3}V)\subset[3k+4, 7k+9]$.
	Assume to  the contrary that $3k+5\in \mathsf L(U^{3k+3}V)$. Then
	$3k+6=3(k+1)+3\in \mathsf L(U^{3k+4}V)$, a contradiction to Part 1. Thus $\mathsf L(U^{3k+3}V)=\{3k+4\}\cup [3k+6, 7k+9]$.

4. 
	We have that
$k\{3,5,7\}+\{2,4\}\subset \mathsf L(U^{3k})+\mathsf L(UW)\subset \mathsf L(U^{3k+1}W)$
	 and, if $k \ge 1$, that
	 $(k-1)\{3,5,7\}+\{2,4\}+\{3,5,6\}\subset \mathsf L(U^{3k-3})+\mathsf L(U^2)+\mathsf L(U^2W)\subset \mathsf L(U^{3k+1}W)$.
	 It follows that $\{3k+2\}\cup [3k+4,7k+4]\subset \mathsf L(U^{3k+1}W)$.

	Since $\mathsf D(G)=7$ and $e_3(-e_3)$ is the only atom of length $2$ dividing $U^{3k+1}W$, we obtain $\mathsf L(U^{3k+1}W)\subset[3k+2, 7k+4]$.
	It suffices to show $3k+3\not\in \mathsf L(U^{3k+1}W)$. Assume to  the contrary that there exist $W_1,\ldots, W_{3k+3}\in \mathcal A(G)$  such that $U^{3k+1}W=W_1\cdot\ldots\cdot W_{3k+3}$.
	If there exists $i\in [1,3k+3]$, say $i=1$, such that $(e_1-e_2+e_3)(-e_3)\t W_1$, then $W_1=W$ and hence  $3k+2\in \mathsf L(U^{3k+1})$, a contradiction.	
	Otherwise, after renumbering if necessary, we may assume that $(e_1-e_2+e_3)\t W_1$ and  $(-e_3)\t W_2$, whence
	$W_1\in \{(e_1-e_2+e_3)e_1^2e_2e_3^2, (e_1-e_2+e_3)(e_1+e_2+e_3)e_1e_3\}$ and  $W_2\in \{(-e_3)e_3, (-e_3)(e_1+e_2+e_3)e_1^2e_2^2\}$.
	If $W_1=(e_1-e_2+e_3)e_1^2e_2e_3^2$ or $W_2=(-e_3)(e_1+e_2+e_3)e_1^2e_2^2$, then $W_1W_2=WW'$ for some $W'\in \mathcal A(G)$ and hence  $3k+2\in \mathsf L(U^{3k+1})$, a contradiction.	
	Thus $W_1=(e_1-e_2+e_3)(e_1+e_2+e_3)e_1e_3$ and  $W_2=(-e_3)e_3$, whence
	\[
	|U^{3k+1}W|=21k+12=|W_1\cdot\ldots\cdot W_{3k+3}|\le 4+2+7\cdot (3k+1)=21k+13\,.
	\]
	It follows that there is exactly one $i\in [3,3k+3]$, say $i=3$, such that $|W_3|=6$ and $|W_j|=7$ for all $j\in [4, 3k+3]$. Since there is only one atom over $G_0$ having length $7$, we obtain that $W_j=U$ for all $j\in [4, 3k+3]$ and hence $UW=W_1W_2W_3$, a contradiction to $\mathsf L(UW)=\{2,4\}$.
\end{proof}

\smallskip
\subsection{On the elasticity of length sets $L$  with $\max L - 1 \notin L$}~

\begin{proposition} \label{prop:upper-hole-endpoint}
Let $L \in \mathcal L (G)$ with
$\max  L  -1\notin L $.  Then
\begin{equation}
\max  L  \le \frac{7}{3} \min  L  + \frac{8}{3} \,.\label{eq:cap-endpoint-bound}
\end{equation}
Consequently, if $L =[m,\nu]\cup\{\nu+2\}$, then
\begin{equation}
 \nu\le
 \begin{cases}
  7k,&m=3k,\\
  7k+3,&m=3k+1,\\
  7k+5,&m=3k+2.
 \end{cases}\label{eq:cap-preliminary-bound}
\end{equation}
\end{proposition}

\begin{proof}
Let $B \in \mathcal B (G)$ with $\mathsf L (B)=L$, and let $m=\min\mathsf L(B)$,
$M=\max\mathsf L(B)$. By  \eqref{shift-by-zeros}, we may suppose that $0 \nmid B$.
For a  factorization of $B$ with maximal length $M$,  by
Lemma~\ref{lem:upper-hole-affine-cap}, at most eight atoms occurring in this factorization  have length two.  It follows that $3M-8=2\times 8+3(M-8)\le |B|\le7m$, whence the assertion follows.
\end{proof}

\begin{lemma}\label{lem:oriented-cap-top-gap}
Let $(e_1,e_2,e_3)$ be a basis of $G$, and let   $e_0 = e_1+e_2+e_3$.
If
\[
 S=\prod_{i=0}^3 \big( e_i(-e_i) \big)^{\lambda_i}
   \prod_{i=0}^3(e_i^3)^{c_i},\quad
 \text{ where }\lambda_i\in [0,2] \text{ and } c_i\in \N_0 \text{ for all }i\in [0,3]\,,
\]
then
\[
 \max\mathsf L(S)=\sum_{i=0}^3(\lambda_i+c_i),
 \qquad \text{and} \qquad \max\mathsf L(S)-1\notin\mathsf L(S).
\]
\end{lemma}

\begin{proof}
Let $a = \sum_{i=0}^3 \lambda_i$ and $r = \sum_{i=0}^3 c_i$, so that
$|S| = 2a+3r$. Let $G_0=\{e_0,e_1,e_2,e_3\}$. Since $G_0 \cup (-G_0)$ is a cap (see Subsection~\ref{2.5}), the atoms of length two
dividing $S$ are the $e_i(-e_i)$ with $i \in [0,3]$, and the atoms of length
three dividing $S$ are the $g^3$ with $g \in G_0 \cup (-G_0)$; as
$\mathsf v_{-e_i}(S) = \lambda_i \le 2$, these are the $e_i^3$.  In
particular, a factorization of $S$ has at most $a$ atoms of length two.  If
$z'$ is a factorization of $S$ of length $\ell$, with $a'$ atoms of length
two, and if $W_1,\ldots,W_{\ell-a'}$ are the remaining atoms occurring in $z'$, then
\[
 2a+3r = |S| = 2a' + 3(\ell - a') + \delta ,
 \qquad\text{that is,}\qquad
 \ell = r + \frac{2a + a'-\delta}{3}\,,
\]
where $\delta=\sum_{i=1}^{\ell-a'}(|W_i|-3)$.
Thus $\ell \le a+r$, because $a' \le a$ and $\delta \ge 0$, with equality for the displayed factorization of $S$,
which proves the first assertion.

Suppose that $\ell = a+r-1$.  Then $a' = a-3+\delta$, so that
$\delta \in [0,3]$.  The sequence $S$ has exactly $a$ terms in $-G_0$, and the
$a'$ atoms of length two contain $a'$ of them, so that the remaining $\ell-a'$
atoms contain exactly $3-\delta$ terms in $-G_0$.  Every atom of
length three has the form $e_i^3$ and has no term in $-G_0$, and every atom of
length at least four is an atom as in Lemma~\ref{lem:cap-atoms}.2.
\begin{itemize}
\item If $\delta = 0$, then the remaining $\ell-a'$
atoms contain exactly $3$ terms in $-G_0$, while all those atoms have the form $e_i^3$, a contradiction.
\item If $\delta = 1$, then the remaining $\ell-a'$
atoms contain exactly $2$ terms in $-G_0$, while all those atoms but one have the form $e_i^3$, and the
      exceptional one has length four and hence one or three terms in $-G_0$,
      a contradiction.
\item If $\delta = 2$, then the remaining $\ell-a'$
atoms contain exactly $1$ term in $-G_0$, while the remaining atoms of length at least four are
      either one of length five, with $0$, $2$, $3$ or $5$ terms in $-G_0$, or
      two of length four, with an even number of terms in $-G_0$, a contradiction.
\item If $\delta = 3$, then no remaining atom has a term in $-G_0$, so all of
      them lie in $\mathcal A (G_0)$ and have odd length by
      Lemma~\ref{lem:cap-atoms}.1, implying that $\delta$ is even, a contradiction.
\end{itemize}
Hence $\max\mathsf L(S)-1 \notin \mathsf L(S)$.
\end{proof}

\begin{lemma}\label{le-1}
	Let $(e_1,e_2,e_3)$ be a basis of $G$, let $e_4=-e_1-e_2-e_3$, and let
	$f_i=\varepsilon_ie_i$ for all $i\in [1,4]$, where $\varepsilon_i\in \{-1,1\}$.
	If
	\[
	S=\prod_{i=1}^4\bigl((-e_i)e_i\bigr)^2f_i^{3c_i}
	\]
	is a product of atoms of length $7$, where $c_1,c_2,c_3,c_4$ are non-negative integers such that $c_1+c_2+c_3+c_4\ge 4$,
	then
	\[
	c_1=c_2=c_3=c_4=1, \ \text{ and }\ \varepsilon_1+\varepsilon_2+\varepsilon_3+\varepsilon_4=0\,.
	\]
\end{lemma}

\begin{proof} 
For each $i\in [1,4]$, let $A_i=e_i\prod_{j\in [1,4]\setminus\{i\}}(-e_j)^2$.
	By  Lemma~\ref{lem:cap-atoms}.2 (applied with the basis $(e_1,e_2,e_3)$), the
	atoms of length $7$ whose support is contained in
	$\{\pm e_1, \pm e_2, \pm e_3, \pm e_4\}$ are precisely
	$A_1,A_2,A_3,A_4,-A_1,-A_2,-A_3,-A_4$.
	Suppose
	\[
	S=\prod_{i=1}^4 A_i^{x_i}(-A_i)^{y_i},
	\qquad
	X=\sum_{i=1}^4 x_i,\quad 
	Y=\sum_{i=1}^4 y_i,\quad \text{ and }
	N=X+Y.
	\]
	Then
	\[
	\mathsf v_{e_i}(S)=x_i+2(Y-y_i),
	\qquad
	\mathsf v_{-e_i}(S)=y_i+2(X-x_i),
	\qquad \text{ and } \qquad
	7N=|S|=16+3\sum_{i=1}^4 c_i\ge 28\,,
	\]
	whence $N\equiv1\pmod3$ and $N\ge 4$.
	
	Set
	\[
	P=\{i\in [1,4]\colon \varepsilon_i=1\}\qquad \text{ and }\qquad Q=\{i\in [1,4]\colon \varepsilon_i=-1\}.
	\]
	By simultaneously negating all terms if necessary, we may assume $|P|\le 2$.
	
	If $|P|=0$, then 
	\[
	8=\sum_{i=1}^4 \mathsf v_{e_i}(S)=\sum_{i=1}^4 (x_i+2(Y-y_i))=X+6Y\,,
	\]
	whence
	\[
	N=X+Y=8-5Y\in\{8,3\},
	\]
	contrary to $N\equiv1\pmod3$.
	
	If $|P|=1$, then since any three of $\{e_1,e_2,e_3,e_4\}$ form a basis of $G$, by changing basis if necessary, we may assume that  $P=\{4\}$, whence
	\[
	2=\mathsf v_{-e_4}(S)=y_4+2(x_1+x_2+x_3)\qquad \text{ and } \qquad
	6=\sum_{i = 1}^3 \mathsf v_{e_i}(S)=x_1+x_2+x_3+4Y+2y_4\,,
	\]
	whence $x_1+x_2+x_3\in \{0,1\}$ and $4Y=2+3(x_1+x_2+x_3)\in \{2,5\}$, a contradiction.
	
	Thus $|P|=|Q|=2$, which implies that $\varepsilon_1+\varepsilon_2+\varepsilon_3+\varepsilon_4=0$. In view of the structure of $A_i$, 
	every atom $A_i$ or $-A_i$ contains either two or five terms from
	$\{-f_1,-f_2,-f_3,-f_4\}$, whence
	\[
	8=\sum_{i=1}^4 \mathsf v_{-f_i}\left(\prod_{j=1}^4\bigl((-e_j)e_j\bigr)^2\right)=\sum_{i=1}^4 \mathsf v_{-f_i}(S)=\sum_{i=1}^4 \mathsf v_{-f_i}\left(\prod_{j=1}^4 A_j^{x_j}(-A_j)^{y_j}\right)\ge  2N\,.
	\]
	In view of  $N\ge 4$,
	we obtain $N=4$ and hence $c_1+c_2+c_3+c_4=4$.
	
	Assume that there exists $i\in [1,4]$ such that $c_i=0$.
	 Since any three of $\{e_1,e_2,e_3,e_4\}$ form a basis of $G$, by changing basis if necessary, we may assume that $c_4=0$.
	  Since each of the atoms
	$\pm A_j$ with $j \ne 4$ has two terms in $\{e_4,-e_4\}$ and $\pm A_4$ has
	one, it follows from $c_4=0$ that $$2(N-x_4-y_4)+x_4+y_4= \mathsf v_{e_4}\left(\prod_{i=1}^4 A_i^{x_i}(-A_i)^{y_i}\right)+\mathsf v_{-e_4}\left(\prod_{i=1}^4 A_i^{x_i}(-A_i)^{y_i}\right) =\mathsf v_{e_4}(S)+\mathsf v_{-e_4}(S) = 4\,,$$ that is, $x_4+y_4=4=N$, so that
	$S=A_4^{x_4}(-A_4)^{y_4}$, and $x_4=\mathsf v_{e_4}(S)=2$ and
	$y_4=\mathsf v_{-e_4}(S)=2$.  But then
	$4+3c_1=\mathsf v_{e_1}(S)+\mathsf v_{-e_1}(S)=2x_4+2y_4=8$, a contradiction. Thus $c_i\ge 1$ for all $i\in [1,4]$, giving
	\[
	c_1=c_2=c_3=c_4=1. \qedhere
	\]
\end{proof}

\begin{proposition} \label{prop:saturated-three-q-plus-one}
Let $L \in \mathcal L (G)$ and let $k\in\N$.
\begin{enumerate}
\item If
      \[
       \min  L =3k+1,
       \quad\max L =7k+5,
       \quad \text{and} \quad 7k+4\notin  L ,
      \]
      then $k=1$ and
      $L  = [4,10]\cup\{12\}$.
\item If
      \[
       \min  L =3k,
       \quad\max  L =7k+2,
       \quad \text{and} \quad 7k+1\notin  L,
      \]
      then $k=2$ and 
      $L  = [6,14]\cup\{16\}$.
\end{enumerate}
\end{proposition}

\begin{proof}
Let $B \in \mathcal B (G)$ with $\mathsf L (B)=L$.

1.  Set $m = 3k+1$ and $M = 7k+5$. Let $x, z \in \mathsf Z (B)$ with $|x|=\min \mathsf L (B)$ and $|z|=\max \mathsf L (B)$.  Denote by $a$ the number of atoms of length two occurring in $z$.  Then $|B| \le 7m$ and $|B| \ge 2a+3(M-a) = 3M-a$, while
$a \le 8$ by Lemma~\ref{lem:upper-hole-affine-cap}. Since
$3M = 21k+15 = 7m+8$, all these estimates are equalities, whence 
\[
0\nmid B,\quad |B| = 7m = 21k+7,
 \qquad \text{ and }\qquad a = 8,
\]
every atom occurring in $x$ has length seven, and the remaining $M-8 = 7k-3$
atoms occurring in $z$ have length three.  By 
Lemma~\ref{lem:upper-hole-affine-cap}, there is a basis $(e_1,e_2,e_3)$ of
$G$ such that, with $e_0 = e_1+e_2+e_3$ and $e_4 = -e_0$, the product of the
atoms of length two occurring in $z$ is
\[
D = \prod_{i=1}^{4}\bigl((-e_i)e_i\bigr)^2 = \prod_{i=0}^{3}\bigl((-e_i)e_i\bigr)^2 ,
\]
the sequence of Subsection~\ref{sec:software:doubled-cap}.
Let $T$ be one of the remaining $7k-3$ atoms of length three.  Since
$M - 1 \notin \mathsf L (B)$, it follows from  Lemma~\ref{lem:hereditary-gap}.2 that 
$8 \notin \mathsf L (DT)$, so that $T = (\varepsilon e_i)^3$ for some
$i \in [1,4]$ and some $\varepsilon \in \{1,-1\}$ by
Lemma~\ref{lem:comp:doubled-cap-length-three}, since $(\varepsilon e_0)^3 = (-\varepsilon e_4)^3$.  If both $e_i^3$ and
$(-e_i)^3$ occurred among these atoms, then $\bigl((-e_i)e_i\bigr)^3 = e_i^3(-e_i)^3$ would give a factorization of $B$ of length $M+1$, a contradiction.
Thus there are
$\varepsilon_i \in \{1,-1\}$ and $c_i \in \N_0$ such that
\[
 B = D\prod_{i=1}^{4}(\varepsilon_i e_i)^{3c_i},
 \qquad \sum_{i=1}^{4}c_i = 7k-3\ge 4.
\]
Since all atoms occurring in $x$ have length seven, Lemma \ref{le-1} gives $c_1=c_2=c_3=c_4=1$ and $\varepsilon_1+\varepsilon_2+\varepsilon_3+\varepsilon_4=0$, whence
$k=1$,   two of $\varepsilon_1,\varepsilon_2,\varepsilon_3,\varepsilon_4$
are equal to $1$, and two are equal to $-1$.  For every permutation $\pi$ of
$[1,4]$, the automorphism of $G$ with $e_i \mapsto e_{\pi (i)}$ for
$i \in [1,3]$ maps $e_4 = -e_1-e_2-e_3$ to $e_{\pi (4)}$ and fixes $D$.
Hence we may assume $\varepsilon_2=\varepsilon_3=1$ and
$\varepsilon_1=\varepsilon_4=-1$, that is, $B = D(-e_1)^3e_2^3e_3^3(-e_4)^3 = De_0^3(-e_1)^3e_2^3e_3^3$, whose length set is
$[4,10]\cup\{12\}$ by Lemma~\ref{lem:comp:doubled-cap-seeds}.

2.  For $k \in \N$, we say that $B$ satisfies $(\mathrm H_k)$ if
\[
 \min\mathsf L(B)=3k,
 \qquad \max\mathsf L(B)=7k+2,
 \qquad 7k+1\notin\mathsf L(B) .
\]

We prove by induction on $k$ the following assertion: if $B$ satisfies
$(\mathrm H_k)$, then $k=2$, and there is a basis $(e_1,e_2,e_3)$ of $G$
such that $B = D\bigl(e_0^3(-e_1)^3e_2^3e_3^3\bigr)^2$, where, as in Part 1, $e_0=e_1+e_2+e_3$ and $D=\prod_{i=0}^{3}\bigl((-e_i)e_i\bigr)^2$.
The length set
$[6,14]\cup\{16\}$ of this zero-sum sequence is given by
Lemma~\ref{lem:comp:doubled-cap-seeds}.

Let $B$ satisfy $(\mathrm H_k)$, and choose $x\in \mathsf Z (B)$ with
$|x| = 3k$.  If $A^3\nmid x$ for every $A\in \mathcal A(G)$ of length $7$, then
Lemma~\ref{lem:comp:doubled-cap-cubes} gives the assertion.

Suppose now that there is an atom
$A = g_1 \cdot \ldots \cdot g_7$ of length $7$ such that $A^3 \mid x$.  Since $g_1^3, \ldots, g_7^3$ are atoms with $g_1^3 \cdot \ldots \cdot g_7^3 = A^3 \mid B$,
applying Lemma~\ref{A:short-atom} seven times (after removing  $g_i^3$
from $B$, the maximum decreases by one and its predecessor is still not a
length), we obtain some $z\in \mathsf Z (B)$ with $|z| = 7k+2$ such that $g_1^3 \cdot \ldots \cdot g_7^3 \mid z$,
so $B' = BA^{-3} \in \mathcal B (G)$ has the factorizations
\[
 x' = x A^{-3}
 \qquad \text{and} \qquad
 z' = z \bigl(g_1^3 \cdot \ldots \cdot g_7^3\bigr)^{-1}
\]
of lengths $3k-3$ and $7k-5 = 7(k-1)+2$.  Moreover, $3(k-1)=\min \mathsf L(B')$, $7(k-1)+2=\max\mathsf L(B')$, and $7(k-1)+1\not\in \mathsf L(B')$, because
adjoining the three atoms $A$ resp.\ the seven atoms $g_1^3, \ldots, g_7^3$ to a factorization of $B'$ of length $\ell$
gives a factorization of $B$ of length $\ell+3$ resp.\ $\ell+7$.
If $k = 1$, then $|x'| = 0$ forces $B' = 1$, which
contradicts $|z'| = 2$.  If $k \ge 2$, then $B'$ satisfies
$(\mathrm H_{k-1})$, and by the induction hypothesis, we obtain that $k-1 = 2$ and there is
a basis $(e_1,e_2,e_3)$ of $G$ such that
$B' = D\bigl(e_0^3(-e_1)^3e_2^3e_3^3\bigr)^2$ with $e_0$ and $D$ as above. 
 Let $Y = e_0^3(-e_1)^3e_2^3e_3^3$, so that $B = B'A^3 = DY^2A^3$, and let
$B^*=BY^{-1} = DYA^3$.  Applying Lemma~\ref{A:short-atom} to the four atoms $(-e_1)^3,e_2^3,e_3^3,e_0^3$ as above gives $\max\mathsf L(B^*) = \max\mathsf L(B) - 4 = 19$ and
$18 \notin \mathsf L(B^*)$.  Moreover $|B^*| = 16+12+21 = 49$ gives
$\min\mathsf L(B^*) \ge 7$, and $4 \in \mathsf L(DY)$ by
Lemma~\ref{lem:comp:doubled-cap-seeds} gives $7 \in \mathsf L(B^*)$.  Hence
\[
\min \mathsf L(B^*)=7=3\cdot 2+1, \quad \max\mathsf L(B^*)=19=7\cdot 2+5, \quad\text{and}\quad 18\notin\mathsf L(B^*)\,,
\]
contradicting Part 1, which forces $k=1$ for such a zero-sum sequence.
\end{proof}

\begin{proposition}\label{P-D}
Let $(e_1,e_2,e_3)$ be a basis of $G$,  let $e_0=e_1+e_2+e_3$, and let $U=e_1^2e_2^2e_3^2e_0$,  
	\begin{align*}
	&	D_6= \big(e_1(-e_1) \big) \big(e_2(-e_2) \big)^2 \big(e_3(-e_3) \big)^2 \big(e_0(-e_0) \big), \quad
	D_7= \big(e_1(-e_1) \big)^2 \big(e_2(-e_2) \big)^2 \big(e_3(-e_3) \big)^2 \big(e_0(-e_0) \big),\\
	& D_8= \big(e_1(-e_1) \big)^2 \big(e_2(-e_2) \big)^2 \big(e_3(-e_3) \big)^2 \big(e_0(-e_0) \big)^2, \quad 
	D_{10}=D_8e_1^3e_2^3, \quad
	D_{11}=D_8e_1^3e_2^3e_3^3\,.
	\end{align*}
Then for every $k\in \N_0$, we have
\begin{enumerate}
	\item $\mathsf L(D_6U^{3k})=[3k+2, 7k+4]\cup \{7k+6\}$;
	\item $\mathsf L(D_7U^{3k})=[3k+2, 7k+5]\cup \{7k+7\}$;
	\item $\mathsf L(D_{10}U^{3k})=[3k+4, 7k+8]\cup \{7k+10\}$;
	\item $\mathsf L(D_{11}U^{3k})=[3k+4, 7k+9]\cup \{7k+11\}$.
\end{enumerate}
\end{proposition}

\begin{proof}
Let $D\in \{D_6,D_7, D_{10}, D_{11}\}$. It follows from Lemma \ref{lem:oriented-cap-top-gap} that $\max\mathsf L(DU^{3k})=\max\mathsf L(D)+7k$ and $\max\mathsf L(D)+7k-1\not\in \mathsf L(DU^{3k})$.
Since $\min \mathsf L(D)+3k\ge \min \mathsf L(DU^{3k})\ge \left\lceil \frac{|D|+21k}{7}\right\rceil =3k+\left\lceil \frac{|D|}{7}\right\rceil$,
we obtain $\min \mathsf L(D)+3k=\min \mathsf L(DU^{3k})$, because
$\min \mathsf L(D) = \lceil |D|/7 \rceil$ by Lemma \ref{lem:comp:doubled-cap-seeds}.
Since $1 \in \Delta \big(\mathsf L(D) \big)$ by Lemma \ref{lem:comp:doubled-cap-seeds}
and $D \mid DU^{3k}$, we have $1 \in \Delta \big(\mathsf L(DU^{3k}) \big)$, and the
assertions follow from Theorem \ref{thm:mixed-shape}.
\end{proof}

\newpage
\section{Proof of Theorem \ref{1.1}} \label{6}
\smallskip

\begin{proof}[Proof of Theorem \ref{1.1}]
Let $G = C_3 \oplus C_3 \oplus C_3$ and let $H$ be a transfer Krull monoid over $G$. By Subsection \ref{2.4}, we have $\mathcal L (H) = \mathcal L (G)$, whence it suffices to consider the monoid $\mathcal B (G)$. 
We fix a basis $(e_1,e_2,e_3)$ of $G$ and set
\begin{align*}
&U_1=e_1^2e_2^2e_3^2(e_1+e_2+e_3), \quad  U_2=e_1^2e_2^2e_3(e_1+e_3)(e_2+e_3)\,,\\
&V_1=e_2^2(e_1+e_2+e_3)^2e_1(-e_2+e_3),\quad  V_2=(-e_1)^2e_3^2(-e_1-e_2-e_3)(e_2-e_3),\\
& V_3=(-e_2)^2(-e_3)^2e_1(-e_1-e_2-e_3),\quad
V_4=e_1^2e_2^2(e_1+e_2+e_3)^2(-e_1-e_2+e_3),\\
&V_5=(-e_1)^2e_3^2(-e_1-e_2-e_3)^2(e_1-e_2), \quad
 V_6=(-e_3)^2(-e_2)^2(-e_1-e_2-e_3)(-e_1+e_2)(-e_1-e_2),\\
&V_7=e_2^2e_3^2(e_1+e_2+e_3)^2(e_1-e_2-e_3),\quad
V_8=e_1^2(-e_3)^2(-e_1-e_2-e_3)^2(-e_2+e_3), \ \text{and} \\
& V_9=(-e_1)^2(-e_2)^2(-e_3-e_2)(-e_1-e_2-e_3)(-e_3+e_2)\,.
\end{align*}
Then $U_1, U_2, V_4,V_5,V_6,V_7,V_8,V_9$ are  atoms of length $7$ and $V_1,V_2,V_3$ are atoms of length $6$ with \begin{align*}
V_1V_2V_3
= & \, (-U_1) U_1 \Big((e_1+e_2+e_3)(-e_1-e_2-e_3)\Big) \Big((-e_2+e_3)(e_2-e_3) \Big)\\
= & \, \Big(e_2(-e_3)(-e_2+e_3)\Big)  \Big( (-e_2)e_3(e_2-e_3)\Big) \cdot \\
 & \quad \Big(e_1(-e_1)\Big)^2 \Big( e_2(-e_2)\Big) \Big( e_3(-e_3)\Big)   \Big( (e_1+e_2+e_3)(-e_1-e_2-e_3) \Big)^2\,,
\end{align*}
whence $\big( 2+\mathsf L \big((-U_1)U_1 ) \big) \cup\{8\} \subset \mathsf L(V_1V_2V_3)$.
We claim that
\[
\mathsf L \big((-U_1)U_1 \big)=[2,5]\cup\{7\}, \ \mathsf L \big((-U_2)U_2 \big)=[2,7], \mathsf L(V_4V_5V_6)=\mathsf L(V_7V_8V_9)=[3,10], \ \text{and} \ \mathsf L(V_1V_2V_3)=[3,9] \,.
\]
It is easy to see that $\mathsf L \big((-U_1)U_1 \big)=[2,5]\cup\{7\}$, whence  $\mathsf L(V_1V_2V_3)=[3,9]$. Since $V_7V_8V_9=V_1V_2V_3\big((-e_1-e_2-e_3)(e_1-e_2-e_3)(-e_2-e_3)\big)$, we obtain that $\mathsf L(V_7V_8V_9)=[3,10]$. It follows by changing basis that $\mathsf L(V_4V_5V_6)=\mathsf L(V_7V_8V_9)=[3, 10]$.

Let $W$ be an atom dividing $(-U_2)U_2$. Then $2+|U_2|-|W|\in \mathsf L \big((-U_2) U_2\big)$. Since $W_1=e_1e_3(-e_1-e_3)$,
$W_2=e_1(-e_2)(-e_1-e_3)(e_2+e_3)$, $W_3=(-e_1)(-e_2)e_3(e_1+e_3)(e_2+e_3)$, $W_4=(-e_1)e_2^2e_3(e_1+e_3)(e_2+e_3)$ are all atoms, we have
\[
\{3,4,5,6\}\subset \mathsf L \big((-U_2) U_2 \big)\,,
\]
whence $[2,7]=\mathsf L \big((-U_2)U_2 \big)$.
Note that $\{L\in \mathcal L(G)\colon \Delta (L)=\emptyset\}=\mathcal L_1$. We proceed to finish the proof by a series of claims.

\begin{claim}\label{k5-A1}
	$[2k, \nu] \in \mathcal L (G)$ for all $\nu \in [2k, 7k]$.
\end{claim}

\begin{proof}[Proof of \ref{k5-A1}] We proceed by induction on $k$.
	Note that $\mathsf L \big((-U_2)U_2 \big)=[2,7]$, \ $\mathsf L \big( (-W_4)W_4 \big)=[2,6]$,
	and $[2,5], [2,4], [2,3]\in \mathcal L(C_3^2)\subset \mathcal L(G)$.
	Thus the assertion holds for $k=1$.
	Suppose the assertion holds for $k \ge 1$. If $\nu \in [2k, 7k]$, then there exists $A_{\nu} \in \mathcal B (G)$ with $\mathsf L (A_{\nu}) = [2k, \nu]$, and hence $\mathsf L (0^2 A_{\nu}) = [2k+2, \nu+2]$. Thus it remains to show that
	\[
	[2k+2, 7k+3], \ [2k+2, 7k+4], \ [2k+2, 7k+5], \ [2k+2, 7k+6], \ \text{and} \ [2k+2, 7k+7] \in \mathcal L (G)\,.
	\]
	If $U_2, W_1, W_2, W_3$, and $W_4$ are as above, then
	\begin{align*}
	\mathsf L \big( (-U_2)^{k-1}U_2^{k-1}V_1V_2V_30 \big) & = [2k+2, 7k+3]\,, \\
	\mathsf L \big( (-U_2)^k U_2^k (-W_2) W_2 \big) & = [2k+2, 7k+4] \,,\\
	\mathsf L \big( (-U_2)^k U_2^k (-W_3)W_3 \big) & = [2k+2, 7k+5] \,, \\
	\mathsf L \big( (-U_2)^k U_2^k (-W_4) W_4 \big) & = [2k+2, 7k+6] \,, \quad \text{and} \\
	\mathsf L \big( (-U_2)^{k+1} U_2^{k+1} \big) & = [2k+2, 7k+7] ,	\end{align*}
	as desired.  Indeed, in each of these five cases the displayed interval is
	contained in the length set, because the length set of a product contains
	the sum of the length sets of its displayed divisors and
	$\mathsf L \big((-U_2)U_2\big) = [2,7]$, $\mathsf L (V_1V_2V_3) = [3,9]$,
	$\mathsf L (0) = \{1\}$ and $\mathsf L \big((-W_j)W_j\big) \supset
	\{2, |W_j|\}$; and the displayed interval contains the length set, because every atom except  $0$ has length at
	least two  and at most $\mathsf D (G) = 7$, so that
	$\lceil |B|/7 \rceil \le \min \mathsf L (B)$ and
	$\max \mathsf L (B) \le \lfloor |B|/2 \rfloor$ for every
	$B \in \mathcal B (G^{\bullet})$, and the two bounds coincide with the
	endpoints of the displayed interval. Note that for the first case, we show $\mathsf L \big( (-U_2)^{k-1}U_2^{k-1}V_1V_2V_3 \big)  = [2k+1, 7k+2]$ and then shift the length set by one. \qedhere (\ref{k5-A1})
\end{proof}

\begin{claim}\label{k5-A2} $[2k+1, \nu] \in \mathcal L (G)$ for all $\nu \in [2k+1, 7k+3]$.
\end{claim}

\begin{proof}[Proof of \ref{k5-A2}] If  $k \in \N$, $\nu \in [2k, 7k]$, and $A_{\nu} \in \mathcal B (G)$ with $\mathsf L (A_{\nu}) = [2k, \nu]$, then $\mathsf L (0 A_{\nu}) = [2k+1, \nu+1]$. Since $\mathsf L \big((-U_2)^{k-1}U_2^{k-1}V_4V_5V_6 \big) = [2k+1, 7k+3] \in \mathcal L (G)$ and $\mathsf L \big((-U_2)^{k-1}U_2^{k-1}V_1V_2V_3 \big)=[2k+1, 7k+2]$, the assertion follows.\qedhere (\ref{k5-A2})
\end{proof}

\begin{claim}\label{k5-S1}
	$\big\{L\in \mathcal L(G)\colon \Delta (L)=\{1\} \big\}=\mathcal L_2$.
\end{claim}

\begin{proof}[Proof of \ref{k5-S1}]
	For every $L \in \mathcal L (G)$ with $\Delta (L)=\{1\}$, we have   $\rho(L)\le 7/2$ and hence $L$ equals one of the sets given in $\mathcal L_2$. It remains to show that every set $L$ given in $\mathcal L_2$ can be realized as a length set in $\mathcal L (G)$. Letting $k \in \N$, this follows from  \ref{k5-A1} and \ref{k5-A2}.\qedhere (\ref{k5-S1})
\end{proof}

\begin{claim}\label{k5-S2}  $\big\{L\in \mathcal L(G)\colon \Delta (L)=\{2\} \big\}=\mathcal L_3$.
\end{claim}

\begin{proof}[Proof of \ref{k5-S2}]
	Let $G_1=\{0,e_1,e_2,e_3,e_1+e_2+e_3\}$. Then $\Delta(G_1)=\{2\}$ and $U_1\in \mathcal A(G_1)$. Indeed, every atom of $\mathcal B (G_1)$ has odd length (the atom $0$ has length one, and the others by  Lemma~\ref{lem:cap-atoms}.1), whence $\ell \equiv |B| \mod 2$ for every $B \in \mathcal B (G_1)$ and every $\ell \in \mathsf L (B)$ by  Lemma~\ref{lem:parity-formula}.1.  Let $U_3=e_1e_2e_3(e_1+e_2+e_3)^2$. Then $\mathsf L(U_1U_3)=\{2,4\}$ and $\mathsf L(U_1^3)=\{3,5,7\}$.
	
	Let $y\in \N_0$ and $k\in \N$. If $k=2t$ is even, then
	\[
	\mathsf L \big(0^y  U_1^{3t} \big)=y+t\{3,5,7\}=y+3t+2\cdot [0,2t]=y+\left\lceil\frac{3k}{2}\right\rceil+2\cdot [0,k]\,.
	\]
	If $k=2t-1$ is odd, then
	\[
	\mathsf L \big(0^y  (U_1^3)^{t-1}U_1U_3 \big)=y+(t-1)\mathsf L(U_1^3)+\mathsf L(U_1U_3)=y+3t-1+2\cdot [0, 2t-1]=y+\left\lceil\frac{3k}{2}\right\rceil+2\cdot [0,k]\,.
	\]
	Every set in $\mathcal L_3$ has this form: if $L = m + 2 \cdot [0,\nu]$ with
	$m = 3k$ and $\nu \le 2k$, $m = 3k+1$ and $\nu \le 2k$, or $m = 3k+2$ and
	$\nu \le 2k+1$, then $\lceil 3\nu/2 \rceil \le m$, and
	$L = y + \lceil 3\nu/2 \rceil + 2 \cdot [0,\nu]$ with $y = m - \lceil 3\nu/2 \rceil \in \N_0$.
	Therefore, we have shown that $\mathcal L_3\subset \{L\in \mathcal L(G)\colon \Delta (L)=\{2\}\}$.

	Let $L = \mathsf L (B)$ with $B \in \mathcal B(G)$ and $\Delta(L)=\{2\}$. To complete \ref{k5-S2}, it suffices to show $L\in \mathcal L_3$.
	By Theorem~\ref{thm:lower-gap-endpoint}, we have
	$\rho(L)\le 7/3$. If $\min L=3k$ for some $k\in \N$, then $\max L\le 7k$ and hence $L\in \mathcal L_3$. If $\min L=3k+2$ for some $k\in \N_0$, then $\max L\le 7k+4$ and hence $L\in \mathcal L_3$. Suppose $\min L=3k+1$ for some $k\in \N$. Then $\max L\le 7k+2$. If $\max L=7k+2$, then $\max L-\min L=4k+1$ is odd, a contradiction to $\Delta(L)=\{2\}$. Thus $\max L\le 7k+1$ and hence $L\in \mathcal L_3$.\qedhere (\ref{k5-S2})
\end{proof}

\begin{claim}\label{k5-A3}
$\mathcal L_4\subset\big\{ L\in \mathcal L(G)\colon \Delta (L)=[1,2] \big\}$.
\end{claim}

\begin{proof}[Proof of \ref{k5-A3}]
Since $[4,10]\cup\{12\}, [6, 14]\cup\{16\}\in \mathcal L(G)$ by Lemma \ref{lem:comp:doubled-cap-seeds}, 	
it suffices to show that
\begin{align*}
	\mathcal L(G)\ \supset \
	&\big\{[3k+2, \nu]\cup\{\nu+2\} \colon k \in \N_0, \nu\in [3k+3, 7k+5] \big\} \ \sqcup \\
	\quad\quad  &\big\{[3k+3, \nu]\cup\{\nu+2\} \colon k \in \N_0, \nu\in [3k+4, 7k+6] \big\} \ \sqcup\\
	\quad\quad &\big\{[3k+4, \nu]\cup\{\nu+2\} \colon k \in \N_0, \nu\in [3k+5, 7k+9] \big\}\,.
\end{align*}	

We proceed by induction on $k$. 
If $k=0$, then Lemma
\ref{lem:comp:doubled-cap-seeds} implies that
$[2,3]\cup \{5\}, [2,4]\cup \{6\}, [2,5]\cup \{7\}\in \mathcal L(G)$,  and $[4,8]\cup \{10\}, [4,9]\cup \{11\}\in \mathcal L(G)$. Since the shifts of length sets are  also length sets, we obtain  
\[
\big\{[2, \nu]\cup\{\nu+2\} \colon \nu\in [3, 5] \big\} \cup \big\{[3, \nu]\cup\{\nu+2\} \colon \nu\in [4, 6] \big\} \cup\big\{[4, \nu]\cup\{\nu+2\} \colon \nu\in [5, 9] \big\}\subset \mathcal L(G)\,.
\]

Suppose the claim  holds for $k=t$, where $t\in \N_0$. Then $[3t+4, \nu]\cup\{\nu+2\}\in \mathcal L(G)$ for all $\nu\in [3t+5, 7t+9]$, which implies that
$[3t+5, \nu]\cup\{\nu+2\}\in \mathcal L(G) $ for all $\nu\in [3t+6, 7t+10]$. By Proposition \ref{P-D}, we have that
$[3t+5, \nu]\cup\{\nu+2\}\in \mathcal L(G) $ for all $\nu\in [7t+11, 7t+12]$. Thus
\begin{center}
$[3(t+1)+2, \nu]\cup\{\nu+2\}\in \mathcal L(G) $ for all $\nu\in [3(t+1)+3, 7(t+1)+5]$
\end{center}
 and hence 
\begin{center}
$[3(t+1)+3, \nu]\cup\{\nu+2\}\in \mathcal L(G) $ for all $\nu\in [3(t+1)+4, 7(t+1)+6]$.
\end{center}
Moreover, we have $[3t+7, \nu]\cup\{\nu+2\}\in \mathcal L(G) $ for all $\nu\in [3t+8, 7t+14]$, and Proposition \ref{P-D} (applied to $D_{10}$ and $D_{11}$ with $k=t+1$) gives the values $\nu = 7t+15$ and $\nu = 7t+16$, so that
\begin{center}
\hspace{2cm}$[3(t+1)+4, \nu]\cup\{\nu+2\}\in \mathcal L(G) $ for all $\nu\in [3(t+1)+5, 7(t+1)+9]$. \qedhere (\ref{k5-A3})
\end{center}
\end{proof}

\begin{claim}\label{k5-S3} $\mathcal L_5\subset \big\{ L\in \mathcal L(G)\colon \Delta (L)=[1,2] \big\}$.
\end{claim}

\begin{proof}[Proof of \ref{k5-S3}]
We proceed by induction on $k$ to show that
	\begin{align*}
	\mathcal L(G)\supset \mathcal L_5
		&=\big\{\{3k+2\}\cup[3k+4,\nu] \colon k \in \N_0, \nu\in [3k+5, 7k+5] \big\} \ \sqcup \\
		&\quad\quad  \big\{\{3k+3\}\cup[3k+5,\nu] \colon k \in \N_0, \nu\in [3k+6, 7k+7] \big\} \ \sqcup\\
		&\quad\quad \big\{\{3k+4\}\cup[3k+6,\nu] \colon k \in \N_0, \nu\in [3k+7, 7k+9] \big\}\,.
		\end{align*}
	 If $k=0$, then Lemma
	\ref{le3.5} (applied with $k=0$) implies
	$\{2\}\cup [4,5]\in \mathcal L(G)$, $\{3\}\cup [5,7]\in \mathcal L(G)$, and $\{4\}\cup [6,9]\in \mathcal L(G)$. Since the shifts of length sets are  also length sets, we obtain  $\{3\}\cup [5,6]$, $\{4\}\cup [6,7]$, and $\{4\}\cup [6,8]$ as well, whence the assertion holds for $k=0$.
	
	Suppose the assertion  holds for $k=t$, where $t\in \N_0$. Then $\{3t+4\}\cup [3t+6, \nu]\in \mathcal L(G) $ for all $\nu\in [3t+7, 7t+9]$, which implies that
	$\{3t+5\}\cup [3t+7, \nu]\in \mathcal L(G) $ for all $\nu\in [3t+8, 7t+10]$. Again by Lemma \ref{le3.5} (applied with $k=t+1$), we have that
	$\{3t+5\}\cup [3t+7, \nu]\in \mathcal L(G) $ for
	$\nu\in [7t+11, 7t+12]$. Thus
	\begin{center}
	$\{3(t+1)+2\}\cup [3(t+1)+4, \nu]\in \mathcal L(G) $ for all $\nu\in [3(t+1)+5, 7(t+1)+5]$
\end{center} and hence $\{3t+6\}\cup [3t+8, \nu]\in \mathcal L(G) $ for all $\nu\in [3t+9, 7t+13]$ by shifting. It follows from Lemma
	\ref{le3.5} that
	\begin{center}
	$\{3(t+1)+3\}\cup [3(t+1)+5, \nu]\in \mathcal L(G) $ for all $\nu\in [3(t+1)+6, 7(t+1)+7]$.
	\end{center} Since $\{3t+7\}\cup [3t+9, \nu]\in \mathcal L(G) $ for all $\nu\in [3t+10, 7t+15]$ by shifting, it follows from Lemma
	\ref{le3.5} again that 
	\begin{center}
	\hspace{2cm}	$\{3(t+1)+4\}\cup [3(t+1)+6, \nu]\in \mathcal L(G) $ for all $\nu\in [3(t+1)+7, 7(t+1)+9]$.	
		\qedhere (\ref{k5-S3})
	\end{center}
\end{proof}

\begin{claim}\label{k5-S4}  $\mathcal L_4\cup \mathcal L_5\supset \big\{ L\in \mathcal L(G)\colon \Delta (L)=[1,2] \big\}$.
\end{claim}

\begin{proof}[Proof of \ref{k5-S4}]
Let $L\in \mathcal L(G)$ with $\Delta(L)=[1,2]$. Then Theorem
\ref{thm:mixed-shape} implies that  \begin{center}
	either $L=[\min L, \max L-2]\cup \{\max L\}$ or $L=\{\min L\}\cup [\min L+2, \max L]$.
\end{center}

Suppose $L=[\min L, \max L-2]\cup \{\max L\}$.  Proposition
\ref{prop:upper-hole-endpoint} gives $\max L - 2 \le 7k, 7k+3, 7k+5$
according as $\min L = 3k, 3k+1, 3k+2$, and in the first two cases
Proposition~\ref{prop:saturated-three-q-plus-one} excludes the largest value
unless $L = [4,10] \cup \{12\}$ (for $\min L = 4$) or
$L = [6,14] \cup \{16\}$ (for $\min L = 6$). Assume that $L\not\in \left\{[4,10] \cup \{12\}, [6,14] \cup \{16\}\right\}$. Then
  $\min L = 3k+2$,
$3k+3$ or $3k+4$ with $k \in \N_0$ implies $\max L - 2 \le 7k+5$, $7k+6$, or $7k+9$, respectively.  It follows that
$L\in\mathcal L_4$.

Suppose $L=\{\min L\}\cup [\min L+2, \max L]$. It follows from
Theorem \ref{thm:lower-gap-endpoint} that $L\in \mathcal L_5$. \qedhere (\ref{k5-S4})
\end{proof}
This completes the proof of Theorem \ref{1.1}.
\end{proof}

\bibliographystyle{amsplainurl}
\bibliography{length-sets}

\end{document}